\documentclass[12pt, reqno]{amsart}
\usepackage{amsthm}
\usepackage{amssymb}
\usepackage{amsmath}

\def\1{\hbox{1\kern-.35em\hbox{1}}}

\newtheorem{theorem}{Theorem}[section]
\newtheorem*{theorem*}{Theorem}
\newtheorem{lemma}[theorem]{Lemma}
\newtheorem{proposition}[theorem]{Proposition}
\newtheorem*{proposition*}{Proposition}
\newtheorem{corollary}[theorem]{Corollary}
\newtheorem{definition}[theorem]{Definition}
\newtheorem{remark}[theorem]{Remark}

\newtheorem{example}[theorem]{Example}

\numberwithin{equation}{section}

\newcommand{\bea}{\begin{eqnarray}}
\newcommand{\eea}{\end{eqnarray}}
\newcommand{\be}{\begin{eqnarray*}}
\newcommand{\ee}{\end{eqnarray*}}

\newcommand{\Z}{{\mathbb Z}}

\newcommand{\C}{{\mathbb C}}

\newcommand{\fb}{{\mathfrak b}}
\newcommand{\fg}{{\mathfrak g}}
\newcommand{\fh}{{\mathfrak h}}

\newcommand{\ad}{{\rm ad}}

\newcommand{\Hom}{{\rm Hom}}

\newcommand{\U}{{\rm U}}

\def\ker#1{{\rm Ker}(#1)}
\def\Im#1{{\rm Im}(#1)}

\def\ptl{\partial}

\def\LL{{\cal L}}

\def\a{\alpha}
\def\KK{{\mathcal K}}
\def\LL{{\mathcal L}}

\def\b{\beta}
\def\d{\delta}
\def\D{\Delta}
\def\g{\gamma}
\def\G{\Gamma}
\def\l{\lambda}
\def\L{\Lambda}

\def\ad#1{{\rm ad\sc\,}(#1)}

\def\Si{\Sigma}
\def\si{\sigma}
\def\sc{\scriptstyle}
\def\ssc{\scriptscriptstyle}
\def\dis{\displaystyle}
\def\cl{\centerline}

\def\ol{\overline}

\def\wt{\widetilde}
\def\wh{\widehat}
\def\rar{\rightarrow}
\def\rrar{\rightarrowtail}
\def\llar{\leftarrowtail}

\def\erar{\rb{3pt}{\mbox{$\ {}^{\sc\ \,e}_{\dis\rrar}\ $}}}
\def\frar{\rb{3pt}{\mbox{$\ {}^{\sc\ \,f}_{\dis\rrar}\ $}}}
\def\elar{\rb{3pt}{\mbox{$\ {}^{\sc\ \,e}_{\dis\llar}\ $}}}
\def\flar{\rb{3pt}{\mbox{$\ {}^{\sc\ \,f}_{\dis\llar}\ $}}}
\def\drar{\mbox{${\sc\,}\rrar\!\!\!\!\rar{\sc\,}$}}

\def\dlar{\mbox{${\sc\,}\lar\!\!\!\!\llar{\sc\,}$}}

\def\link{\mbox{${\ssc\,}\lar\!\!\!\!-\!\!\!\!\rar{\ssc\,}$}}

\def\LRar{\Longrightarrow}
\def\lar{\leftarrow}
\def\D{\Delta}

\def\Lra{\Longleftrightarrow}
\def\bs{\backslash}

\def\rb{\raisebox}
\def\vs{\vspace*}

\def\Z{\mathbb{Z}}

\def\C{\mathbb{C}}

\def\gl{{\mathfrak g}}

\def\glss{{{\mathfrak g}_{0}^{\rm ss}}}
\def\vp{\varphi}
\def\es{\epsilon}
\def\deg#1{[#1]}
\def\level{\ell}
\def\top{{\it top\sc\,}}
\def\bottom{{\it bottom\sc\,}}
\def\amax{{\alpha_{\rm max}}}
\def\amin{{\alpha_{\rm min}}}
\def\refequa#1{(\ref{#1})}
\def\stl{\stackrel}
\def\ob{\overbrace}
\def\ETa{{2\amax}}
\def\ETb{{-2\amin}}
\def\ETc{{\eta^{(1)}}}
\def\ETd{{\eta^{(2)}}}

\def\cp#1{{\dot #1}}
\def\equa#1#2{
\begin{equation}\label{#1}#2\end{equation}}
\def\equan#1#2{$$#2$$}
\numberwithin{equation}{section}

\begin{document}

\title[Cohomology of Lie superalgebras]
{Cohomology of Lie superalgebras ${\mathfrak{sl}}_{m|n}$
and ${\mathfrak{osp}}_{2|2n}$ }

\author[Yucai Su]{Yucai Su}
\address{Yucai Su, Current Address:  School of Mathematics and
Statistics, University of Sydney, NSW 2006, Australia; \newline
\indent Permanent Address: Department of Mathematics,
Shanghai Jiaotong University, Shanghai 200030, China.}
\email{yucai@maths.usyd.edu.au}
\author[R. B. ZHANG]{R. B. ZHANG}
\address{R. B. ZHANG, School of Mathematics and Statistics,
University of Sydney, NSW 2006, Australia.}
\email{rzhang@maths.usyd.edu.au}

\begin{abstract} We explicitly compute the first and second cohomology
groups of the classical Lie superalgebras ${\mathfrak{sl}}_{m|n}$
and ${\mathfrak{osp}}_{2|2n}$ with coefficients in the finite dimensional
irreducible modules and the Kac modules.
We also show that the second cohomology groups of these Lie superalgebras
with coefficients in the respective universal enveloping algebras
(under the adjoint action) vanish. The latter result in particular
implies that the universal enveloping algebras ${\rm U}({\mathfrak{sl}}_{m|n})$
and ${\rm U}({\mathfrak{osp}}_{2|2n})$ do not admit any non-trivial
formal deformations of Gerstenhaber type.

\thanks{{\em 2000 Mathematics Subject Classification.
Primary  17B37, 20G42, 17B10.}}
\end{abstract}
\maketitle

\tableofcontents

%
\section{Introduction}\label{Introduction}
%

We investigate the Lie superalgebra cohomology of the type I basic classical
Lie superalgebras \cite{K77}, namely, the special linear superalgebra
${\mathfrak{sl}}_{m|n}$, and the orthosymmplectic superalgebra
${\mathfrak{osp}}_{2|2n}$. Lie superalgebra cohomology
was extensively studied by Fuks, Leites \cite{FL}, and others
(see \cite{Fu} for a review). For any basic classical simple Lie superalgebra
${\mathfrak g}$, the cohomolgy groups $H^i({\mathfrak g}, V)$
were computed \cite{FL, Fu} for all $i$ when the coefficient module $V$ is $\C$
(even though relatively little seems to be known about these cohomolgy groups when the
coefficient module is non-trivial). Variations of these cohomolgy
groups, especially the relative cohomology groups and the cohomology
groups of odd nilpotent subalgebras, have also been
studied in depth because of their importance in the contexts
of the Bott-Borel-Weil theory \cite{Pe} and the Kazhdan-Lusztig theory for Lie
supergroups \cite{Se96, Br}.

A motivation of this investigation and earlier work of one of
us with Scheunert \cite{SZ98, SZ99}
comes from the theory of quantum supergroups \cite{Ma89,Zh93, Zh98}, the foundation of
which lies in the deformation theory \cite{Ge} of universal enveloping algebras
of Lie superalgebras. Recall that
the formal deformations of an associative algebra is classified by
the second Hochschild cohomology group with coefficients in the
algebra itself (regarded as a bi-module) \cite{Ge}. In the case of the
universal enveloping algebra of a Lie superalgebra, this
Hochschild cohomology group can be shown to be isomorphic to the
second Lie superalgebra cohomology group with coefficient module
being the universal enveloping algebra under the adjoint action of
the Lie superalgebra. Similarly, the first Lie superalgebra
cohomology group with coefficients in the universal enveloping
algebra controls the deformations of the co-algebra structure of
the universal enveloping algebra.

One result of the present paper is Theorem
\ref{theo7.1} and (\ref{theo8.5}) of Theorem \ref{last-theo},
which states that for $\gl$ being ${\mathfrak{sl}}_{m|n}$
or ${\mathfrak{osp}}_{2|2n}$,
$H^1(\gl, \U(\gl))\ne 0$, but $H^2(\gl, \U(\gl))=0$. The
vanishing of the second cohomology group implies
that ${\rm U}({\mathfrak{sl}}_{m|n})$ and ${\rm U}({\mathfrak{osp}}_{2|2n})$
are rigid in the sense of \cite{Ge}. Therefore,
the Drinfeld versions of the quantized universal enveloping
algebras of ${\mathfrak{sl}}_{m|n}$ and ${\mathfrak{osp}}_{2|2n}$
defined with any choice of Borel subalgebras are isomorphic to the
universal enveloping algebras themselves over the power series ring.
(This was proved for the special case of ${\mathfrak{sl}}_{m|1}$ in
\cite{SZ99}.) Also $H^1(\gl, \U(\gl))\ne 0$ implies that the
co-algebra structure of $\U(\gl)$ admits non-trivial deformations,
a fact which is known from specific examples.

Another main result of this paper is the computation of the first
and second Lie superalgebra cohomology groups of
${\mathfrak{sl}}_{m|n}$ and ${\mathfrak{osp}}_{2|2n}$
with coefficients in the finite
dimensional Kac modules and the finite dimensional irreducible
modules, which is summarized in Theorems \ref{theo3.1},
\ref{theo4.1}, \ref{theo5.1}, \ref{theo6}, and \ref{last-theo}.
This result is of intrinsic interest to the understanding
of extensions of modules of these Lie superalgebras, and also
extensions of the Lie superalgebras themselves.
As a matter of fact, in his foundational paper \cite{K77} on
the theory of Lie superalgebras,
Kac posed the problem of determining the first cohomology groups
of basic classical simple Lie superalgebras with coefficients in
the finite dimensional irreducible modules (see also \cite{Kac2, Kac3}).
Part of the paper solves
the problem for the Lie superalgebras ${\mathfrak{sl}}_{m|n}$ and
${\mathfrak{osp}}_{2|2n}$. As we have alluded to earlier, when the
coefficient module $V$ is not $\C$, little seems to be known about
the cohomolgy groups $H^i({\mathfrak g}, V)$ for the basic classical
Lie superalgebras;  the  main results (Theorems
\ref{theo3.1}, \ref{theo4.1}, \ref{theo5.1}, \ref{theo6},
\ref{theo7.1}, and \ref{last-theo}) of the present paper appear to be new.

The computations of the cohomology groups are carried out in this
paper at an elementary level by exploring long exact
sequences of cohomology groups arising from short exact sequences
of modules, and also by using some elements of the Hochschild spectral
sequence associated with the maximal even subalgebras of the Lie
superalgebras. The computations also rely heavily on detailed knowledge
on structures of Kac modules. We may mention that the
analysis of structures of Kac modules is a technical and difficult
problem. This renders Subsection \ref{Technical-lemmas} rather technical.

The organization of the paper is as follows. Sections
\ref{Preliminaries} to \ref{enveloping} treat the Lie superalgebra
cohomology groups of the special linear superalgebra in detail.
Section \ref{Preliminaries} provides some necessary background material on
${\mathfrak{sl}}_{m|n}$. Sections \ref{1-Kac modules} and \ref{2-Kac module}
respectively present the computations of the first and second cohomology
groups of ${\mathfrak{sl}}_{m|n}$ with coefficients in finite dimensional
Kac modules. Section \ref{1-irr-modules} is devoted to the computation of the first
cohomology groups of ${\mathfrak{sl}}_{m|n}$ with coefficients in finite dimensional
irreducible modules, where we make use of the concepts of atypicality matrices,
northeast chains (NE) of a weight, and $n$-,
$q$-, $c$-relationships of atypical roots, which are all explained in
Appendix \ref{Appendix}. In Section \ref{2-irr-modules} we calculate
the second cohomology groups of ${\mathfrak{sl}}_{m|n}$ with coefficients in finite dimensional
irreducible modules. This section is divided into three subsections.
Subsection \ref{graph} introduces the notion of primitive
weight graphs, which is very useful for studying
the structure of indecomposable ${\mathfrak{sl}}_{m|n}$-modules such as
Kac modules. In Subsection \ref{Technical-lemmas} we analyse
structures of some Kac modules of ${\mathfrak{sl}}_{m|n}$ and establish
a series of technical lemmas needed for proving the main result on the
second cohomology groups of ${\mathfrak{sl}}_{m|n}$ with
coefficients in finite dimensional
irreducible modules. Subsection \ref{2-irr-main} proves the main result
of Section \ref{2-irr-modules}.  Section \ref{enveloping} treats the
 second cohomology groups of ${\mathfrak{sl}}_{m|n}$ with
coefficients in the universal enveloping algebra. Finally, in
Section \ref{last-section} we present the results
on the cohomology groups of ${\mathfrak{osp}}_{2|2n}$,
while omitting most of the technical details.
%
%
%
\section{Preliminaries on the special linear superalgebra}
\label{Preliminaries}
%
%
\subsection{The special linear superalgebra}
We present some background material on the special linear superalgebra
here and refer to \cite{K77, Sc, Ma97} for more details. For general
notions of graded vector spaces and graded algebraic structures we refer to
the classic paper \cite{MM} by Milnor and Moore.

We shall work over the
complex number field $\C$ throughout the paper. Given a
$\Z_2$-graded vector space $W=W_{\bar 0}\oplus W_{\bar 1}$, we
call $W_{\bar 0}$ and $W_{\bar 1}$ the even and odd subspaces,
respectively. The elements of $W_{\bar 0}\cup W_{\bar 1}$ will be
called homogeneous. Define a map $[\ ]: W_{\bar 0}\cup W_{\bar 1}
\rightarrow \Z_2$ by $[w] = {\bar i}$ if $w \in W_{\bar i}$. For
any two $\Z_2$-graded vector spaces $V$ and $W$, the space of
morphisms $\Hom_{\C}(V, W)$ is also $\Z_2$-graded with $\Hom_{\C}(V, W)_{\bar
k} =\sum_{\bar{i}+\bar{j}\equiv \bar{k} (\rm{mod} 2)}
\Hom_{\C}(V_{\bar i}, W_{\bar j})$. We write ${\rm End}_{\C}(V)$ for
$\Hom_{\C}(V, V)$.

Let $\C^{m|n}$ be the $\Z_2$-graded vector space with even
subspace $\C^m$ and odd subspace $\C^n$. Then ${\rm
End}_{\C}(\C^{m|n})$ with the $\Z_2$-graded commutator forms the
general linear superalgebra. To describe its structure, we choose
a homogeneous basis $\{ v_a \,|\, a\in {\bf I}\}$, for $\C^{m|n}$,
where ${\bf I}=\{1, 2, \ldots , m+n\}$, and $v_a$ is even if $a\le
m$, and odd otherwise. The general linear superalgebra  relative
to this basis of $\C^{m|n}$ will be denoted by
${\mathfrak{gl}}_{m|n}$. Let $E_{a b}$ be the matrix unit, namely,
the $(m+n)\times(m+n)$-matrix with all entries being zero except
that at the $(a, b)$ position which is $1$. Then $\{E_{a b}\,
|\,a,b\in{\bf I}\}$ forms a homogeneous basis of
${\mathfrak{gl}}_{m|n}$, with $E_{a b}$ being even if $a, b\le m$,
or $a, b> m$, and odd otherwise. For convenience, we let ${\bf I}_1=\{1,...,m\}$ and
${\bf I}_2=\{\cp 1,...,\cp n\}$, where we  have written $\cp \nu=\nu+m$.
Then ${\bf I}={\bf I}_1\cup{\bf I}_2$.
Define the map $[\ ]: {\bf I}\rightarrow \Z_2$ by
$[a]=\big\{\begin{array}{l
l}
               \bar{0},  & \mbox{if} \ a\in{\bf I}_1, \\
               \bar{1},  & \mbox{if} \ a\in{\bf I}_2.
              \end{array} $
Now the commutation relations of the general linear superalgebra
${\mathfrak{gl}}_{m|n}$ can be written as
\be [E_{a b}, \ E_{c d}] &=& E_{a d}\delta_{b c}
 - (-1)^{([a]-[b])([c]-[d])} E_{c b}\delta_{a d}.
\ee

The upper triangular matrices form a Borel subalgebra $B$ of
${\mathfrak{gl}}_{m|n}$,
which contains the Cartan subalgebra $H$ of diagonal matrices.
Let $\{\epsilon_a \,|\, a\in{\bf I}\}$ be the basis
of $H^*$ such that $\epsilon_a(E_{b b})=\delta_{a b}$.
The supertrace induces a bilinear form $(\;
, \: ): H^*\times H^* \rightarrow \C$ on $H^*$ such
that $(\epsilon_a, \epsilon_b)=(-1)^{[a]} \delta_{a b}$.
Relative to the Borel subalgebra $B$,
the roots of ${\mathfrak{gl}}_{m|n}$
can be expressed as $\epsilon_a-\epsilon_b, \:\, a\ne b$, where
$\epsilon_a-\epsilon_b$ is even if $[a]+[b] = \bar{0}$ and odd
otherwise. The set of the positive roots is $\D^+=\D_0^+\cup\D_1^+$,
with the set $\D_0^+$ of positive
even roots and  the set $\D_1^+$ of the positive odd roots
 respectively given by
\begin{eqnarray*}
\D_0^+&\!\!\!=\!\!\!&\{\a_{i,j}=\es_i-\es_j,\,
\a_{\nu,\eta}=\es_\nu-\es_\eta
\,|\,1\le i<j\le m,\,\cp 1\le\nu<\eta\le\cp n\},
\\
\D_1^+&\!\!\!=\!\!\!&\{\a_{i,\nu}=\es_i-\es_\nu\,|\,i\in{\bf I}_1,\,
\nu\in{\bf I}_2\}.
\end{eqnarray*}
We define a total order on $\D_1^+$ by
\equa{0.0}
{
\a_{i,\nu}<\a_{j,\eta}\,\ \ \Lra\,\ \
\nu-i<\eta-j\mbox{ or }\nu-i=\eta-j\mbox{ but }i>j.
}
Then $\amin=\a_{m,\cp 1},\,\amax=\a_{1,\cp n}$
are respectively the minimal and maximal roots in $\D_1^+$.

Throughout the paper, we shall denote by
$\fg$ the special linear superalgebra ${\mathfrak{sl}}_{m|n}$,
which is the subalgebra of
${\mathfrak{gl}}_{m|n}$ consisting of supertraceless matrices.
Since ${\mathfrak{sl}}_{m|n}$ is isomorphic to ${\mathfrak{sl}}_{n|m}$,
we shall assume that $1\le n\le m$. We choose the Borel subalgebra
$\fb=B\cap\fg$ for $\fg$, which contains the Cartan subalgebra
$\fh=H\cap \fg$. We identify the dual space $\fh^*$ of $\fh$ with
the subspace $\sum_{a, b}\C(\epsilon_a-\epsilon_b)$ of $H^*$
spanned by the roots of ${\mathfrak{gl}}_{m|n}$. Then the roots of
$\fg$ coincide with those of ${\mathfrak{gl}}_{m|n}$, and relative
to $\fb$, a root $\a$ is positive if and only if $\a\in\D^+$.

The special linear superalgebra admits a $\Z$-grading
$\fg=\fg_{-1} \oplus \fg_{0}\oplus \fg_{+1}$, where $\fg_0=\fg_{\bar 0}
\cong\mathfrak{sl}(m) \oplus\C\check\rho_1\oplus\mathfrak{sl}(n)$ and
$\check\rho_1$ is defined in (\ref{rho1}).
Also,  $\fg_{\pm
1}\subset \fg_{\bar 1}$, with $\fg_{+1}$ being the nilpotent
subalgebra spanned by the odd positive root spaces, and $\fg_{-1}$
that spanned by the odd negative root spaces.
A basis of $\gl$ is given by
\equan{gl}
{
\{E_{ab},E_{aa}-(-1)^{[b]}E_{bb}\,|\,a,b\in
{\bf I},\, a\ne b\}.
}
We shall denote
\begin{eqnarray*}
&&
e_\a=E_{ab},f_\a=E_{ba}\mbox{ \ if \ }\a=\a_{a,b}\in\D^+=\D_0^+\cup\D_1^+,
\\ &&
h_i=E_{ii}-E_{i+1,i+1},\
h_0=E_{mm}+E_{\cp 1\cp 1},\
h_\nu=E_{\nu\nu}-E_{\nu+ 1,\nu+ 1},
\end{eqnarray*}
for $i\in{\bf I}_1\bs\{m\},\ \nu\in{\bf I}_2\bs\{\cp n\}.$

An element in $\fh^*$ is called a weight. A weight $\l$ is integral if
$(\lambda, \, \alpha)\in\Z$ for all roots, and
dominant if ${2(\lambda, \, \alpha)}/{(\alpha,\, \alpha)}\ge0$
for all positive even roots $\a$ of $\fg$.
A weight $\l\in\fh^*$ can be written in terms of $\es$-basis
\equa{weight1}
{\l=(\l_1,...,\l_m\,|\,\l_{\cp 1},...,\l_{\cp n})=
\mbox{$\sum\limits_{a\in {\bf I}}$}\l_a\es_a
\mbox{ such that }
\mbox{$\sum\limits_{i=1}^m$}\l_i+
\mbox{$\sum\limits_{\nu=\cp 1}^{\cp n}$}\l_\nu=0,
}
or in terms of Dynkin labels
\equan{weight2}
{
\l=[a_1,...,a_{m-1};a_0;a_{\cp 1},...a_{\cp n- 1}],
}
where
\equan{weight3}
{
a_i=\l_i-\l_{i+1},\ a_0=\l_m+\l_{\cp 1},\
a_\nu=\l_{\nu}-\l_{\nu+ 1},
}
for $i\in{\bf I}_1\bs\{m\},\ \nu\in{\bf I}_2\bs\{\cp n\}.$
We call $\l_p$ the $p$-th coordinate of $\l$ for $p\in {\bf I}$, and
$a_p$ the $p$-th Dynkin label of $\l$ for
$p\in {\bf I}\cup \{0\}\bs\{m,\cp n\}$.

The following weights will appear frequently in the remainder of the
paper:
\equa{notations1}
{
\mu^{(i,j)}=
(\stl{i-1}{\ob{0,...,0}},j+1,\stl{m-n+j}{\ob{1,...,1}},0,...,0\,
|\,0,...,0,\stl{j}{\ob{-1,...,-1}},-j-1-m+n,\stl{i-1}{\ob{0,...,0}}),
}
for $1\le i\le n,\,0\le j\le n+1-i$;
\equa{notations2}
{
\biggl\{
\begin{array}{ll}
\mu^{(j)}=\mu^{(1,j)},&
\mu_{\pm}^{(j)}=\mu^{(j)}\pm\amax,
\vs{4pt}\\
\ETc=(1,1,0,...,0\,|\,0,...,0,-1,-1),&
\ETd=(0,...,0,-1,-1\,|\,1,1,0,...,0),
\end{array}
}
where $\ETc,\,\ETd$ can occur only when $n\ge2$;
\begin{eqnarray}
\label{rho0}
\rho_0
&\!\!\!=\!\!\!&\frac{1}{2}\big(\mbox{$\sum\limits_{i=1}^m$}(m-2i+1)\es_i+
\mbox{$\sum\limits_{\nu=1}^n$}(n-2\nu+1)\es_{\cp \nu}\big)
\nonumber\\&\!\!\!=\!\!\!&\frac{1}{2}(m-1,m-3,...,-m+1\,|\,n-1,n-3,...,-n+1),
\\
\label{rho1}
\rho_1
&\!\!\!=\!\!\!&\frac{1}{2}\big(n\mbox{$\sum\limits_{i=1}^m$}\es_i
-m\mbox{$\sum\limits_{\nu=1}^n$}\es_{\cp\nu}\big)
=
\frac{1}{2}(n,...,n\,|\,-m,...,-m),\\ \nonumber
\rho&=&\rho_0-\rho_1,
\end{eqnarray}
where $\rho_0$ (resp.~$\rho_1$) is half the sum of positive even (resp.~odd) roots.

For every integral dominant weight $\l$, we denote by $L^{(0)}_\l$
the finite-dimensional irreducible $\fg_0$-module with highest
weight $\l$. Extend it to a $\gl_0\oplus\gl_{+1}$-module by
putting $\gl_{+1}L^{(0)}_\l=0$. Then the Kac module $V_\l$ is the
induced module \equa{Kac module} { V_\l={\rm
Ind}_{\gl_0\oplus\gl_{+1}}^{\gl}L^{(0)}_\l\cong U(\gl_{-1})
\otimes_{\C}L^{(0)}_\l. } Denote by $L_\l$ the irreducible module
with highest weight $\l$ (which is the unique irreducible
quotient module of $V_\l$) and we always fix a highest weight
vector $v_\l$.

For any finite-dimensional highest weight $\gl$-module $V$, we can
decompose $V$ into a direct sum of $\gl_0$-submodules with respect
to its level (an element $x\in\gl$ is said to have level $i$ if
$x\in\gl_i$ for $i=-1,0,1$, this defines a level structure on
$V$): \equa{level0} {
V=\mbox{$\sum\limits_{\level\in\Z}$}V^\level. } Set \equan{level1}
{ \top={\rm max}\{\level\in\Z\,|\,V^\level\ne0\}, \ \ \bottom={\rm
min}\{\level\in\Z\,|\,V^\level\ne0\}. }
When it is necessary to indicate the module $V$, we denote them by
$\top(V)$ and $\bottom(V)$.
 Then $\top-\bottom\le mn$.
In most cases, we shall specify the highest weight vector $v_\l$
to have level zero, then $\top=0$. But in some cases, we shall
shift level so that a vector with weight $0$ has level $0$.

\subsection{Lie superalgebra cohomology}
In this subsection we explain some basic concepts of
Lie superalgebra cohomology. The material can be found in many sources,
say, \cite{SZ98}. For $p\ge1$
and a finite-dimensional $\gl$-module $V$, let $C^p(\gl,V)$ (the
set of $p$-cochains)  be the $\Z_2$-graded vector space of all
$p$-linear mappings $\vp$ of $\gl^p=\gl\times\cdots\times\gl$ into
$V$ satisfying
$$
\vp(x_1,...,x_i,x_{i+1},...,x_p)=-(-1)^{[x_i][x_{i+1}]}
\vp(x_1,...,x_{i+1},x_i,...,x_p) \mbox{\ (super-skew-symmetry)}
$$
for $1\le i\le p-1$
(recall that we denote the degree of an element $x$ by $[x]\in\Z_2$).
Set $C^0(\gl,V)=V$. We define the differential operator $d:C^p(\gl,V)\to
C^{p+1}(\gl,V)$ by
\begin{eqnarray}
\label{differential}
&&\!\!\!\!\!\!\!\!(d\vp)(x_0,...,x_p)
\nonumber\\
&&=
\mbox{$\sum\limits_{i=0}^p$}(-1)^{i+[x_i]([\vp]+[x_0]+\cdots+[x_{i-1}])}x_i
\vp(x_0,...,\hat x_i,...,x_p)
\nonumber\\&&
+\mbox{$\sum\limits_{i<j}$}(-1)^{j+[x_j]([x_{i+1}]+\cdots+[x_{j-1}])}
\vp(x_0,...,x_{i-1},[x_i,x_j],x_{i+1},...,\hat x_j,...,x_p),
\end{eqnarray}
for $\vp\in C^p(\gl,V)$ and $x_0,...,x_p\in\gl$, where the sign
$\hat{\ }$ means that the element below it is omitted. It can be
verified that $d^2=0$. Set
$$
Z^p(\gl,V)=\ker{d|_{C^p(\gl,V)}},\ \ B^p(\gl,V)=\Im{d|_{C^{p-1}(\gl,V)}},\ \
H^p(\gl,V)=Z^p(\gl,V)/B^p(\gl,V).
$$
Elements in $Z^p(\gl,V)$ are called $p$-cocycles, elements  in
$B^p(\gl,V)$ are $p$-coboundaries. Two
elements of $Z^p(\gl,V)$ are said to be cohomologous if their
difference lies in $B^p(\gl,V)$. For $\vp\in Z^p(\gl,V)$, we
denote by $\ol\vp$ its residue class modulo $B^p(\gl,V)$. The
space $H^p(\gl,V)$ is called the $p$-th cohomology group.

Let $U,V,W$ be three $\gl$-modules such that
\equa{short-exact}
{
0\to U\rb{2pt}{\mbox{$\ ^{\ f}_{\dis\to}\ $}} V\rb{2pt}
{\mbox{$\ ^{\ g}_{\dis\to}\ $}}W\to0
}
is a short exact sequence, where $f,g$ are homogenous $\gl$-module homomorphisms.
Then there exists a long exact sequence
\equa{long-exact}
{
\cdots\to H^p(\gl,U)\rb{3pt}{\mbox{$\,\ ^{\ f^p}_{\dis\!\!\!-\!\!\!\to}\ $}} H^p(\gl,V)
\rb{3pt}{\mbox{$\ ^{\ g^p}_{\dis-\!\!\!-\!\!\!\to}\ $}}H^p(\gl,W)\to H^{p+1}(\gl,V)\to\cdots,
}
where the maps $f^p,\,g^p$ can be defined easily from $f,\,g$
 (cf.~\cite[(2.50)]{SZ98}).

For $x\in\gl$, we define $i_x:C^p(\gl,V)\to C^{p-1}(\gl,V),\,
\theta_x:C^p(\gl,V)\to C^p(\gl,V)$ by
\begin{eqnarray}
\label{i-x}
&&\!\!\!\!\!\!\!\!(i_x\vp)(x_1,...,x_{p-1})=(-1)^{[x][\vp]}\vp(x,x_1,...,x_{p-1}),
\\[4pt]
\label{theta-x}
&&\!\!\!\!\!\!\!\!(\theta_x\vp)(x_1,...,x_p)
\nonumber\\&&
=x\vp(x_1,...,x_p)
\nonumber\\&&
-\mbox{$\sum\limits_{i=1}^p$}(-1)^{[x]([\vp]+[x_1]+\cdots+[x_{i-1}])}
\vp(x_1,...,x_{i-1},[x,x_i],x_{i+1},...,x_p).
\end{eqnarray}
One can verify that
\equa{ix-thetax}
{
d i_x+i_x d=\theta_x,\ \ \ d\theta_x=\theta_x d.
}
Note that $\theta:x\mapsto\theta_x$ defines a $\gl$-module structure on $C^p(\gl,V)$
such that $Z^p(\gl,V)$, $B^p(\gl,V)$ are submodules by (\ref{ix-thetax}).
Since $\gl_0$ is a reductive Lie algebra, we can decompose
\equa{decompose-zp}
{
Z^p(\gl,V)=Z^p_0(\gl,V)\oplus B^p(\gl,V)
}
as a direct sum of $\gl_0$-submodules. For any $\gl$-module $V$, we denote
\equa{gl0-invariant}
{V^{\gl_0}=\{v\in V\,|\,\gl_0v=0\}.
}
Elements in $V^{\gl_0}$ are called $\gl_0$-invariant elements.
\begin{proposition}
\label{prop2.1}
Every $p$-cocycle is cohomologous to a $\gl_0$-invariant $p$-cocycle. More precisely,
\equa{equ-prop2.1}
{
Z^p_0(\gl,V)\subset Z^p(\gl,V)^{\gl_0}.
}
\end{proposition}
\begin{proof}
Let $\vp\in Z^p_0(\gl,V)$.
For all $x\in\gl_0$, we have, by (\ref{ix-thetax}),
$$
\theta_x\vp=di_x\vp+i_xd\vp=di_x\vp\in B^p(\gl,V)\cap Z^p_0(\gl,V)=\{0\}.
$$
This proves (\ref{equ-prop2.1}). Now the first statement follows from
(\ref{decompose-zp}).
\end{proof}

%
\section{First cohomology groups with coefficients in Kac modules} %
%
\label{1-Kac modules}
Let $\vp\in Z^1(\gl,V_\l)$ be a $1$-cocycle, by Proposition
\ref{prop2.1}, we can suppose $\vp\in Z^1(\gl,V_\l)^{\gl_{0}}$.
Denote $\vp^{(0)}=\vp|_{\gl_0}$. For $x\in\gl_{0},\,z\in\gl$, by
(\ref{differential}) and (\ref{theta-x}), we have \equa{1.1} {
0=d\vp(x,z)=
\theta_x\vp(z)-(-1)^{[z][\vp]}z\vp(x)=-(-1)^{[z][\vp]}z\vp(x), }
which implies that $\gl$ acts trivially on $\vp(x)$.

We shall divide the study into two cases according to the highest
weights of the Kac modules.

We first assume that $\l\ne2\rho_1$. Then $V_\l$ does not contain
a trivial $\gl$-submodule. So, $\vp^{(0)}=0$. By \refequa{level0},
for any $v\in V_\l$, we can uniquely write \equan{1.1+} {
v=\mbox{$\sum\limits_{\level=0}^{mn}$}v_{-\level},\mbox{ \ where \
} v_{-\level}\in V_\l^{-\level}. } Now (\ref{differential}) gives
\equa{1.2} {
f_\a\vp(f_\b)+f_\b\vp(f_\a)=(-1)^{\deg{\vp}}d\vp(f_\a,f_\b)=0
\mbox{ \ for \ }\a,\b\in\D_1^+. }  As $\U(\gl_{-1})$ is a
Grassmann algebra generated by the odd root vectors, and acts
freely on the Kac module $V_\l$ (cf.~(\ref{Kac module})), we can
take Grassmannian differentiation $\frac{\ptl}{\ptl f_\a}$ with
respect to $f_\a$, $\a\in \D_1^+$, to obtain \equa{1.3} {
\vp(f_\b)-f_\a\mbox{$\frac{\ptl}{\ptl
f_\a}$}\vp(f_\b)+\d_{\a,\b}\vp(f_\a) -f_\b\mbox{$\frac{\ptl}{\ptl
f_\a}$}\vp(f_\a)=0. } Sum over $\a\in\D_1^+$, we obtain \equa{1.4}
{
(mn+1)\vp(f_\b)-\mbox{$\sum\limits_{\level=0}^{mn}$}\level\vp(f_\b)_{-\level}
=f_\b v',\mbox{ \ where \ }
v'=\mbox{$\sum\limits_{\a\in\D_1^+}\frac{\ptl}{\ptl
f_\a}$}\vp(f_\a). } Note that $v'$ has degree $[\vp]$. Denote
\equan{1.5} {
v_{-\level}=(mn+1-\level)^{-1}(-1)^{\deg{\vp}}v'_{-\level+1},\mbox{
and } v=\mbox{$\sum\limits_{\level=0}^{mn}$}v_{-\level}. } Then by
\refequa{1.4}, \equa{1.6} { \vp(f_\b)=(-1)^{\deg{\vp}}f_\b v=d
v(f_\b)\mbox{ for all }\b\in\D_1^+. } Replacing $v$ by
$v-v_{-mn}$, we still have \refequa{1.6}, and we can suppose
\equa{1.6+} { v_{-mn}=0. } By \refequa{1.6}, and the
$\gl_0$-invariant property of $\vp$, we have \equan{1.6*} { f_\b
xv=xf_\b v-[x,f_\b]v=(-1)^{\deg{\vp}}(x\vp(f_\b)-\vp([x,f_\b]))=0
\mbox{ \ for \ }x\in\gl_0. } This together with equation
\refequa{1.6+} shows that $ \gl_0v=0. $ Thus by replacing $\vp$ by
$\vp-dv$, we still have that $\vp$ is $\gl_0$-invariant, and
\equa{1.6*3} {\vp(f_\b)=0\mbox{ \ for \ }\b\in\D_1^+. } Now using
the fact that for all $\a, \b\in\D_1^+$,  \equa{1.7} { 0
=(-1)^{\deg{\vp}}d\vp(e_\a,f_\b)= e_\a\vp(f_\b)+f_\b\vp(e_\a)
=f_\b\vp(e_\a), } we obtain that either $\vp(e_\a)\equiv0$, or for
some $\a$, $0\ne\vp(e_\a)\in V_\l^{-mn}$ (the bottom level) and
$V_\l^{-mn}\cong\gl_{+1}$ as $\gl_0$-modules.

In the former case, $H^1(\gl,V_\l)=0$. In the latter case,
\equan{lambda} {\l=2\rho_1+\amax=(n+1,n,...,n\,|\,-m,...,-m,-m-1),
}  and we only need to consider such $1$-cocycles $\phi$ that
satisfy the  conditions \equa{1.8}{
\left\{\begin{array}{l l l}
\phi|_{\gl_0\oplus\gl_{-1}}&=&0, \\
\phi(e_\a)&=& \mbox{\,the  image of $e_\a$ in $V_\l^{-mn}$
under the }\\
&&\mbox{$\gl_0$-module isomorphism $\gl_{+1}\cong V_\l^{-mn}$}.
\end{array}\right.}
Obviously, $\phi$ is $\gl_0$-invariant.  In fact, all 1-cochains
satisfying \eqref{1.8} are 1-cocycles. To prove this claim, it
only requires to verify the condition \equa{1.8*} {
e_\a\phi(e_\b)+e_\b\phi(e_\a)=0\mbox{ for all }\a,\b\in\D_1^+. }
Since the left hand side of \refequa{1.8*} is in
$V_\l^{\bottom+1}$, condition (\ref{1.8*}) holds if and only if
\equa{1.8**} { f_\g(e_\a\phi(e_\b)+e_\b\phi(e_\a))=0\mbox{ for all
}\a,\b,\g\in\D_1^+. } The left hand side of \refequa{1.8**} is
equal to
\begin{eqnarray}
\label{1.8*3}
&&
[f_\g,e_\a]\phi(e_\b)+[f_\g,e_\b]\phi(e_\a)
-e_\a f_\g\phi(e_\b)-e_\b f_\g\phi(e_\a)
\nonumber\\&&=
\phi([[f_\g,e_\a],e_\b]+[[f_\g,e_\b],e_\a])=0.
\end{eqnarray}
Thus $\phi$ is a 1-cocycle and it is non-trivial because
$\phi(e_\a)$ is in the bottom level of $V_\l$ which cannot be
written in the form $e_\a v$ for any $v\in V_\l$. Using the
$\gl_0$-invariant property of $\vp$, and the fact that $\gl_{+1}$
is irreducible as a $\gl_{0}$-module, we easily see that the space
of all the 1-cocycles satisfying \eqref{1.8} is 1-dimensional.
Thus \equa{1.9} { H^1(\gl,V_\l) \cong \C \mbox{ \ if \
}\l=2\rho_1+\amax. } Note that the above discussion actually
provides an explicit construction of the non-trivial 1-cocycles.

\bigskip
Now we consider the remaining case with $\l=2\rho_1$. This time
the Kac module $V_\l$ is free over $\U(\gl_{-1})$ of rank $1$. Let
$v_\l$ be a fixed non-zero $\gl$-highest weight vector of $V_\l$.
Then $\C\mbox{$\prod_{\a\in\D_1^+}$}f_\a v_\l$ forms a
$1$-dimensional $\gl$-submodule of $V_\l$, where the product
$\prod_{\a\in\D_1^+}f_\a$ is ordered with respect to the ordering
in \refequa{0.0}.

Consider any $\vp\in Z^1(\gl, V_\l)^{\gl_0}$. By equation
\eqref{1.1}, this $\vp$ must satisfy \equa{1.2.1} { \vp|_\glss=0
\mbox{ and
}\vp(2\check\rho_1)=c'\mbox{$\prod\limits_{\a\in\D_1^+}$}f_\a v_\l
\in\C\mbox{$\prod\limits_{\a\in\D_1^+}$}f_\a v_\l, } for some
$c'\in\C$, where $\glss$ denotes the semi-simple part of $\gl_0$.
Similar as before, we can suppose that \refequa{1.6*3} holds. Then
\equa{1.2.2} {
0=d\vp(e_\a,f_\b)=(-1)^{\deg{\vp}}f_\b\vp(e_\a)-\vp([e_\a,f_\b])
=(-1)^{\deg{\vp}}f_\b\vp(e_\a)-\d_{\a,\b}\vp(h_\a), } where the
last equality follows from (\ref{1.2.1}). Note that
$h_\a=(h_\a-\frac{1}{mn}(2\check\rho_1))
+\frac{1}{mn}(2\check\rho_1)$ and
$h_\a-\frac{1}{mn}(2\check\rho_1)\in\glss$, we thus obtain
\equan{1.2.3} {
f_\b\vp(e_\a)=(-1)^{\deg{\vp}}\d_{\a,\b}\frac{1}{mn}\vp(2\check\rho_1).
} In particular $f_\b\vp(e_a)=0$ for all $\b\ne\a$. Thus
$\vp(e_\a)$ has the form $c_\a\prod_{\g\in\D_1^+\bs\{\a\}}f_\g
v_\l$ for some $c_\a\in\C$. Combining this with \refequa{1.2.1},
we obtain \equa{1.2.4} {
\vp(e_\a)=c'c_\a\mbox{$\prod\limits_{\g\in\D_1^+\bs\{\a\}}$}f_\g
v_\l, } where \equa{1.2.4p}{
c_\a=(-1)^{\deg{\vp}+m(\a)}\frac{1}{mn},\mbox{ \ and \ }
m(\a)=\#\{\g\in\D_1^+\,|\,\g<\a\}. }

Now if we define a 1-cochain $\phi'$ by setting \equa{1.2.5} {
\phi'|_{\glss\oplus\gl_{-1}}=0 \mbox{ \ and \ }
\phi'(2\check\rho_1)=\mbox{$\prod\limits_{\g\in\D_1^+}$}f_\g
v_\l,\ \
\phi'(e_\a)=c_\a\mbox{$\prod\limits_{\g\in\D_1^+\bs\{\a\}}$}f_\g
v_\l, }
then one can check that $\phi'$ is a $\gl_0$-invariant $1$-cochain
(note that $\C v_\l$ is a trivial $\glss$-module). To verify that
it is a cocycle, we only need to show that condition (\ref{1.8*})
holds. This is indeed true as follows from (\ref{1.8*3}). Thus
$\phi'$ is a 1-cocycle. Furthermore, it is non-trivial, as can be
seen from the following arguments:  if $\phi'(e_\a)$ has the form
$e_\a v$ for some $v\in V_\l$, then $v$ is in the bottom level,
i.e., $v\in\C\prod_{\g\in\D_1^+}f_\g v_\l$. But
$\C\prod_{\g\in\D_1^+}f_\g v_\l$ forms a $1$-dimensional
$\gl$-submodule of $V_\l$, thus we must have $\phi'(e_\a)=e_\a v=0$, a
contradiction.

\bigskip
The preceding discussions in this section complete the study of
the first cohomology groups with coefficients in Kac modules. We
summarize the results below.

\begin{theorem}
\label{theo3.1} Let $V_\l$ be the finite-dimensional Kac module
with highest weight $\l$. Then \equan{L1} { H^1(\gl,V_\l)\cong
\left\{\begin{array}{ll} \C &\mbox{if \ }\l=2\rho_1+\amax,
\vs{4pt}\\
\C &\mbox{if \ }\l=2\rho_1,
\vs{4pt}\\
0&\mbox{otherwise}.
\end{array}\right.
}
\end{theorem}
%
%

As an immediate consequence of the theorem, we have the following
result.
\begin{corollary}\label{coro3.2}
Suppose $V$ is any weight module {\rm(}not necessary finite dimensional$\sc\,${\rm)}
over $\gl$. If $V$ contains a submodule isomorphic to the Kac module
$V_{2\rho_1}$, then $H^1{\ssc\!}({\ssc_{\!}}\gl,{\!}V{\ssc\!})\!\ne\!0$.
\end{corollary}
\begin{proof}
Suppose $v_{2\rho_1}\in V$ is a highest weight vector which generates the
Kac module $V_{2\rho_1}$. We define a linear map $\phi':\gl\to V$
by (\ref{1.2.5}). Obviously, $\phi'$ is a $\gl_0$-invariant
$1$-cocycle. We claim that it is a non-trivial cocycle. Otherwise
$\phi'=d v$ for some $\gl_0$-invariant vector $v\in V$. This in
turn leads to \equan{Coro1.0} {
\mbox{$\prod\limits_{\a\in\D_1^+}$}f_\a
v_{2\rho_1}=\phi'(2\check\rho_1)=2\check\rho_1 v=0, } which is a
contradiction.
\end{proof}
\begin{corollary}\label{coro3.3}
$H^1(\gl,U(\gl))\ne0$.
\end{corollary}
\begin{proof}
Let $v_{2\rho_1}=\prod_{\a\in\D_1^+}e_\a\in U(\gl)$ with weight
$2\rho_1$. Obviously $v_{2\rho_1}$ is a $\gl$-highest weight
vector of $U(\gl)$, thus generating a highest weight module $V$,
which is a quotient of the Kac module $V_{2\rho_1}$. We claim that
$V=V_{2\rho_1}$. If not, then $V$ does not contain the bottom composition factor
(which is the trivial module)
of $V_{2\rho_1}$, i.e.,
$\prod_{\a\in\D_1^+}\ad{f_\a}v_{2\rho_1}=0$, where ad denotes the
adjoint action. But we have \equan{Cor1.ad_rep} {
\mbox{$\prod\limits_{\a\in\D_1^+}$}\ad{f_\a}v_{2\rho_1}=
\mbox{$\prod\limits_{\a\in\D_1^+}$}\ad{f_\a}\mbox{$\prod\limits_{\a\in\D_1^+}$}e_\a=
\pm\mbox{$\prod\limits_{\a\in\D_1^+}$}h_\a+\cdots\ne0, } where
$h_\a=[e_\a,f_\a]$. A contradiction.
\end{proof}
%
\section{Second cohomology groups with coefficients in Kac modules} %
%
\label{2-Kac module}
We turn to the computation of the second cohomology groups with
coefficients in the Kac modules. By Proposition \ref{prop2.1}, we
only need to consider $\gl_0$-invariant 2-cocycles.

Let $\vp\in Z^2(\gl,V_\l)^{\gl_0}$. Since $\glss$ is semi-simple,
we have $H^2(\glss,V_\l)=0$, i.e., $\vp|_{\glss\times\glss}=d\psi$
for some $\glss$-invariant 1-cochain $\psi\in C^1(\glss,V_\l)$.
Extend $\psi$ to $\psi\in C^1(\gl,V_\l)^{\gl_0}$ by setting
$\psi|_{\gl_{\bar1}\oplus\C\check\rho_1}=0$, and replace $\vp$ by
$\vp-d\psi$ we can suppose \equan{3.vp1} {
\vp|_{\glss\times\glss}=0. } Then from \equan{3.vp2} {
\begin{array}{ll}
0
\!\!\!\!&
=d\vp(x_1,x_2,z)
\vs{4pt}\\&
=\theta_{x_1}\vp(x_2,z)-\theta_{x_2}\vp(x_1,z)+
(-1)^{\deg{z}\deg{\vp}}z\vp(x_1,x_2)-\vp([x_1,x_2],z)
\vs{4pt}\\&
=-\vp([x_1,x_2],z),
\end{array}
} where $x_1,x_2\in\glss,\,z\in\gl$, we obtain \equa{3.vp3} {
\vp|_{\glss\times\gl}=0. } {}From \equan{3.rho1} {
0=d\vp(\check\rho_1,\xi_1,\xi_2)=
\theta_{\check\rho_1}\vp(\xi_1,\xi_2)+
(-1)^{\deg{\vp}}(\xi_1\vp(\check\rho_1,\xi_2)+\xi_2\vp(\check\rho_1,\xi_1))
-\vp(\check\rho_1,[\xi_1,\xi_2]), } where
$\xi_1,\xi_2\in\gl_{\bar1}$, we obtain \equa{3.rho2} {
\xi_1\vp(\check\rho_1,\xi_2)+\xi_2\vp(\check\rho_1,\xi_1)=0, }
(note that $[\xi_1,\xi_2]\in\gl_0=\glss\oplus\C\check\rho_1$).
This together with equation (\ref{3.vp3}) shows that
$\psi':z\mapsto\vp(\check\rho_1,z)$,  $z\in \gl$,  is a
$1$-cocycle.

First assume $\l\ne2\rho_1, 2\rho_1+\amax$. Then it follows from
Theorem \ref{theo3.1} that $\psi'$ is a $1$-coboundary, i.e.,
there exists $v\in V_\l$ of degree $[\vp]$ such that
$\vp(\check\rho_1,z)=(-1)^{\deg{z}\deg{\vp}}zv$ for $z\in\gl$ and
such that $\gl_0v=0$ by (\ref{3.vp3}). We define a
$\gl_0$-invariant 1-cochain $\psi$ of degree $[\vp]$ by setting
$\psi|_{\glss\oplus\gl_{\bar1}}=0$ and $\psi(\check\rho_1)=v$.
Then by replacing $\vp$ by $\vp-d\psi$, we have \equa{3Case1.1} {
\vp|_{\gl_0\times\gl}=0. }

Now $d\vp(f_\a, f_\b, f_\g)=0$, $\a,\b,\g\in\D_1^+$, leads to
\equa{3Case1.2} {
f_\a\vp(f_\b,f_\g)+f_\b\vp(f_\a,f_\g)+f_\g\vp(f_\a,f_\b)=0 \mbox{
\ for \ }\a,\b,\g\in\D_1^+. } By taking Grassmannian
differentiation $\frac{\ptl}{\ptl f_\a}$ and argue as in the
derivations of equations (\ref{1.2})--(\ref{1.6*3}), we can show
that $\vp$ is cohomologous to a $\gl_{0}$-invariant $2$-cocycle
which satisfies \eqref{3Case1.1} and vanishes on
$\gl_{-1}\times\gl_{-1}$. Thus we can assume that \equa{3Case1.2+}
{\vp(f_\a,f_\b)=0, \quad \forall \,\a, \b\in\D_1^+. } Under this
condition (\ref{differential}) gives \equa{3Case1.3} {
f_\b\vp(e_\a,f_\g)+f_\g\vp(e_\a,f_\b)=0, \quad \forall \,\a, \b,\g
\in\D_1^+.} {}Again the same arguments as in the derivations of
(\ref{1.2})--(\ref{1.6*3}) renders $\vp$ satisfying the following
equation \equa{3Case1.3+} { \vp(e_\a,f_\b)=0, \quad \forall \,\a,
\b\in\D_1^+. } Then (\ref{3Case1.2+}) and (\ref{3Case1.3+}) gives
\equa{3Case1.4} { f_\a\vp(e_\b,e_\g)=0. } Thus $\vp(e_\a, e_\b)$
is in the bottom level of $V_\l$. By super-skew-symmetry,
$\vp|_{\gl_{+1}\times\gl_{+1}}$ is in fact a $\gl_0$-invariant map
from $\gl_{+1}\wedge\gl_{+1}$ to the bottom level of $V_\l$. (Here
$\wedge$ means symmetric tensor product as $\gl_{+1}$ is odd.)
Thus the problem of finding non-trivial $2$-cocycles is now
reduced to the determination of such maps. Note that
\equan{gl1-gl1} {
\gl_{+1}\wedge\gl_{+1}=L^{(0)}_{\ETa}\oplus
L^{(0)}_{\ETc}, \mbox{ as $\gl_0$-modules}, }
(cf.~notations (\ref{notations2})). If $\l\ne 2\rho_1+\ETa,
2\rho_1+\ETc$, the space of such maps is zero, thus
$H^2(\gl,V_\l)=0$.

For the remaining two cases with $\l=2\rho_1+\ETa$ and
$\l=2\rho_1+\ETc$ respectively, the space $S$ of
$\gl_0$-invariant maps from $\gl_{+1}\wedge\gl_{+1}$ to the bottom
level of $V_\l$ is $1$-dimensional. Let $\omega$ be the generator
of this space $S$. Now we can construct a $2$-cocycle $\phi_{2,1}$ as
follows:  set \equa{3Case1.5} {
\phi_{2,1}|_{\gl_0\times\gl_{\bar1}\oplus\gl_{+1}\times\gl_{-1}}=0,\quad
\phi_{2,1}|_{\gl_{+1}\times\gl_{+1}}=\omega. } Then
$H^2(\gl,V_\l)=\C\ol{\phi_{2,1}}$.

Next suppose $\l=2\rho_1$. We consider $\vp\in Z^2(\gl,
V_\l)^{\gl_0}$ satisfying \eqref{3.vp3}. The above arguments show that $\vp|_{\fg_{+1}\wedge\fg_{+1}}=0$
and so $\vp$ is
uniquely determined by $i_{2\check\rho_1}\vp$, which is a
$\gl_0$-invariant $1$-cochain. With the help of equation
(\ref{3.rho2}) we can show that $i_{2\check\rho_1}\vp$ is closed.
Thus by Theorem \ref{theo3.1} there exists a constant $c\in \C$
and a $\gl_0$-invariant $v\in V_\l$ ($v$ necessarily has weight
$0$) such that \equa{3Case2.1} { i_{2\check\rho_1}\vp=c\phi'+ d v,
} where $\phi'$ is defined by equation \eqref{1.2.5}. As
$\vp(2\check\rho_1, 2\check\rho_1)=0$, while
$\phi'(2\check\rho_1)\ne 0$, the constant $c$ must vanish.
Therefore,
$$i_{2\check\rho_1}\vp=d v.$$ We can always express $v$ as
$i_{2\check\rho_1}\psi$ for some $\gl_0$-invariant $1$-cochain
$\psi$ such that $d \psi$ satisfies equation \eqref{3.vp3}.
Therefore,
$$ i_{2\check\rho_1}(\vp-d\psi)=0, $$
which implies that $\vp=d\psi$. Hence $H^2(\gl,V_{2\rho_1})=0$.

Finally we consider the case with $\l=2\rho_1+\amax$. Any $\vp\in
Z^2(\gl, V_\l)^{\gl_0}$ satisfying \eqref{3.vp3} is uniquely
determined by $i_{2\check\rho_1}\vp$. We can argue as in the
preceding paragraph to show that
$$i_{2\check\rho_1}\vp = c \phi + d v, $$
where $c\in\C$, and $\phi$ is defined by \eqref{1.8}. Here $v\in
V_\l$ is a $\gl_0$-highest weight vector of weight $0$. But $V_\l$
does not have such elements (the difference of levels of weight
$\l$ and weight $0$ is bigger than $mn$), we have $v=0$. Thus
\equa{3Case3.1} {i_{2\check\rho_1}\vp = c \phi.} As
(\ref{3Case1.2}) still holds, we can assume that $\vp$ satisfies
(\ref{3Case1.2+}). Now $d\vp(e_\a, f_\b, f_\g)=0$, for $\a, \b,
\g\in\D_1^+$, leads to \equa{3Case3.2} {
f_\b\vp(e_\a,f_\g)+f_\g\vp(e_\a,f_\b)=\pm\vp([e_\a,f_\b],f_\g)\pm\vp([e_\a,f_\g],f_\b),
} where the sign on the right hand side depends on the parity of
$\vp$. Because of \eqref{1.8} and \eqref{3.vp3}, the right hand
side vanishes identically. Thus again we can suppose that equation
(\ref{3Case1.3+}) holds. In this case, we have from $d\vp(f_\a,
e_\b, e_\g)=0$, for $\a, \b, \g\in\D_1^+$, the following equation
\equa{3Case3.3} {f_\a\vp(e_\b,e_\g)-(-1)^{\deg{\vp}}
(\vp([f_\a,e_\b],e_\g) + \vp([f_\a,e_\g],e_\b))=0. }

Choose a basis $\{e'_\a\,|\,\a\in\D_1^+\}$ for $V_\l^\top$ such
that $\prod_{\tau\in\D_1^+}f_\tau e'_\a$ is the image of $e_\a$
under the $\gl_0$-module isomorphism $\gl_{+1}\cong V_\l^\bottom$.
Using \eqref{1.8} and \eqref{3.vp3} we can derive from
\eqref{3Case3.3} the following result \equan{3Case3.4} {
\vp(e_\b,e_\g)=c\big(c_\b\mbox{$\prod\limits_{\tau\in\D_1^+\bs\{\b\}}$}f_\tau
 e'_\g+c_\g\mbox{$\prod\limits_{\tau\in\D_1^+\bs\{\g\}}$}f_\tau
 e'_\b\big),
} where $c_\b$, $\b \in \D_1^+$ are the constants defined by
equation \eqref{1.2.4p}. Thus if we define \equa{define-phi22}
{\biggl\{\begin{array}{l l}
&\phi_{2,2}|_{\glss\times\gl\oplus\gl_{-1}\times\gl_{+1}}=0,\quad
i_{2\check\rho_1} \phi_{2,2}=\phi, \\[4pt]
&\phi_{2,2}(e_\b,e_\g)=c_\b\mbox{$\prod\limits_{\tau\in\D_1^+\bs\{\b\}}$}f_\tau
 e'_\g+c_\g\mbox{$\prod\limits_{\tau\in\D_1^+\bs\{\g\}}$}f_\tau
 e'_\b,\end{array}
} we obtain $H^2(\gl,V_\l)=\C\ol{\phi_{2,2}},$ as it can be easily
shown that $\phi_{2,2}$ is a non-trivial $2$-cocycle. Therefore, we
obtain the following result.

\begin{theorem}\label{theo4.1}
For the finite-dimensional Kac module $V_\l$, we have
\equan{Prop3}
{
H^2(\gl,V_\l)=\left\{\begin{array}{ll}
\C\ol{\phi_{2,1}}&\mbox{if \ }\l=2\rho_1+\ETa,\,2\rho_1+\ETc,
\vs{4pt}\\
\C\ol{\phi_{2,2}}&\mbox{if \ }\l=2\rho_1+\amax,
\vs{4pt}\\
0&\mbox{otherwise,}
\end{array}\right.
} where $\phi_{2,1},\,\phi_{2,2}$ are defined by $(\ref{3Case1.5})$
and  $(\ref{define-phi22})$. \hfill$\Box$
\end{theorem}

%
%
\section{First cohomology groups with coefficients in irreducible modules} %
%
\label{1-irr-modules}
In this section, we compute $H^1(\gl,L_\mu)$ for
$\gl={\mathfrak{sl}}_{m|n}$, where $L_\mu$ is the
finite-dimensional irreducible $\gl$-module with highest weight
$\mu$. We shall need the notion of primitive vectors, which we
recall:

\begin{definition}
\label{primitive-weight} For any $\gl$-module $V$, a non-zero
$\gl_0$-highest weight vector $v\in V$ of weight $\l$ is called a
primitive vector and $\l$ a primitive weight if $v$ generates an
indecomposable $\gl$-submodule and if there exists a
$\gl$-submodule $W$ of $V$ such that $v\notin W$ but $\gl_{+1}v\in
W$. If we can take $W=0$, then $v$ is called a strongly primitive
vector or a $\gl$-highest weight vector and $\l$ a strongly
primitive weight or a $\gl$-highest weight.
\end{definition}

Let us now begin the computation of the cohomology group. It is
known from \cite{FL, Fu} that $H^1(\gl,\C)=0$. (This can also be
easily proved by a direct computation by using Proposition
\ref{prop2.1}.)  Thus we suppose $\mu\ne0$.

Consider an arbitrary $\vp\in Z^1(\gl,L_\mu)^{\gl_0}$. We shall
write \equan{components}{\vp|_{\gl_0}=\vp^{(0)}, \quad
\vp|_{\gl_{\bar 1}}=\vp^{(1)}. } Equation \refequa{1.1} remains
valid in the present case, which implies that the image of
$\vp^{(0)}$ is a $\gl$-submodule of $L_\mu$ consisting of
invariants. Thus we have $\vp^{(0)}=0$.

We now turn to the consideration of $\vp^{(1)}$. We shall separate
this into several cases.

\subsection{The case with all $\vp(f_\a)\in L_\mu^\top$
and some $\vp(f_\a)\ne 0$.} This implies
that $\mu=-\amin$. Then from \equan{2.1} {
0=(-1)^{\deg{\vp}}d\vp(e_\a,f_\b)=e_\a\vp(f_\b)+f_\b\vp(e_\a)=f_\b\vp(e_\a),
\quad \a, \b\in\D_1^+, } it follows that non-zero $\vp(e_\a)$ must be in
the bottom level $L_\mu^{bottom}$ of $L_\mu$. This requires the
$\gl_0$-highest weight $\mu_\bottom$ of $L_\mu^{bottom}$ to be
$\amax$, because of the $\gl_0$-invariance of $\vp$. However, by
using \cite[Proposition 3.5]{VZ} or \cite[Theorem 3.5]{S}, we can
easily obtain $\mu_\bottom=(-1,...,-1,-n\,|\,m,1,...,1)\ne\amax$.
Thus $\vp(e_\a)=0$, for all $\a\in\D_1^+$.

Let us now define a $\gl_0$-invariant 1-cochain $\phi_1$ in the
following way: fix a $\gl_0$-module isomorphism $J: \gl_{-1}\cong
L_\mu^{\top}$, and let \equan{2.2} {
\phi_1|_{\gl_0\oplus\gl_{+1}}=0, \quad  \phi_1(f_\a)= J(f_\a).} To
verify that $\phi_1$ is a cocycle, we need to show that the
following condition \equan{2.3} {
f_\a\phi_1(f_\b)+f_\b\phi_1(f_\a)=0 \mbox{ for all
}\a,\b\in\D_1^+, } is satisfied, which is equivalent to \equa{2.4}
{ e_\g(f_\a\phi_1(f_\b)+f_\b\phi_1(f_\a))=0 \mbox{ for all
}\a,\b,\g\in\D_1^+. } Easy manipulations similar to the derivation
of \refequa{1.8*3} can show that \refequa{2.4} indeed holds. We
can also easily show that $\phi_1$ is a non-trivial $1$-cocycle.
Clearly, $\vp$ must be a scalar multiple of $\phi$. Therefore
\equa{2.5} { H^1(\gl,L_\mu)=\C\ol{\phi_1}\mbox{ \ if \
}\mu=-\amin. }

\subsection{The case with all $\vp(e_\a)\in L_\mu^\bottom$ and some $\vp(e_\a)\ne 0$ }
This implies
that
the $\fg_0$-highest weight $\mu_{\bottom}$ of $L_\mu^\bottom$ is $\amax$, thus
$\mu=\mu^{(n-1)}$ (cf.~notations (\ref{notations2})) by
\cite[Proposition 3.5]{VZ} or \cite[Theorem 3.5]{S}. Similar as
above, we can construct a non-trivial $\gl_0$-invariant 1-cocycle
$\phi'_1$ in the following way: fix a $\gl_0$-module isomorphism
$J': \gl_{+1}\cong L_\mu^{\bottom}$, and set \equa{2.2*} {
\phi'_1|_{\gl_0\oplus\gl_{-1}}=0,  \quad \phi'_1(e_\a)= J'(e_\a).}
Then \equa{2.5*} { H^1(\gl,L_\mu)=\C\ol{\phi'_1}\mbox{ \ if \
}\mu=\mu^{(n-1)}. }

\subsection{The remaining case} Finally we assume that
$\vp(e_\a)\notin L_\mu^\bottom$ and $\vp(f_\a)\not\in
L_\mu^\top$ for any non-zero $\vp(e_\a)$ and $\vp(f_\a)$.
In particular the given condition implies
that $\mu\ne-\amin, \  \mu^{(n-1)}$.

In order for $H^1(\gl,L_\mu)\ne0$, all the central elements of
$U(\gl)$ contained in $\gl U(\gl)$ must act on $L_\mu$ by
zero \cite[Prop. 2.2]{SZ98}. Thus by equations (3.33)--(3.35)
of \cite{SZ98}, $\mu$ has the following form:
there exists some $k$ with $0\le k\le n$ such that
\begin{eqnarray}
\label{mu_form}
\mu
\!\!\!\!&=\!\!\!&
(\mu_1,...,\mu_m\,|\,\mu_{\cp 1},...,\mu_{\cp n})
\nonumber\\&=\!\!\!&
\biggl(
\mu_1,{\sc...},\mu_k,
\stl{m-n+j_k}{\ob{k,{\sc...},k}},
\stl{j_{k-1}}{\ob{k\!-\!1,{\sc...},k\!-\!1}},{\sc...},
\stl{j_1}{\ob{1,{\sc...},1}},\stl{j_0}{\ob{0,{\sc...},0}}\,\biggr|
\nonumber\\&&
\,\,\,\,\stl{j_0}{\ob{0,{\sc...},0}},\stl{j_1}{\ob{-1,{\sc...},-1}},
{\sc...},\stl{j_k}{\ob{-k,{\sc...},-k}},
-\mu_k\!-\!m\!+\!n,{\sc...},-\mu_1\!-\!m\!+\!n
\biggr),
\end{eqnarray}
where \equa{mu_form1} { \mu_1\ge...\ge \mu_k\ge k,\ \ \
j_0,j_1,...,j_k\ge0,\ \ \ j_0+j_1+\cdots+j_k=n-k. }
If $k=0$, we
shall regard the set $\{\mu_1,...,\mu_k\}$ as empty, and $\mu=0$.
Since we have assumed $\mu\ne0$, we have $k\ge1$.

Needless to say, not all the $L_\mu$ with highest weights $\mu$
belonging to the list (\ref{mu_form}) have non-trivial first cohomology.
We now device ways to eliminate all the weights with trivial first cohomology.

Note that an irreducible module $L_\mu$ can always be embedded in a
unique Kac module $V_{\L_\mu}$ as the minimal submodule, where $\L_\mu$
is uniquely determined by $\mu$, see \cite[Theorem 3.2]{S} or \cite[Conjecture 4.1]{VZ}
(which was proved in \cite[Main Theorem]{Br}).
Denote by $\mu^*$ the highest weight of the dual
module $L^*_\mu$ of $L_\mu$. A weight $\mu$ is self-dual if
$\mu^*=\mu$.

For any weight $\l=(\l_1,...,\l_m\,|\,\l_{\cp 1},...,\l_{\cp n})$
as in (\ref{weight1}), we define its level to be
\equan{weight's_level}
{\mbox{$
\level(\l)=\sum\limits_{i=1}^m\l_i.
$}}
Then $\top(L_\mu)=\level(\mu),\,\bottom(L_\mu)=\level(\mu_\bottom)$. By \cite[Theorem 3.2]{S},
\cite[Conjecture 4.1]{VZ} and \cite[Main
Theorem]{Br}, $\mu_\bottom=\L_\mu-2\rho_1$, thus $\bottom(L_\mu)=\level(\L_\mu)-mn$.

There exists an automorphism $\omega: \gl\rightarrow\gl$ of $\gl$
which interchanges $\C e_\a$'s and $\C f_\a$'s:
\equa{Auto-omega}{
\omega(e_\a)=-(-1)^i f_\a,\ \omega(f_\a)=-e_\a,\  \omega(h)=-h
\mbox{ for }\a\in\D_i^+,\,i\in\Z_2,\,h\in\fh.
}
Using this automorphism we may
define a new action of $\gl$ on the space $L_\mu$.  Under this new action,
the module becomes $L_\mu^*$. This in particular implies that
$H^p(\gl,L_\mu)\cong$ $H^p(\gl,L_\mu^*),\,p\ge0$. Therefore by considering $\mu^*$
instead of $\mu$ if necessary, we
can assume that $\top(L_\mu)\le -\bottom(L_\mu)$.
 This implies
\equa{level_condition} { \top(L_\mu)\le\frac{1}{2}(mn+1),
\quad \top(L_\mu)+\level(\L_\mu)\le mn. }

If $k=n$, then $\mu=(n,...,n\,|\,-m,...,-m)=2\rho_1$, which does not
satisfy (\ref{level_condition}). Thus $1\le k\le n-1$. This in particular
proves that $H^1(\gl,L_\mu)=0$ if $n=1$.

We can now suppose $n\ge2$. We shall find further conditions on
the highest weight $\mu$.


Regard $\U(\gl_{+1})$ as a $\gl_0$-module under the adjoint
action. We define a $\gl_0$-module homomorphism
$\si:U(\gl_{-1})\to L_\mu$ by \equa{homo} {
\si\big(\mbox{$\prod\limits_{\a\in S}$}e_\a\big)=
\mbox{$\prod\limits_{\a\in S\bs\{\a_S\}}$}e_\a\vp(e_{\a_S}), }
where $S$ is any subset of $\D_1^+$ with $\a_S$ being the largest
element,  and the product is in the order (\ref{0.0}). Using
$e_\a\vp(e_\b)=-e_\b\vp(e_\a)$, one can verify that \refequa{homo}
indeed defines a $\gl_0$-module homomorphism, which is non-zero
(and thus surjective) so
long as $\vp(e_{\a_{max}})\ne 0$.

By assumption, $\vp(\gl_{+1})\not\subset L_\mu^{bottom}$,
in particular $\vp(e_{\a_{max}})\ne 0$, thus the
top level $L_\mu^{top}$ of $L_\mu$ must be a $\gl_0$-submodule of
$\sigma(\U(\gl_{+1}))$. This in particular implies that $\mu$
coincides with some $\gl_0$ highest weight of $\U(\gl_{+1})$.
It is known from \cite[4.1.1]{Ho} that a $\gl_0$-highest weight in
$U(\gl_{+1})\cong\wedge(\gl_{+1})$ (the exterior algebra of the
vector space $\gl_{+1}$) has the following form \equa{gl1_highest}
{ \mu=(\mu_1,\mu_2,...,\mu_m\,|\,-\mu'_{\cp n},-\mu'_{\cp n-
1},...,-\mu'_{\cp 1}), } where $(\mu_1,\mu_2,...,\mu_m)$ is a
partition of $\level(\mu)$ satisfying \equa{3.1} {0\le
\mu_m\le...\le\mu_1\le n, } and $(\mu'_{\cp 1},\mu'_{\cp
2},...,\mu'_{\cp n})$ is the transpose partition
$(\mu_1,\mu_2,...,\mu_m)^T$ of $(\mu_1,\mu_2,...,\mu_m)$, i.e.,
\equa{3.2} { \mu'_\nu=\#\{i\,|\,1\le i\le m,\,\mu_i\ge \nu\}. }
Using these conditions,  we obtain that $\mu$
has the form (\ref{mu_form}) such that  (\ref{mu_form1}) holds,
and \equa{mu_condition} { \mu_1=j_1+\cdots+j_k+k,\ ...,\
\mu_k=j_k+k. } Note that all such $\mu$ have maximal atypicality
(i.e., it is $n$-fold atypical, cf.~\cite{SZ98}) and $\amin$ is the
first atypical root.

To get further conditions on $\mu$, we need information from Appendix \ref{Appendix}.
The concepts of atypicality matrices, northeast chains (NE) of a weight, and $n$-,
$q$-, $c$-relationships of atypical roots, etc. to be used below are all explained
in the Appendix.

By \refequa{Lambda_mu} and Lemma \ref{lemm5.1}(1), we obtain
\equan{Lambda_mu_level} {
\level(\L_\mu)=\top+\#NE_\mu\ge\top+\#P_\mu=mn+(m+n)\mu_m-\top. }
This together with condition (\ref{level_condition}) and Lemma
\ref{lemm5.1}(2) shows that \equa{Further} { \mu_m=0\mbox{ and
$\g_1$ is $c$- or $q$-related to any other atypical roots.} } When
$\mu$ has the form (\ref{mu_form}) satisfying (\ref{mu_form1}) and
(\ref{mu_condition}), condition (\ref{Further}) becomes
\equa{mu_condition1} { j_s+j_{s+1}+\cdots+j_k\le n-k-s\mbox{ for }
s=1,...,k. } This condition in particular shows that \equa{k} {
1\le k\le \frac{n}{2}, } because if $k>\frac{n}{2}$, we would take
$s=n-k+1$ in (\ref{mu_condition1}), then the left hand side of
(\ref{mu_condition1}) is $\ge0$, but the right hand side is $-1$,
a contradiction.

Let us summarize the preceding discussion into the following
Lemma.
\begin{lemma}
In order for $H^1(\gl, L_\mu)\ne 0$, the highest weight $\mu$ of $L_\mu$
must be of the form $(\ref{mu_form})$ and satisfy the conditions
$(\ref{mu_form1})$, $(\ref{mu_condition})$, $(\ref{mu_condition1})$ and
$(\ref{k}).$
\end{lemma}
For convenience, we call such a weight $\mu$ a
{$k$-fold permissible weight}.

Identify a positive odd root $\a_{i,\cp\nu}$ with the
$(i,\nu)$-position in the atypicality matrix $A(\mu)$. We observe
that the number of elements of $NE_\mu$ in $i$-th row is
\equan{Num_NE_r} { \#\{(i,\nu)\in NE_\mu\,|\,1\le\nu\le n\}
=n-\mu_i-\#\{\nu\,|\,\mu'_{\cp\nu}\ge m+1-i\}=n-\mu_i-\mu_{m+1-i},
} (recall (\ref{gl1_highest}) and  (\ref{3.2})), and the number of
elements of $NE_\mu$ in $\nu$-th column is \equan{Num_NE_c} {
\#\{(i,\nu)\in NE_\mu\,|\,1\le i\le m\} =m-\mu'_{\cp n+1-
\nu}-\#\{i\,|\,\mu_i\ge \nu\}=m-\mu'_{\cp n+
1-\nu}-\mu'_{\cp\nu}. } Thus, from this and (\ref{Lambda_mu}),
we have \equa{Lambda_mu1} { \L_\mu=2\rho_1-\mu^R, } where
$\mu^R=(\mu_m,...,\mu_1\,|\,-\mu'_{\cp 1},...,-\mu'_{\cp n})$ is
the reverse of $\mu$. Equation (\ref{Lambda_mu1}) shows that the
$\gl_0$-highest weight $\mu_\bottom$ of $L_\mu^\bottom$ is
$-\mu^R$ and thus the lowest weight of $L_\mu$ is $-\mu$. This
proves

\begin{lemma}\label{lemm5.2}
Suppose $H^1(\gl,L_\mu)\ne0$ and
$\mu\ne-\amin,\mu^{(n-1)}$. Then $L_\mu$ is self-dual, i.e., $\mu^*=\mu$.
\end{lemma}

Since $\L_\mu\ne 2\rho_1, \ 2\rho_1+\amax$ by (\ref{Lambda_mu1}),
Theorem \ref{theo3.1} shows that
$ 
{H^1(\gl,V_{\L_\mu})=0. }$  {}From the short exact sequence
$ 
{0\to L_\mu\to V_{\L_\mu}\to V_{\L_\mu}/L_\mu\to 0, } $ we obtain
a long exact sequence of cohomology groups as in
(\ref{long-exact}).
In particular, since $H^0(\gl,V_{\L_\mu})=(V_{\L_\mu})^{\gl}=0$
(where for a $\gl$-module $V$, we denote $V^\gl=\{v\in V\,|\,\gl
v=0\}$, the set of $\gl$-invariant elements of $V$), and
$H^1(\gl,V_{\L_\mu})=0$, we obtain \equa{H_mu} {
H^1(\gl,L_\mu)\cong
H^0(\gl,V_{\L_\mu}/L_\mu)=(V_{\L_\mu}/L_\mu)^\gl. }

\subsection{Main result on first cohomology groups}
Now we are ready to prove

\begin{theorem} \label{theo5.1}
Let $L_\l$ be the finite-dimensional irreducible $\gl$-module with
highest weight $\l$. Then ${\rm dim\,}H^1(\gl,L_\mu)\le1$ and
$H^1(\gl,L_\mu)\ne0$ if and only if $\mu$ is one of the following
$n+1$ weights: $-\amin$, $\mu^{(j)}$, $j=0,...,n-1$
{\rm(}note that $\mu^{(j)}$ is self-dual if
$j<n-1$ and $(\mu^{(n-1)})^*=-\amin)$.
\end{theorem}
\begin{proof}
By \cite[Conjecture 4.1]{VZ} and
\cite[Main Theorem]{Br}, we have that, the multiplicity of
$L_\l$ in the composition series of $V_{\L_\mu}$ is
$a_{\L_\mu,\l}=[V_{\L_\mu}:L_\l]\le1.$ In particular, taking
$\l=0$, we have \equan{Prop2.1} { {\rm
dim\,}(V_{\L_\mu}/L_\mu)^\gl\le[V_{\L_\mu},L_0]\le1. } We want to
use \cite[Conjecture 4.1]{VZ} and \cite[Main Theorem]{Br} to prove \equa{Prop2.2} {
[V_{\L_\mu},L_0]=1 } for the weights listed in the theorem. Recall
(\ref{rho0}) and (\ref{rho1}),
and that the $n$-fold atypical weight $0$ has the following
atypical roots \equan{Prop2.4} { \g_1^0=\a_{m,\cp 1}, \
\g_2^0=\a_{m-1,\cp 2},...,\g_n^0=\a_{m+1-n,\cp n}, } with the
corresponding data (see \cite[Conjecture 4.1]{VZ}) \equan{Prop2.5}
{ (k_1,k_2,...,k_n)=(m+n-1,m+n-3,...,m-n+1), }
where $k_s=\#NE_0(s)$ is the number of the $s$-th northeast
chain $NE_0(s)$ of weight $0$ (cf.~Definition \ref{ne-chain} in Appendix \ref {Appendix}).
We take
\\\cl{$\theta=(\theta_1,\theta_2,...,\theta_n)=
(\stl{j_0}{\ob{1,{\sc...},1}},0,
\stl{j_1}{\ob{1,{\sc...},1}},0,...,0,\stl{j_{n-k}}{\ob{1,{\sc...},1}}).
$}\\
Then
\\\cl{$ \l_\theta=\dot
d\big(\mbox{$\sum\limits_{s=1}^n$}\theta_s k_s\g^0_s\big)
=d\big(\mbox{$\sum\limits_{s:\theta_s=1}\,
\sum\limits_{\a\in NE_0(s)}$}\a\big)
=2\rho_1-(\mu_m,...,\mu_1\,|\,-\mu'_{\cp 1},...,-\mu'_{\cp
n})=\L_\mu, $}\\
 where in general for a weight $\l$, we define $\dot
d(\l)=d(\l+\rho)-\rho$, and $d(\l)$ is the unique dominant weight
in $W\l$ (where $W$ is the Weyl group of $\gl$). Thus by
\cite[Conjecture 4.1]{VZ} and \cite[Main Theorem]{Br}, (\ref{Prop2.2}) holds.

Let $v_0\in V_{\L_\mu}$ be a primitive vector of weight $0$
corresponding to the composition factor $L_0$. Let
$V:=U(\gl)v_0=U(\gl_{\bar1})v_0$. Assume that
$\gl_{+1}v_0\not\subset L_\mu$. Let \equa{Let_tau} { \mbox{
$v_\tau\in U(\gl_{+1})v_0\bs L_\mu$ be a primitive vector of
maximal level with weight $\tau$. } } Then $\tau$ must also have
the form (\ref{mu_form}). Since $\gl_{+1}U(\gl_{+1})v_0$ is
isomorphic to a quotient of $U(\gl_{+1})$ as a $\gl_0$-module,
$\tau$ also has the form (\ref{gl1_highest}) satisfying condition
(\ref{3.2}) and (\ref{mu_condition}). Thus as proved before, we
have \equa{lambda_form} { \mbox{either }\tau=\mu^{(n-1)}\mbox{ or
$\tau$ is a self-dual weight.} }

Let $-\si$ be any anti-primitive weight in $U(\gl_{-1})v_0$ (i.e.,
$-\si$ is the lowest weight of a composition factor). Then as
before, $\si$ is self-dual or $\si=\mu^{(n-1)}$. Since the lowest
weight of $V_{\L_\mu}$ is $-\mu$, we must have that $\mu-\si$ is a
sum of distinct positive odd roots. Thus $\si$ has to be a
self-dual weight (this also shows that there does not exist a
primitive weight $\l$ in $\gl_{-1}U(\gl_{-1})v_0$). Let $V'$ be
the submodule of $V$ generated by self-dual primitive vectors.
Then \equa{f-v-0} { \gl_{-1}v_0\subset V' } First suppose
$\tau=\mu^{(n-1)}$. Since there is a unique way to write $\tau$ as
a sum of distinct positive odd roots, we must up to a non-zero
scalar have \equa{v_l1} {\mbox{$
v_\tau=\prod\limits_{\a\in\G_\tau}e_\a v_0, $}} where
$\G_\tau=\{\a_{1,\cp n},\a_{2,\cp n}..., \a_{m,\cp n},\a_{m,\cp
n- 1},...,\a_{m,\cp 1}\}$. Note that $\L_\mu-\tau$ can also be
uniquely written as a sum of distinct positive odd roots. Thus we
can up to scalars uniquely write \equa{v_l} {\mbox{$
v_\tau=\prod\limits_{\a\in \G'_\tau}f_\a v_{\L_\mu}+..., $}} where
$\G'_\tau$ is the unique subset of $\D_1^+$ such that
$\L_\mu-\tau=\sum_{\a\in\G'_\tau}\a$ and where the missing terms
are in $U(\gl_{-1})U(\gl_0^-)\gl_0^-v_{\L_\mu}$ (where $\gl_0^-$
is the negative part of the triangular decomposition of $\gl_0$).
Applying any $f_{\a_{i,\cp \nu}}$ to (\ref{v_l1}) for $i>1$ and
$\nu<n$, since $f_{\a_{i,\cp\nu}}$ commutes with
$e_\a,\,\a\in\G_\tau$, by (\ref{f-v-0}), we see that $f_{\a_{i,\cp
\nu}}v_\tau\in V'$. Clearly, by (\ref{v_l}), we have
$w=f_{\a_{i,\cp \nu}}v_\tau\ne0$ for some $i>1$ and $\nu<n$ (since
$\L_\mu\ne2\rho_1$, we must have that $\G'_\tau$ is a proper
subset of $\{\a_{i,\cp \nu}\,|\,i>1,\,\nu<n\}$). By applying
$e_\a,\a\in\D_0^+$ to $w$ until it becomes a $\gl_0$-highest
weight vector, we obtain that $f_{\a_{j,\cp \eta}}v_\tau$ is a
$\gl_0$-highest weight vector of weight $\tau-\a_{j,\cp \eta}$ in
$V'$ for some $j>1,\,\eta<n$. Therefore, there exists some
primitive weight $\si$ of $V'$ such that $\si-(\tau-\a_{j,\cp
\eta})$ must be a sum of distinct positive odd roots
(cf.~\cite[Lemma 5.2]{SHK}), which is impossible, because
$\si_1\le n-1$ and $\tau_1=n$ and $\si-(\tau-\a_{j,\cp
\eta})=(-1,...\,|\,...)$ (the first coordinate is $-1$) cannot be
a sum of distinct positive odd root. Thus $\tau$ is a self-dual
weight.

By considering the lowest weights, we see that $\mu-\tau$ is a sum
of distinct positive odd roots. Thus
\begin{lemma}
\label{lemm5.3}
A weight $\tau$ appeared in $(\ref{Let_tau})$
must be a $k'$-fold permissible weight for some $k'\le k$ and
$\mu-\tau$ is a sum of distinct positive
odd roots.
\end{lemma}

Now we continue the proof of Theorem \ref{theo5.1} and divide it into 2 cases.
\vskip6pt
{\bf Case 1}:  Suppose $\mu=\mu^{(j)}$ with $j\le n-2$.
\vskip6pt
By Lemma \ref{lemm5.3}, $\tau$ must be $1$-fold permissible weight and thus is
of the form
$\tau=\mu^{(\ell)}$ for some $\ell<j$. Let $V^*_{\L_\mu}$ denote
the dual module of the Kac module $V_{\L_\mu}$. Since
$-\mu=(\L_\mu-2\rho_1)^R$ (see (\ref{Lambda_mu1})) is the lowest
weight of $V_{\L_\mu}$ and the lowest weight vector can be
generated by all vectors, we see that $\mu$ is the highest weight
of $V^*_{\L_\mu}$ and the highest weight vector can generate every
vector, i.e., $V^*_{\L_\mu}=V_\mu$. Since $\tau$ is self-dual, we
have \equa{dual.0} { \mbox{$\tau$ is a primitive weight of
$V_{\L_\mu}$}
 \Lra
\mbox{$\tau^*\!=\!\tau$ is a primitive weight of
$V^*_{\L_\mu}\!=\!V_\mu$}. } We shall again use \cite[Conjecture
4.1]{VZ} and \cite[Main Theorem]{Br} to determine possible $\tau$ such that
$a_{\mu,\tau}=[V_\mu,L_\tau]$ $=1$. Since $\tau=\mu^{(\ell)}$, we
see that all atypical roots of $\tau$ are \equa{atypical_tau} {
\g_1^\tau=\a_{m,\cp 1},\g_2^\tau=\a_{m-1,\cp
2},...,\g_{n-1}^{\tau}= \a_{m-n+2,\cp n- 1},
\g_n^{\tau}=\a_{1,\cp n}, } with the corresponding data
$(k_1,...,k_n),\,k_s=\#NE_\tau(s)$ (cf.~Definition \ref{ne-chain})
being \equa{tau_data} { \left\{
\begin{array}{l}
k_s=m+n+1-2s\mbox{ \ if \ }1\le s\le n-\ell-2,
\vs{4pt}\\
k_{n-\ell-1}=1,
\vs{4pt}\\
k_s=m+n-1-2s\mbox{ \ if \ }n-\ell\le s\le n-1,
\vs{4pt}\\
k_n=1.
\end{array}
\right.
}
Let $\theta=(\theta_1,...,\theta_n)\in\Z_2^n$. We want to prove
\equa{theta_condition}
{
\l_\theta\ne\mu\mbox{ \ if  $\theta_s=1$ for some $s\ne n-\ell-1,n$}.
}
Denote $\ol\tau=\mu^{(\ell)}+\sum_{s=1}^n\theta_sk_s\g_s^\tau+\rho$.
Recall (\ref{notations1}), we obtain
that the $\cp n$-th coordinate of $\ol\tau$
is $-(m-n+\ell+1)-k_n+\rho_{\cp n}$, i.e.,
\equan{1_at_least}
{
\ol\tau_{\cp n}=-(m-n+\ell+1)-k_n+\rho_{\cp n},
}
and by
(\ref{atypical_tau}) and (\ref{tau_data}),
\equan{i_at_least}
{
\ol\tau_{\cp s}=
\biggl\{\begin{array}{ll}
0-(m+n+1-2s)+\rho_{\cp s}&\mbox{ if }s<n-\ell-1,
\mbox{ or }
\vs{4pt}\\
-1-(m+n-1-2s)+\rho_{\cp s}&\mbox{ if }s\ge n-\ell.
\end{array}
} Suppose $\ol\tau_{\cp n},\,\ol\tau_{\cp s}$ are respectively the
$\cp p$-th, $\cp q$-th coordinates of $d(\ol\tau)$ (the unique
dominant weight in the Weyl chamber of $\ol\tau$). Then $q<p$
since $\ol\tau_{\cp n}<\ol\tau_{\cp s}$. Then one can easily see
that the $\cp p$-th, $\cp q$-th coordinates of $\l_\theta$ are
respectively $\ol\tau_{\cp n}-\rho_{\cp p}=-(m-n+\ell+1)-k_n-n+p$
$\le-2$ and $\ol\tau_{\cp s}-\rho_{\cp q}\le -2$ (note that
$\rho_{\cp s}-\rho_{\cp q}=q-s$ by (\ref{rho0}) and (\ref{rho1})).
Thus we have (\ref{theta_condition}). So, we have proved that if
$\l_\theta=\mu$, then
$\theta=(0,...,0,\theta_{n-\ell-1},0,...,0,\theta_n)$. One can
easily check that the latter can happen if and only if $j\ge1$ and
$\ell=j-1$ and $\theta_{n-\ell-1}=\theta_n=1$. In particular, if
$j=0$, then no $\tau$ in (\ref{Let_tau}) can occur, i.e.,
$(V_{\L_\mu}/L_\mu)^\gl=\C$. In this case, Theorem \ref{theo5.1} follows
from (\ref{H_mu}).

So suppose $j\ge1$. Then we just proved that $\tau$, if occurred
in (\ref{Let_tau}), must be the form $\tau=\mu^{(\ell)}$ with
$\ell=j-1$. Consider the Kac module $V_\mu$. Suppose the
corresponding primitive vector of weight $\tau$ in $V_\mu$ is
$v'_\tau$. By (\ref{dual.0}), $0$ is also a primitive weight of
$V_\mu$, and suppose $v'_0$ is a corresponding primitive vector.
We claim that \equa{Claim.0} { v'_0\mbox{ cannot be generated by
}v'_\tau. } (In fact by a conjecture of \cite{HKV}, $v'_\tau$ is
generated by $v'_0$). Take \equa{define-Lambda} {
\begin{array}{ll}
\!\!\!\!\L\!
=\!\mbox{$\sum\limits_{s=1}^j$}\mu^{(s,j-s)}
\!=\!
(\stl{m-n+j+1}{\ob{j{\sc\!}+{\sc\!}1,...,j{\sc\!}+{\sc\!}1}},
0,...,0\,|\,0,...,0,
\stl{j+1}{\ob{n{\sc\!}-{\sc\!}m{\sc\!}-{\sc\!}j{\sc\!}-{\sc\!}1,
...,n{\sc\!}-{\sc\!}m{\sc\!}-{\sc\!}j{\sc\!}-{\sc\!}1}}).\!\!\!\!
\end{array}}
 By \cite{SHK}, we see that
$\mu$ and $0$ are strongly primitive weights of $V_\L$
corresponding to unlinked codes. Suppose
$v''_\mu$ and $v''_0$ are their corresponding strongly primitive
vectors in $V_\L$.

We prove that $v''_0$ can be generated by $v''_\mu$ as follows (in
fact, the code corresponding to $v''_\mu$ is a sub-code of the
code corresponding to $v''_0$ (see \cite[Theorem 3.7]{S})):
Consider the dual module $V^*_\L=V_\Si$, where $\Si=2\rho_1-\L^R$.
Note that  $L^*_\mu=L_\mu$, $L^*_0=L_0$ and $L^*_\L$ all are
composition factors of $V_\Si$. Let $v_\Si$ be the highest weight
vector of $V_\Si$. Let $v_{\mu_\bottom}\in V_\Si$ be up to scalars the
unique
$\gl_0$-highest weight vector of weight $\mu_\bottom=-\mu^R$
corresponding to the lowest $\gl_0$-highest weight vector of the
composition factor $L_\mu$ of $V_\Si$. First note that \equa{No1.}
{\mbox{$ \Si-0=2\rho_1-\L^R=\sum\limits_{\a\in\G_1}\a, $}} where
$\G_1$ is a unique subset of $\D_1^+$. Thus according to
\cite[Lemma 5.2]{SHK}, we can up to scalars uniquely write $v'''_0$
(a primitive vector of weight $0$ in $V_\Si$)
as
\equa{No1.1} {\mbox{$ v'''_0=\prod\limits_{\a\in\G_1}f_\a
v_\Si+\cdots, $}} where the missing terms are contained in
$U(\gl_{-1})U(\gl^-_{0})\gl^-_{0}v_\Si$. As in \cite{SHK}, we
shall call the first term of the right hand side of (\ref{No1.1})
the {leading term} of $v'''_0$. Also note that \equan{No2.}
{\mbox{$ 0-\mu_\bottom=\mu^R=\sum\limits_{\a\in\G_2}\a, $}} for a
unique subset $\G_2$ of $\D_1^+$. As in \cite[\S5]{SHK} (where
some operators $\chi_J$ are defined), we can construct a weight
vector $v'_{\mu_\bottom}$ of weight $\mu_\bottom$ from $v'''_0$
such that \equa{No2.1} {\mbox{$
v'_{\mu_\bottom}=\prod\limits_{\a\in\G_2}f_\a v'''_0+\cdots, $}}
and such that $v'_{\mu_\bottom}$ is a $\gl_0$-highest weight
vector as long as $v'_{\mu_\bottom}\ne0$. We claim that
$v'_{\mu_\bottom}$ is indeed non-zero: Substitute (\ref{No1.1})
into (\ref{No2.1}), we see that $v'_{\mu_\bottom}$ has the leading
term $\prod_{\a\in\G_1\cup\G_2}f_\a v_\Si$. This is because: for
any $\a=\a_{i,\cp p}\in\G_1$ and any $\b=\a_{j,\cp q}\in\G_2$, we
have either $i>j$ or $p<q$. Thus when we substitute (\ref{No1.1})
into (\ref{No2.1}), we only produce one leading term. Thus
the $\fg_0$-highest weight vector $v_{\mu_\bottom}=v'_{\mu_\bottom}$ of weight
$\mu_\bottom$ in $V_\Si$ can be generated by $v'''_0$. By
taking dual, we obtain that $v''_0$ in $V_\L$ can be generated by $v''_\mu$.

By \cite[Conjecture 4.1]{VZ} and \cite[Main Theorem]{Br}, we can prove that
$\tau=\mu^{(\ell)}$ is not a primitive weight of $V_\L$: Suppose
$\L=\l_\theta$ for some $\theta$, i.e., $\L+\rho=
d(\tau+\sum_{s=1}^n k_s\theta_s\g^\tau_s+\rho)$. This means that
two sets $\{\L_{\cp s}+\rho_{\cp s}\,|\,s=1,...,n\}$ and
$\{\tau_{\cp s}+k_s\theta_s+\rho_{\cp s}\,|\,s=1,...,n\}$ are
equal, and so by (\ref{rho0}), (\ref{define-Lambda}) and
(\ref{tau_data}), two sets \equan{2-sets1} {
\{1,2,...,n-j-1,m+1,m+2,...,m+j+1\} \mbox{ \ \ and} }
\equan{2-set2}
{
\begin{array}{l}
\{\theta_1(m+n-1)+1,\theta_2(m+n-3)+2,...,\theta_{n-j-1}(m-n+2j+3)+n-j-1,
\vs{4pt}\\ \ \ \
\theta_{n-j}+n-j,
\vs{4pt}\\ \ \ \
1+\theta_{n-j+1}(m-n+2j-3)+n-j+1,...,
1+\theta_{n-1}(m-n+1)+n-1,
\vs{4pt}\\ \ \ \
m-n+j+\theta_n+n\}
\end{array}
}
are equal. No matter whether
$\theta_{n-j}=0$ or $1$, the first set does not contain
$\theta_{n-j}+n-j$. Thus $\tau$ is not a primitive weight of $\L$.

Now we define the module homomorphism $\pi:V_\mu\to
u(\gl)v''_\mu\subset V_\L$ from Kac module $V_\mu$ to the highest
weight module $U(\gl)v''_\mu$ as follows: $\pi$ sends $v_\mu$ to
$v''_\mu$. Then $\pi$ sends $v'_\tau$ to zero
(since $\tau$ is not a primitive weight of $V_\L$)
but sends $v'_0$ to
$v''_0$ (since $v''_0$ is generated by $v''_\mu$, there is a
pre-image $v'_0$ of $v''_0$ which is also primitive
in $V_\L$). This
proves (\ref{Claim.0}).

Now consider $V^*_\mu=V_{\L_\mu}$. By (\ref{Claim.0}), we see that
in $V_{\L_\mu}$, $v_\tau$ cannot be generated by $v_0$. This
proves that no $\tau$ in (\ref{lambda_form}) can occur in
(\ref{Let_tau}), i.e., $(V_{\L_\mu}/L_\mu)^\gl=\C$.
By (\ref{H_mu}), this proves Theorem \ref{theo5.1} in this case.

\vskip6pt {\bf
Case 2}:  Suppose $\mu$ is a $k$-fold permissible weight with
$k\ge2$ (see Remark \ref{rema5.2} below). \vskip6pt

Note that we
can uniquely decompose $\mu$ as $\mu=\nu+\eta$, where \equan{mu-1}
{ \nu=(\mu_1,\stl{\mu'_1-1}{\ob{1,...,1}},0,...,0\,|\,0,...,0,
\stl{\mu_1-1}{\ob{-1,...,-1}},-\mu'_1) } is a $1$-fold permissible
weight and $\eta=\mu-\nu$ is a $(k-1)$-fold permissible weight
when restricting it to be a dominant weight of $sl(m-1/n-1)$
($\eta$ is not dominant as a weight of $\gl$). Using
\cite[Conjecture 4.1]{VZ} and \cite[Main Theorem]{Br},
we see that both $\nu$ and $0$ are
primitive weights of $V_\mu$, which correspond to linked codes
(cf.~\cite{HKV}). Let $w_\nu$ and $w_0$ be the corresponding
primitive vectors in $V_\mu$. Using exactly the same arguments in
the paragraph after (\ref{define-Lambda}),
we obtain that
${ \mbox{ $w_0$ can be generated by $w_\nu$. } }$ This means that
in $V^*_\mu=V_{\L_\mu}$, a primitive vector of weight $\nu$ can be
generated by a primitive vector of weight $0$, i.e., $L_0$ is not
a submodule in $V_{\L_\mu}/L_\mu$. Thus
$(V_{\L_\mu}/L_\mu)^\gl=0$. By (\ref{H_mu}),
this completes the proof of Theorem
\ref{theo5.1}. \end{proof}

\begin{remark}\label{rema5.0}
It is also possible to prove Theorem \ref{theo5.1} using the machinery to be developed
in Section \ref{2-irr-modules}{\rm)}.
\end{remark}

\begin{remark}\label{rema5.2}
Suppose $\mu$ is a weight such that
all the central elements of $U(\gl)$ contained in
$\gl U(\gl)$ act  trivially on $L_\mu$.
As in $(\ref{5.2-w1})$, one can prove that
$L_\mu$ contains a copy of $\gl_0$-module $\gl_{+1}$ or $\gl_{-1}$
if and only if $\mu=-\amin,\,\mu^{(j)},\,0\le j\le n-1$ $($this also
shows that $H^1(\gl,L_\mu)=0$ if $\mu$ is a $k$-fold permissible
weight with $k\ge2)$. Thus
\equa{amax-copy}{
H^1(\gl,L_\mu)\ne0\ \ \Lra\ \
\mbox{$L_\mu$ contains a copy of $\gl_0$-module $\gl_{+1}$ or $\gl_{-1}$}.
}
\end{remark}
%
%
\section{Second cohomology groups with coefficients in irreducible modules} %
%
\label{2-irr-modules}
In this section we compute the second cohomology groups of the special
linear superalgebra with coefficients in the finite dimensional irreducible modules.

\subsection{Primitive weight graphs}\label{graph}
We first introduce some concepts, which will be used extensively in the
remainder of the paper.
\begin{definition}\label{defi6.1}
For a $\gl$-module $V$, we denote by $P(V)$ the set of primitive
weights of $V$.  Let $P_0(V)=\{\mu\in P(V)\,|\,L_\mu$ is an
irreducible submodule of $V\}$. A weight in $P_0(V)$ is called a
{lowest} or {bottom primitive weight}. Let
$P'(V)=\{v_\mu\,|\,\mu\in P(V)\}$ be a collection of non-zero
primitive vectors. For a primitive weight $\mu$ of $V$, we denote
${ U(\mu)=U_V(\mu)=U(\gl)v_\mu. }$ If $\mu\ne\nu$ and $v_\nu\in U_V(\mu)$,
we say that $\nu$ is {derived from} $\mu$ and write $\nu\dlar\mu$
or $\mu\drar\mu$. If $\mu\drar\nu$ and there
exists no $\l$ such that $\mu\drar\l\drar\nu$, then we say
$\nu$ is {\it directly derived from} $\mu$ and write
$\mu\rrar\nu$ or $\nu\llar\mu$. If $\mu\rrar\nu$ and $v_\nu\in
U(\gl_{+1})v_\mu$, we also use $\mu\erar\nu$ or $\nu\elar\mu$ to
denote this fact; similarly, we also use $\mu\frar\nu$ or $\nu\flar\mu$
to denote $v_\nu\in U(\gl_{-1})v_\mu$ $($Sometimes for convenience,
we also use symbols $\mu\link\l$ to denote $\mu\rrar\l$ or
$\mu\llar\l)$.

We can associate $P(V)$ with a directed graph, still denoted by $P(V)$,
called the {\it primitive weight graph of $V$}, such that two weights $\l$ and $\mu$ are
connected by a single directed edge $($i.e., the two weights are
linked$)$ pointing toward $\mu$
if and only if $\mu$ is directly derived from $\l$.
A {\it subgraph of $P(V)$} is a
subset $S$ of $P(V)$ together with all edges linking elements of $S$.
A subgraph $S$ of $P(V)$ is {\it closed} if it satisfies the following condition:
For any $\eta\in P(V),\ \mu,\nu\in S$,
\equan{closed}
{
\mu\drar\eta\drar\nu
\ \ \ \ \ \mbox{implies} \ \ \ \ \eta\in S.
}
\end{definition}

It is clear that a module is indecomposable if and only if its
primitive weight graph is connected (in the usual sense). It is
also clear that a subgraph of $P(V)$ corresponds to a subquotient
of $V$ if and only if it is closed. Thus a subgraph with only $2$
weights is always a closed subgraph.

For any subset $S$ (not necessarily a subgraph), we denote by $\ol S$
the smallest closed subgraph which contains $S$.

For any graph $\G$, we denote by $M(\G)$ any module with primitive
weight graph $\G$ if it is indeed a primitive weight graph of a
module. If $\G$ is a subgraph of $P(V)$, then $M(\G)$ is defined, and
we shall also denote
$ 
{
\ol M(\G)=M(\ol \G).
}$
If $\G$ corresponds to a submodule  of $V$ and
we need to indicate $M(\G)$ as a subquotient
module of $V$, we denote $M(\G)$ by $M_V(\G)$.

For a dominant weight $\mu$, we let $P(\mu)=P(V_\mu)$, and set
\equan{4UP}
{
P^\vee(\mu)=\{\l\,|\,\mu\in P(\l)\}.
}
Then $P^\vee(\mu)$ is the set of dominant weights $\l$ such that every
Kac module $V_\l$ has a composition factor $L_\mu$.

We shall say that a weight $\mu$ has {non-zero
$1$-cohomology} if $H^1(\gl,L_\mu)\ne0$.

\begin{remark}\label{rema6.1}
{\rm(1)} Since we will frequently need to
determine primitive weights of Kac modules, it will
be convenient to first use permissible codes defined in \cite{HKV} to
find a possible primitive weight, then use \cite[{\it Conjecture} 4.1]{VZ}
and \cite[{\it Main Theorem}]{Br}
to check if it is a primitive weight.

{\rm(2)} Let $P(V)$ be the primitive weight graph of the module $V$.
The dual primitive weight graph
$P^*(V)$ of $V$ is the graph obtained from $P(V)$
by reversing the directions of all arrows and changing all
weights to their dual weights. Note that $P^*(V)=P(V^*)$,
where $V^*$ denote the dual module of $V$.
If we change the action of $\gl$ on $P(V^*)$ by the
automorphism $\omega$ defined in $(\ref{Auto-omega})$ which exchanges
$\C e_\a$'s and $\C f_\a$'s,
then we obtain another module, called
the {inverse module of} $V$, with graph $\wt P(V)$ obtained from $P(V)$ by
reversing the directions of all arrows
$($recall that using the automorphism $\omega$, the module $L_\mu$ becomes
$L_{\mu^*}=L^*_\mu$ for all $\mu)$.
In particular, we have
\equa{inverse}
{
\exists\,M(\mu\rrar\nu)\,\Lra\,\exists\,M(\mu\llar\nu)\,\Lra
\,\exists\,M(\mu^*\rrar\nu^*)\,\Lra\,\exists\,M(\mu^*\llar\nu^*).
}

{\rm(3)} If $(\ref{inverse})$ occurs and $\mu\ne\nu$, then $\nu\in
P(\mu)\cup P^\vee(\mu)$ since either $M(\mu\rrar\nu)$ or
$M(\mu\llar\nu)$ must be a highest weight module $($we adopt
the convention that a highest weight module is cyclically generated by
a highest weight vector$)$, and similarly,
$\nu^*\in P(\mu^*)\cup P^\vee(\mu^*)$ $($note that it is possible
that $\mu=\nu$, in this case, $M(\mu\rrar\mu)$ is not necessarily
a weight module$)$. {}From this we obtain that if $\mu$ is a
primitive weight of the highest $($resp. lowest$)$ level in an
indecomposable module $P(V)$ such that either every other primitive
weight derives $\mu$ or is derived from $\mu$, then
\equan{Remark2.3} { P(V)\subset P(\mu) \mbox{ \ $($resp. \
$P(V)\subset P^\vee(\mu)\, )$}. }
\end{remark}

\subsection{Technical lemmas}\label{Technical-lemmas}
This subsection contains a series of lemmas which will be used
in establishing Theorem \ref{theo6} in the next subsection. The proofs of
most of the lemmas rely on detailed analysis of structures of Kac
modules, which unfortunately is a matter of a very technical nature.

\begin{lemma}\label{lemm6.1}
$\mu\link0$ implies that $H^1(\gl,L_\mu)\cong\C$.
\end{lemma}

\begin{proof} Say we have a module $V=M(0\rrar\mu)$.
Then we can define a $1$-cocycle $\vp\in Z^1(\gl,L_\mu)$ by
$\vp(x)=x v_0\in L_\mu$ for $x\in\gl$, where $v_0$ is a primitive
vector with weight $0$. Clearly it is non-trivial,
otherwise $V$ is decomposable.
\end{proof}

\begin{lemma}\label{lemm6.2}
{\rm(1)} If there is a short exact sequence $0\to\C\to V\to W\to 0$, then
$H^1(\gl,V)\cong H^1(\gl,W)$.

{\rm(2)} If there is a non-split short
exact sequence $0\to W\to V\to\C\to 0$, then ${\rm dim\,}H^1(\gl,V)={\rm dim\,}H^1(\gl,W)-1$.
\end{lemma}

\begin{proof}
(1) Since $H^1(\gl,\C)=H^2(\gl,\C)=0$, by (\ref{long-exact}), we have
$$
0=H^1(\gl,\C)\to H^1(\gl,V)\to H^1(\gl,W)\to H^2(\gl,\C)=0,
$$
which gives $H^1(\gl,V)\cong H^1(\gl,V)$.

(2) Since the short exact sequence is non-split, we have the exact
sequence
$$
H^0(\gl,W) {\,}^{\ i}_{\dis\to}{\,}  H^0(\gl,V) {\,}^{\ j}_{\dis\to}{\,}H^0(\gl,\C)
\to H^1(\gl,W)\to H^1(\gl,V)\to H^1(\gl,\C),
$$
where $i$ is the identity map, and so $j$ is the zero map. Thus we
have an exact sequence $0\to\C\to H^1(\gl,W)\to H^1(\gl,V)\to0$,
which gives the result.
\end{proof}

\begin{lemma}\label{lemm6.3}
Suppose $\gl\ne sl(2/1)$, and suppose there is a short exact sequence of $\gl$-modules
\equan{4.4.1}
{
0\to\C\to V\to W\to 0\mbox{ \ or \ }0\to W\to V\to \C\to 0.
}
Then $H^2(\gl,V)\cong H^2(\gl,W)$.
\end{lemma}

\begin{proof} Note that $H^3(\gl,\C)=0$ as long as $\gl\ne sl(2/1)$.
We have the following exact sequence (cf.~(\ref{long-exact}))
\equan{4.4.2}
{
\begin{array}{l}
0=H^2(\gl,\C)\to H^2(\gl,V)\to H^2(\gl,W)\to H^3(\gl,\C)=0
\mbox{ \ or }
\vs{4pt}\\
0=H^1(\gl,\C)\to H^2(\gl,W)\to H^2(\gl,V)\to H^2(\gl,\C)=0,
\end{array}
}
which give the result.
\end{proof}

\begin{lemma}\label{lemm6.4}
For any $($finite-dimensional${\sc\,})$ module $V$, we have
\equa{Lemma7.1}
{
{\rm dim\,}H^1(\gl,V)\le \mbox{$\sum\limits_{\mu\in P(V)}$}{\rm dim\,}H^1(\gl,L_\mu).
}
More generally, suppose
\equa{Lemma7.1-}
{
P(V)=\biggl(\bigcup_{i=1}^p P_i\biggl)\bigcup\biggl(\bigcup_{j=1}^q Q_j\biggr)
}
is a disjoint union such that all $P_i,\,Q_j$ are closed subgraph of
$P(V)$ and such that all $Q_j$ do not contain any primitive weight
with non-zero $1$-cohomology, then
\equa{Lemma7.1+}
{
{\rm dim\,}H^1(\gl,V)\le \mbox{$\sum\limits_{i=1}^p$}{\rm dim\,}H^1(\gl,M(P_i)).
}
\end{lemma}

\begin{proof} We prove (\ref{Lemma7.1}) by induction on number of
composition factors of $V$. If $V$ is irreducible, the claim is obvious.
Suppose $V$ is not irreducible, and let $L_\nu$
be an irreducible submodule of $V$. Then the exact sequence $0\to
L_\nu\to V\to V/L_\nu\to 0$ gives \equan{Lemma7.2} {
H^1(\gl,L_\nu) {\,}^{\ i}_{\dis\to}{\,}H^1(\gl,V) {\,}^{\
j}_{\dis\to}{\,}H^1(\gl,V/L_\nu). } Thus as a vector space,
$H^1(\gl,V)\cong i(H^1(\gl,L_\nu))\oplus j(H^1(\gl,V))$, and ${\rm
dim\,}H^1(\gl,V)\le$ ${\rm dim\,}H^1(\gl,L_\nu)+{\rm
dim\,}H^1(\gl,V/L_\nu)\le\sum_{\mu\in P(V)}{\rm
dim\,}H^1(\gl,L_\mu)$.
\end{proof}

\begin{lemma}\label{lemm6.5}
Suppose $V$ is a module without any trivial $($i.e., $1$-dimensional${\sc\,})$
composition factor. Then
\equan{Lemma8}
{\mbox{$
\sum\limits_{\mu\in P_0(V)}{\rm dim\,}H^1(\gl,L_\mu)\le {\rm dim\,}H^1(\gl,V).
$}}
Furthermore, suppose we have $(\ref{Lemma7.1-})$ such that all $P_i$
correspond submodules of $V$ and $V/\oplus_{i=1}^p M(P_i)$ does not
contain a trivial composition factor, then
\equan{Lemma8++}
{\mbox{$
\sum\limits_{i=1}^p{\rm dim\,}H^1(\gl,M(P_i))\le {\rm dim\,}H^1(\gl,V).
$}}
\end{lemma}

\begin{proof} Let $V'=\oplus_{\mu\in P_0(V)}L_\mu\subset V$. Clearly,
$H^1(\gl,V')\cong\oplus_{\mu\in P_0(V)}H^1(\gl,L_\mu)$.
{}From $0\to V'\to V\to V/V'\to0$, by (\ref{long-exact}), we have
$ 
{
0=H^0(V/V')\to H^1(V') {\,}^{\ i}_{\dis\to}{\,}H^1(V),
}$
i.e., $i$ is an injective map.
\end{proof}

We need to use \cite[Conjecture 4.1]{VZ} and \cite[Main Theorem]{Br}
to determine $P^\vee(\mu)$ and $P(\mu)$ for some $\mu$.
Note that the necessary and sufficient condition for the existence
of an element $d$ in the Weyl group of $\gl_0$, which relates
the weight $\mu=(\mu_1,...,\mu_m\,|\,\mu_{\cp1},...,\mu_{\cp n})$ to
a weight $\l=(\l_1,...,\l_m\,|\,\l_{\cp1},...,\l_{\cp n})$ such that
$\l+\rho$ is dominant through the equation
$\l+\rho=d(\mu+\rho)$,
is the following equalities of sets:
\equan{5.0.3}
{
\{\l_1\!+\!m,...,\l_m\!+\!1\}=\{\mu_1\!+\!m,...,\mu_m\!+\!1\},\ \quad
\{\l_{\cp1}\!-\!1,...,\l_{\cp n}\!-\!n\}=\{\mu_{\cp1}\!-\!1,...,\mu_{\cp n}\!-\!n\},
}
where dominance of $\l+\rho$ means that
\equan{5.0.2}
{
\l_1+m\ge\l_2+m-1\ge...\ge\l_m+1,\ \ \
\l_{\cp1}-1\ge\l_{\cp2}-2\ge...\ge\l_{\cp n}-n.
}
Recall notations (\ref{notations1}) and (\ref{notations2}).
For $\tau=\mu^{(\ell)}$, we have (\ref{atypical_tau}) and (\ref{tau_data}).

\begin{lemma}\label{lemm6.6}
 {\rm(1)} Suppose $0\le\ell\le n-1$. Then
\begin{eqnarray}
\label{5.0.1}
P^\vee(\mu^{(\ell)})\!\!\!\!&=
\!\!\!\!&
\{\l_\theta(\mu^{(\ell)})=
\mu^{(1,j_1)}+\cdots+\mu^{(s-1,j_{s-1})}+\mu^{(s,\ell)}+\mu^{(s+1,j_{s+1})}+
\cdots+\mu^{(k,j_k)}
\nonumber\\&&\ \
+\,\theta_{n-\ell-1}\a_{m-n-i-\ell+3,\cp n-\ell- i}
+\theta_n\a_{m+1-i,\cp n+1- i}\,\ |\ \,1\le s\le n-\ell-1,
\nonumber\\&&\ \
s\le k\le n-1,\,j_1>j_2>...>j_{s-1}>\ell+1,\,
\ell>j_{s+1}>...>j_k,
\nonumber\\&&\ \
\theta_{n-\ell-1},\theta_n \in\{0,1\}
\},
\end{eqnarray}
$($note that
when $\theta_{n-\ell-1}=\theta_n=1$, the sum of the last two terms together with
$\mu^{(s,\ell)}$ is equal to $\mu^{(s,\ell+1)}).$

{\rm(2)} Similarly, we have
\begin{eqnarray}
\label{5.0.4}
P^\vee(-\amin)
\!\!\!\!&=\!\!\!\!&
\{
\l_\theta(-\amin)=-\amin+\mu^{(1,j_1)}+\cdots+\mu^{(k,j_k)}+\theta_1\amin
\nonumber\\&&\ \ \ \ \
\,|\,\
0\le k\le n-1,
n-1>j_1>...>j_k\ge0,\, \theta_1=0,1\}.
\end{eqnarray}
\end{lemma}

\begin{proof}
(1) Suppose
$ 
{
\{j\,|\,\theta_j=1,\,j\ne n-\ell-1,n\}
=\{j_1,...,j_{s-1},j_{s+1},...,j_k\},
}$
where $j_1>$ $...>j_{s-1}>\ell+1,\,\ell>j_{s+1}>...>j_k$.
Using (\ref{tau_data}),
by induction on $\#\{j\,|\,\theta_j=1\}$, we obtain (\ref{5.0.1}).
Similarly we have (\ref{5.0.4}).
\end{proof}
For a weight $\tau$, we denote
\equan{5Lemma12.01}
{
P_+(\tau)=\big\{\mu\in P(\tau)\,\big|\,
\mbox{$\sum\limits_{i:\,1\le i\le m,\mu_i<0}$}\mu_i\ge-1
\big\},\ \ \
P_-(\tau)=P(\tau)\bs P_+(\tau).
}

\begin{lemma}\label{lemm6.7}
{\rm(1)} Suppose $0\le\ell\le n-1$. We have
\begin{eqnarray}
\label{5.1.1}
P_+(\mu^{(\ell)})
&\!\!\!=
\!\!\!\!&
\{
\mu^{(\ell)},\,\,\mu^{(\ell)}_-,
\,\,\mu^{(\ell-1)}_+,\,\,
\mu^{(\ell-1)},\,\,
0,\,
\mu^{(\ell)}-\amin,\,\,
\nonumber\\&&\ \
\mu^{(\ell)}_--\amin,\,\,
\mu^{(\ell-1)}_+-\amin,\,\,
\mu^{(\ell-1)}-\amin,\,\,
-\amin\},
\end{eqnarray}
where if $\ell=n-1$ then the last 5 weights do not occur,  and if $\ell=0$ then
the 6 weights with plus or minus subscript or with supscript $(\ell-1)$
do not occur.

{\rm(2)} If $\mu\in P_-(\mu^{(\ell)})$, there is no link $\mu^{(\ell)}\link\mu$
unless $n\ge2,\,\ell=n-2$ and $\mu=\ETd$.

{\rm(3)} If $\mu\in P_-(-\amin)$, there is no link $-\amin\link\mu$ unless
$n\ge2,\,\mu=\ETb,\ETd$.
\end{lemma}

\begin{proof}
(1) Let $\mu\in P_+(\mu^{(\ell)})$. Then $\mu$ satisfies the
condition (\ref{mu_form}) and so \equan{5.2.0} {
0\le\mu_1\le\ell\!+\!1,\;\;\;\; 0\le\mu_2,...,\mu_{m-n+\ell}\le
1,\;\;\;\; \mu_{m-n+\ell+1}\!=\!...\!=\!\mu_{m-1}\!=\!0,\;\;\;\;
\mu_m=0,-1. } {}From this and using \cite[Conjecture 4.1]{VZ}
or \cite[Main Theorem]{Br}, it
is straightforward to show that $\mu$ must be one of weights
in (\ref{5.1.1}).

(2) Suppose $\mu\in P_-(\mu^{(\ell)})$. Say $\mu^{(\ell)}\rrar\mu$.
If  $\mu\ne\ETb,\ETd$,
then $H^2(\gl,L_\mu)=0$ (since $\mu$ cannot be a weight in (\ref{5.00-})
either). Since $\mu$ must have the form (\ref{mu_form}), one can easily check that
the lowest weight of $L_\mu$ (which is lower than $-\mu^{(\ell)}$)
is not a sum of distinct
negative odd roots, this means that $H^1(\gl,M(\mu^{(\ell)}\rrar\mu))=0$
(if $H^1(\gl,M(\mu^{(\ell)}\rrar\mu))\ne0$, then in particular
the lowest weight vector of $L_\mu$ can be
generated by a primitive vector of weight $0$). Then from the short
exact sequence $0\to L_\mu\to M(\mu^{(\ell)}\rrar\mu)\to L_{\mu^{(\ell)}}\to0$,
by (\ref{long-exact}), we obtain that $H^2(\gl,L_\mu)\ne 0$,
contradicting the fact that $H^2(\gl,L_\mu)=0$.
If $\mu=\ETb,\ETd$,
then only the second case can occur and in this case $\ell=n-2$.

The proof of  part (3) is similar.
\end{proof}

\begin{lemma}\label{add-lemma}
Suppose $n\ge2$ and let $\mu=\mu^{(n-1)}$. We have
\equa{4Case3+}
{
\L_\mu\!=\!(n\!-\!1,...,n\!-\!1,0\,|\,0,1\!-\!m,...,1\!-\!m)
=\mu^{(1,n-2)}+\mu^{(2,n-3)}+\cdots+\mu^{(n-1,0)}.
}
Consider the Kac module $V_{\L_\mu}$, $P(\L_\mu)$ has a subgraph
$($not necessarily closed$)$
\equa{4Case3.1}
{
\mu^{(n-2)}{\ }^{\dis\rrar\l_1\rrar}_{\dis\rrar\ 0\ \rrar}{\ }-\amin,
}
where $\l_1=\mu^{(n-2)}-\amin$.
\end{lemma}
\begin{proof}
First, (\ref{4Case3+}) follows from (\ref{Lambda_mu}).
Note that all weights in (\ref{4Case3.1}) are strongly primitive weights
(they in fact correspond to unlinked codes defined in \cite{SHK}).
Since $0$ is not a bottom primitive weight, there must be a weight
linked to $0$ in $U_{V_{\L_\mu}}(0)$, and
this weight has to be $-\amin$ by Lemma \ref{lemm6.1}.
Thus $0\rrar-\amin$. We see from (\ref{5.0.1}) that $\mu^{(n-2)}$
is the only primitive weight of
$V_{\L_\mu}$ other than $\mu$ with non-zero $1$-cohomology. By Lemma
\ref{lemm6.1}, we must have $\mu^{(n-2)}\rrar0$ (since there
must be some weight $\tau$ such that $\tau\erar0$).
Note that the primitive vector $ v_{\l_1}$ can be obtained from the primitive vector
$v_{\mu^{(n-2)}}$, in fact, $ v_{\l_1}=f_\amin v_{\mu^{(n-2)}}$ (indeed it is
non-zero and strongly primitive). Thus we have
$\mu^{(n-2)}\drar \l_1$.
Obviously $\l_1\drar-\amin$ since $-\amin$ is the bottom primitive weight.
If there exists some $\nu$ such that
\equan{AdA1}
{
\mu^{(n-2)}\rrar\nu\drar \l_1
\mbox{ \ \ or \ \ }\l_1\drar\nu\rrar-\amin,
}
then since $\level(\l_1)=
\level(\mu^{(n-2)})-1$, we must have $\level(\nu)< \level(\l_1)$.
By Remark \ref{rema6.1},
we must have
\equan{add111}
{\nu\in P(\L_\mu)\cap P(\mu^{(n-2)})
\cap P(\l_1).
}
By Lemma \ref{lemm6.7}, we have $\nu\in P_+(\mu^{(n-2)})$, but
we see from (\ref{5.0.4}) that
such $\nu$ does not exist. Thus $\mu^{(n-2)}\rrar\l_1\rrar-\amin$.
\end{proof}

\begin{lemma}\label{lemm6.10}
In $V_{\mu^{(j)}}$, we have the following subgraph of $P(\mu^{(j)})$:
\equa{Lemma5.1}
{
\mu^{(j-1)}{\,}
^{^{^{^{\dis\llar\mu^{(j)}_-\llar}}}}
_{_{_{_{\dis\lar\!\!\!\!-\!\!\!\!-\!\!\!\!-
\!\!\!\!-\!\!\!\!-\!\!\!\!-\!\!\!\!-\!\!\!\!-\!\!\!\!\rb{1.2pt}{$\,\ssc<$}}}}}
{}
{^{^{^{^{^{\dis\,\mu^{(j)}}
}}}}_{_{_{_{_{_{_{_{\dis\ 0}}}}}}}}
}
\!\!\!\!\!\!\!\!\!\!{\dis\downarrow}\ \ \
{}
^{^{^{^{\dis\rrar\l_1\rrar}}}}
_{_{_{_{\dis\rb{1.2pt}{$\,\ssc>$}\!\!\!-\!\!\!\!-\!\!\!\!-\!\!\!\!-
\!\!\!\!-\!\!\!\!-\!\!\!\!-\!\!\!\!-\!\!\!\!\rar}}}}
{\,}
-\amin\ ,
}
with $\l_1=\mu^{(j)}-\amin$, where if $j=0$, the part
$\mu^{(j-1)}{\sc\,}^{\llar{\ssc\,}\mu^{(j)}_-{\ssc\,}\llar}
_{\ \,\lar\!\!\!-\!\!\!-\!\!\!-\!\!\!-\!\!\!-\!\!\!-
\!\!\!\!\!\dis\rb{1.2pt}{$\,\ssc<$}}$
 is missing.
\end{lemma}
\begin{proof}
We prove by induction on $j$. First suppose $j=0$. Then (\ref{Lemma5.1}) is reduced to
\equa{Lemma5.1-}
{
\mu^{(0)}{\,}_{\dis\rrar0\,\rrar}^{\dis\rrar\l_1\rrar}-\amin
}
The part $\mu^{(0)}\rrar0$ is clear since it is the dual of the part
$0\rrar\mu^{(0)}$ in $P(\L_{\mu^{(0)}})$ which has been used to prove
$H^1(\gl,L_{\mu^{(0)}})\ne0$. Since $-\amin$ is the only other
primitive weight of $V_{\mu^{(j)}}$ with non-zero $1$-cohomology and
since $0$ is not the bottom primitive weight, $0$ must be linked to some primitive weight $\l$
with arrow pointed to $\l$. By Lemma \ref{lemm6.1}, we must have $0\rrar-\amin$.
The proof of the part $\mu^{(0)}\rrar\l_1\rrar-\amin$ is similar to that of (\ref{4Case3.1}).

Next suppose $j>0$. We observe the following facts:

{\bf Fact 1:} As in case $j=0$, we have $\mu^{(j)}\rrar0$.

{\bf Fact 2:} $\mu^{(j-1)}\llar0$:
Note that $\mu^{(j-1)},-\amin$ are the only primitive weights of $V_{\mu^{(j)}}$
beside $\mu^{(j)}$ which have non-zero $1$-cohomology.
We claim that $0$ is not a strongly primitive weight of $V_{\mu^{(j)}}$.
Otherwise, the bottom primitive weight $\tau$ of $V_{\L_\mu^{(j)}}$ is in $P(V_0)$,
but  $\tau_1=j>0$ (note from \cite[Theorem 3.2]{S} that
$\tau=\mu^{(j)}-\sum_{\a\in SW_{\mu^{(j)}}}\a$, where
$SW{\mu^{(j)}}$ is the set of roots in the southwest chains of $\mu^{(j)}$, there is
only one element of $SW{\mu^{(j)}}$ located on the first row, i.e., the $n$-th
atypical root $\amax$),
a contradiction (using \cite[Proposition 3.5]{VZ}, one can also see that $\tau$ is not in
$P(V_0){\ssc\,}$).
Thus there is a primitive weight $\l$ in $U(0)=U_{V_{\L_\mu}}(0)$
with level higher than $0$.
Thus $\l$ has the form (\ref{gl1_highest}) satisfying (\ref{3.2}).
This $\l$ must be $\mu^{(j-1)}$. Thus $\mu^{(j-1)}\dlar0 $ and $\mu^{(j-1)}$ is the
highest weight in $\ol M(\mu^{(j-1)}\dlar0)$. So we have a highest
weight module $\ol M(\mu^{(j-1)}\drar0)$, which is a quotient
of $V_{\mu^{(j-1)}}$. But by the inductive assumption,  we have $\mu^{(j-1)}\rrar0$
in $V_{\mu^{(j-1)}}$, thus also in $\ol M(\mu^{(j-1)}\drar0)$. Therefore, $\mu^{(j-1)}\llar0$
 in $\ol M(\mu^{(j-1)}\dlar0)$. Thus $\mu^{(j-1)}\llar0$.

{\bf Fact 3:} $0\rrar-\amin$: First note that
$(-\amin)^*=\mu^{(n-1)}$. As $\mu^{(j)}=\sum_{\a\in\G}\a$
for a unique subset $\G\subset\D_1^+$, there is up to scalars a unique
$\gl_0$-highest weight vector of weight $0$ which must be a
primitive vector $v_0$ of weight $0$. Similarly, there is up to scalars a unique
$\gl_0$-highest weight vector $v'$ of weight $-(\mu^{(n-1)})^R$
which must be the one corresponding to the lowest $\gl_0$-highest
weight vector of $L_{-\amin}$, and we have \equa{add.Fact3}
{\mbox{$ v'=\prod\limits_{\a\in\G_1}f_\a v_0, $}} where
$\G_1=\{\a_{1,\cp1},...,\a_{m-1,\cp1},\a_{m,\cp1},\a_{m,\cp2},...,\a_{m,\cp
n}\}$. Indeed the right hand side of (\ref{add.Fact3}) is non-zero
(since $v_0$ has the leading term $\prod_{\a\in\G}f_\a
v_{\mu^{(j)}}$) and is a $\gl_0$-highest weight vector of weight
$-(\mu^{(n-1)})^R$. This means $0\drar-\amin$. We cannot have
$\mu^{(j-1)}\drar-\amin$. Assume otherwise. Since $U(\mu^{(j-1)})$ is a
quotient of $V_{\mu^{(j-1)}}$, by inductive assumption, we must
have $\mu^{(j-1)}\rrar0\rrar-\amin$, contradicting
$\mu^{(j-1)}\llar0$. Thus in $U(0)/U(\mu^{(j-1)})$ which is now a
highest weight module, we still have $0\drar-\amin$. So
$0\rrar-\amin$, since $-\amin$ is the unique primitive weight in
$U(0)/U(\mu^{(j-1)})$ with non-zero $1$-cohomology.

{\bf Fact 4:} $\mu^{(j-1)}\llar\mu^{(j)}_-\llar\mu^{(j)}$:
Note that $\mu^{(j)}_-,\mu^{(j-1)}$ are strongly primitive weights such that
$\mu^{(j-1)}\dlar\mu^{(j)}_-\dlar\mu^{(j)}$
 (in \cite{SHK}, the primitive vector $v_{\mu^{(j-1)}}$ is
 constructed from the primitive vector $v_{\mu^{(j)}_-}$).
If $\mu^{(j-1)}\llar\mu^{(j)}_-\llar\mu^{(j)}$ is not valid,
then there is some $\nu$ such that
\equa{some-nu11}
{
\mu^{(j-1)}\llar\nu\dlar\mu^{(j)}_-\mbox{ \ or \ }
\mu^{(j)}_-\dlar\nu\llar\mu^{(j)}.
}
Such $\nu$ can only be in $P_+(\mu^{(j)})\cap P(\mu^{(j)}_-)$ by
Lemma \ref{lemm6.6} and Remark \ref{rema6.1}.
But by (\ref{5.1.1}), we see that no $\nu$ can satisfy (\ref{some-nu11}).

{\bf Fact 5:} $\mu^{(j)}\rrar\l_1\rrar-\amin$:
The proof is similar to that of (\ref{4Case3.1}), (\ref{Lemma5.1-}) and Fact 4.
\end{proof}

Now consider dual graph of (\ref{Lemma5.1}). Note that
\equan{4P.4}
{
(\mu^{(j)}_-)^*=\mu^{(j-1)}_+, \ \ \
\l^*_1=\mu^{(n-1)}+\mu^{(2,j)}.
}
In $P(\L_{\mu^{(j)}})$, we have
\equa{duallemma5}
{
\mu^{(j-1)}{\,}
_{_{_{_{_{_{\dis\rar\mu^{(j-1)}_+\rar}}}}}}
^{^{^{^{\dis
-\!\!\!\!-\!\!\!\!-\!\!\!\!-\!\!\!\!-\!\!\!\!-\!\!\!\!-\!\!\!\!-\!\!\!\!\rar}}}}
{}
{_{_{_{_{_{_{\dis\,\mu^{(j)}\,}}}}}}
^{^{^{^{^{\dis\ 0}}}}}
}
\!\!\!\!\!\!\!\!\!\!{\dis\downarrow}\ \ \
{}
_{_{_{_{_{_{\dis\lar\l^*_1\lar}}}}}}
^{^{^{^{\dis\lar
\!\!\!\!-\!\!\!\!-\!\!\!\!-\!\!\!\!-\!\!\!\!-\!\!\!\!-\!\!\!\!-\!\!\!\!-}}}}
{\,}
(-\amin)^*\ .
}
\begin{lemma}\label{lemm6.11}
Denote $M_1= M(\mu^{(j-1)}\rrar\mu^{(j-1)}_+)$, $M_2= M(\l_1^*\llar(-\amin)^*)$. We have
$H^1(\gl,M_1)= H^1(\gl,M_2)=0 $.
\end{lemma}

\begin{proof}
The arguments for the proof of $H^1(\gl,M_1)=0$ are similar to those given after (\ref{4.Add4+}).

To prove $ H^1(\gl, M_2)=0$, first we consider the Kac module $V_{\l_1}$.
Since we have a highest weight module with graph $\l_1\rrar-\amin$ by
(\ref{Lemma5.1-}), we must also have this in $V_{\l_1}$. It is
straightforward to verify that $-\amin$ is the only primitive weight of $V_{\l_1}$
with non-zero $1$-cohomology and $0$ is not a primitive weight of $V_{\l_1}$.
Thus in the dual Kac module $V^*_{\l_1}=V_{2\rho_1-(\l_1)^R}$, we have
$\l_1^*\llar(-\amin)^*$ (so $M_2$ is a submodule of $V^*_{\l_1}$)
and $(-\amin)^*$ is the only primitive weight with non-zero $1$-cohomology
while $0$ is not a primitive weight. Thus
$H^0(\gl,V^*_{\l_1}/ M_2)=0$.
{}From the exact sequence $0\to M_2\to V^*_{\l_1}\to V^*_{\l_1}/ M_2\to0$,  we obtain
\equan{4P.Add1}
{
0=H^0(\gl,V^*_{\l_1}/M_2)\to H^1(\gl,M_2)\to H^1(\gl,V^*_{\l_1})=0.
}
Thus we have the lemma.
\end{proof}

\subsection{Computation of second cohomology groups}\label{2-irr-main}
With the technical preparations in the last subsection, we can now
compute the second cohomology groups with coefficients in the finite
dimensional irreducible modules.

Suppose $H^2(\gl,L_\mu)\ne0$. Then $\mu\ne0$.
Let $\vp\in Z^2(\gl,L_\mu)^{\gl_0}$.
As before, we can assume
\equa{4.00}
{
\vp|_{\glss\times\gl}=0\mbox{ \ and \ }\vp|_{\check\rho_1\times\gl}\mbox{ is a $1$-cocycle}.
}

\subsubsection{The case with $f_\a\vp(e_\b,e_\g)=0$ for $\b,\g\in\D_1^+$
and $\vp(e_\b,e_\g)\ne0$ for some $\b,\g\in\D_1^+$}

Then $\vp(e_\b,e_\g)\in L^\bottom_\mu$. This means that
$\mu_\bottom=\ETa$ or
$\ETc$ (in the latter case $n\ge2$). We obtain (cf.~\cite[Proposition 3.5]{VZ})
\equa{4.1.1}
{
\mu=\mu^{(n-1)}_+,\ \ \ \mbox{ or \ \ }\mu=\mu^{(n-1)}_-\mbox{ with } n\ge2.
}
As in (\ref{3Case1.5}), we have $H^2(\gl,L_\mu)=\C\ol{\phi_{2,1}}$.

\subsubsection{The case with $e_\a\vp(f_\b,f_\g)=0$ for $\b,\g\in\D_1^+$
and $\vp(f_\b,f_\g)\ne0$ for some $\b,\g\in\D_1^+$}
In this case,
\equa{4.2.1}
{
\mu=\ETb=(\mu^{(n-1)}_+)^*,\ \ \ \mbox{ or \ \ }
\mu=\ETd=(\mu^{(n-1)}_-)^* \mbox{ with }n\ge2.
}
and so we have $H^2(\gl,L_\mu)\cong\C$.

\subsubsection{The case with $\mu=\mu^{(n-1)}$ or its dual $\mu=-\amin$}
It suffices to consider the case
$\mu=-\amin$.
First suppose $n=1$. Then $\L_\mu=0$ and $V_0$ has two
composition factors $L_\mu$ and $L_0=\C$. Since $H^1(\gl,L_0)=H^2(\gl,V_0)=0$,
from $0\to L_\mu\to V_0\to L_0\to0$, we have
\equan{4Case3}
{
0=H^1(\gl,L_0)\to H^2(\gl,L_\mu)\to H^2(\gl,V_0)=0,
}
Thus we obtain $H^2(\gl,L_\mu)=0$.
So suppose $n\ge2$.
Then by (\ref{Lemma7.1+}) and (\ref{4Case3.1}), we obtain
\equan{4.aaa}
{
{\rm dim\,}H^2(\gl,L_\mu)={\rm dim\,}H^1(\gl,V_{\L_\mu}/L_\mu)
\le {\rm dim\,}H^1(\gl, M(\mu^{(n-2)}\rrar\mu_+^{(n-3)})).
}
We claim that $H^1(\gl, M(\mu^{(n-2)}\rrar\mu_+^{(n-3)}))=0$.
Otherwise, there exists a non-split extension
of $ M(\mu^{(n-2)}\rrar\mu_+^{(n-3)})$:
\equan{4.Add4}
{
0\to M(\mu^{(n-2)}\rrar\mu_+^{(n-3)})\to\wh W\to\C\to0,
}
i.e., we have an indecomposable module with graph
\equa{4.Add4+}
{0\rrar\mu^{(n-2)}\rrar\mu^{(n-3)}_+\mbox{ \ \ or \ \ }
\mu^{(n-2)}\rrar\mu^{(n-3)}_+\llar0.
}
Since $H^1(\gl,L_{\mu^{(n-3)}_+})=0$, the second case cannot occur.
For the first case, we obtain that the lowest $\gl_0$-highest
weight vector $v_{\l}\in L_{\mu^{(n-3)}_+}$ is of weight
\equan{4.Add4++}
{\l=-(\mu^{(n-2)}_+)^R=(0,-1,...,-1,-n\,|\,m,1,...,1,0)
}
(which is lower than the lowest $\gl_0$-highest weight of $L_{\mu^{(n-2)}}$)
and can be generated by a primitive  vector $v_0$ of weight $0$
such that $v_\l\in U(\gl_{-1})v_0$. But $-\l$ cannot be written
as a sum of distinct positive odd roots. This is a contradiction.
Thus we have the claim. So we obtain
\equan{4.3.1}
{
{\rm dim\,}H^2(\gl,L_\mu)=0 \mbox{ if }
\mu=\mu^{(n-1)}\mbox{ or }-\amin.
}

\subsubsection{The case with $\mu=\mu^{(j)}$ with $j\le n-2$ $($thus $n\ge2)$}
From (\ref{Lemma7.1+}) we obtain
\equan{4P.Add2}
{
\begin{array}{ll}
{\rm dim\,}H^2(\gl,L_{\mu^{(j)}})
\!\!\!\!&
=
{\rm dim\,}H^1(\gl,V_{\L_{\mu^{(j)}}}/L_{\mu^{(j)}})
\vs{4pt}\\ &
\le
{\rm dim\,}H^1(\gl,M_1)+
{\rm dim\,}H^1(\gl,M_2),
\end{array}
}
where $M_1,\,M_2$ are as in Lemma \ref{lemm6.11}.
By Lemma \ref{lemm6.11}, the far right hand side of the above
inequality vanishes. Thus $${\rm dim\,}H^2(\gl,L_{\mu^{(j)}})=0.$$

Before considering the next case with $\mu=\mu_\pm^{(\ell)}$, we first
assume that $\mu$ is not any of the weights considered the earlier cases.
Then $H^1(\gl,V_\mu)=$ $H^2(\gl,V_\mu)=0$ by Theorems \ref{theo3.1} and \ref{theo4.1}.
Suppose $H^2(\gl,L_\mu)\ne0$. In (\ref{4.00}), we can further suppose
$\vp(\check\rho_1,\gl)=0.$
We shall further suppose that $\vp(e_\a,f_\b)\ne0$ for some $\a,\b$ (otherwise
it is already considered above). Thus $L_\mu$ must contain a $\gl_0$-submodule
isomorphic to an irreducible $\gl_0$-submodule
of $\gl_{+1}\wedge\gl_{-1}$, i.e., contain a $\gl_0$-highest weight
which is one of the following
\equa{5.00-}
{
(1,0,...,0,-1{\sc\,}|{\sc\,}0,...,0),\,
(0,...,0{\sc\,}|{\sc\,}1,0,...,0,-1),\,
(1,0,...,0,-1{\sc\,}|{\sc\,}1,0,...,0,-1),\,  0.
\!\!}
Thus
\equa{5.00}
{
\mbox{$\mu$
 minus a weight in (\ref{5.00-}) is a sum of distinct positive odd roots}.
}
Since $H^2(\gl,L_\mu)\ne0$,
there exists the exact sequence
\equa{5.A1}
{
0\to L_\mu\to V_1\to V_2\to \C\to0.
}
This implies that
 $V_1$ has graph $\mu\llar\tau$ for some $\tau$, and $V_2$ has graph $\tau\llar0$.
Thus
\equa{5.A2}
{
\tau=-\amin,\,\mu^{(\ell)},\,0\le \ell\le n-1\mbox{ \ \  and \ \ }\mu\in P(\tau)\cup P^\vee(\tau).
}
Therefore $\mu$ is a weight in (\ref{5.0.1})--(\ref{5.1.1}).

First we give the duals of some useful weights:
\equa{5.dual1}
{
(\mu^{(\ell)}_+)^*=
\biggl\{\begin{array}{ll}
\mu^{(n-2)}-\amin&\mbox{if \ \ }\ell=n-2,
\vs{4pt}\\
\mu^{(\ell+1)}_-&\mbox{if \ }0\le\ell\le n-3,
\end{array}
}
and $(\mu^{(n-1)}_-)^*=\ETd$ (this weight
has been considered in
(\ref{4.1.1})), and
\equa{5.dual2}
{
(\mu^{(\ell)}-\amin)^*=
\biggl\{\begin{array}{ll}
\mu^{(n-2)}_+&\mbox{if \ \ }\ell=n-2,
\vs{4pt}\\
\mu^{(n-1)}+\mu^{(2,\ell)}&\mbox{if \ }0\le\ell\le n-3,
\end{array}
}
\equa{5.dual3}
{
(\mu^{(\ell)}\!+\!\mu^{(2,j)})^*\!=\!\biggl\{\begin{array}{ll}
\mu^{(j)}-\amin\!\!\!&\mbox{if }\ell=n-1,\,j\le n-3,
\vs{4pt}\\
\mu^{(\ell)}\!+\!\mu^{(2,j)}\!\!\!&\mbox{if }0\!\le\!\ell\!\le\! n\!-\!2,\,
0\!\le\! j\!\le\!\ell\!-\!1,\,
(\ell,j)\!\ne\!(n\!-\!2,n\!-\!3),
\end{array}
}
and
\equan{5.daul+}
{
\begin{array}{l}
(\mu^{(n-1)}+\mu^{(2,n-2)})^*=(0,...,0,-2,-2\,|\,2,2,0,...,0),
\vs{4pt}\\
(\mu^{(n-2)}+\mu^{(n-3)})^*=\mu^{(n-3)}_++\ETd,
\end{array}
}
(these two weights do
not meet (\ref{5.00})).

\subsubsection{The case with $\mu=\mu_\pm^{(\ell)}$}
\begin{lemma}\label{lemm6.12}
If $\mu=\mu^{(\ell)}_+,\,0\le\ell\le n-2$, then ${\rm dim\,}H^2(\gl,L_\mu)=1$.
\end{lemma}
\begin{proof}
By (\ref{5.0.1}) and (\ref{5.0.4}), we see that $\mu^{(\ell)}$ is the only primitive
weight of $V_{\mu_+^{(\ell)}}$ with non-zero $1$-cohomology.
In $V_{\mu_+^{(\ell)}}$, we
have $\mu^{(\ell)}_+\rrar\mu^{(\ell)}$ since
by (\ref{duallemma5}) and by Remark \ref{rema6.1},
we see that a highest weight module with graph $\mu^{(\ell)}_+\rrar\mu^{(\ell)}$
exists (thus in $V_{\mu_+^{(\ell)}}$, it must be so).
By considering the dual Kac module $(V_{\mu_+^{(\ell)}})^*$, it is straightforward
to see that ${\rm dim\,}H^1(\gl,(V_{\mu_+^{(\ell)}})^*/L_{\mu_+^{(\ell)}}^*)=1$.
Thus ${\rm dim\,}H^2(\gl,L_{\mu_+^{(\ell)}})={\rm dim\,}H^2(\gl,L^*_{\mu_+^{(\ell)}})=
1$.
\end{proof}

\subsubsection{The case with $\mu=\mu^{(\ell)}+\mu^{(2,\ell-1)},\,
1\le\ell\le n-3$}

\begin{lemma}\label{lemm6.13}
If $\mu=\mu^{(\ell)}+\mu^{(2,\ell-1)},\,
1\le\ell\le n-3$, then $H^2(\gl,L_\mu)=0$.
\end{lemma}
\begin{proof}
First we prove that there is a subgraph of $P(\mu)$,
\equa{5.Lemma14.1}
{
\mu^{(\ell-2)}{\,}
^{^{^{^{\dis\llar\l_1\llar}}}}
_{_{_{_{\dis\lar\!\!\!\!-\!\!\!\!-\!\!\!\!-\!\!\!\!-\!\!\!\!-\!\!\!\!-
\!\!\!\!-\!\!\!\!-\!\!\!\!\rb{1.2pt}{$\,\ssc<$}}}}}
{}
{^{^{^{^{^{
^{^{\dis\mu}_{\dis\downarrow}}
_{\dis\,\mu^{(\ell)}
}}
}}}}_{_{_{_{_{_{_{_{\dis\ 0}}}}}}}}
}
\!\!\!\!\!\!\!\!\!\!{\dis\downarrow}\ \ \
{}
^{^{^{^{\dis\rrar\l_2\rrar}}}}
_{_{_{_{\dis\rb{1.2pt}{$\,\ssc>$}\!\!\!-\!\!\!\!-\!\!\!\!-
\!\!\!\!-\!\!\!\!-\!\!\!\!-\!\!\!\!-\!\!\!\!-\!\!\!\!\rar}}}}
{\,}
-\amin\ ,
}
where $\l_1=\mu^{(\ell)}+\mu^{(2,\ell-2)},\,\l_2=\mu^{(\ell)}-\amin$,
and where if $\ell=1$, the part $\mu^{(\ell-2)}{}^{\llar\l_1\llar}_{\lar\!\!\!-
\!\!\!-\!\!\!-\!\!\!-\!\!\!-\!\!\!-
\!\!\!\!\!\dis\rb{1.2pt}{$\,\ssc<$}}$
 is missing.

We assume that $\ell\ge2$ as the case $\ell=1$ can be regarded as a special
case (cf.~the proof of Lemma \ref{lemm6.10}).
We observe the following Facts:

{\bf Fact 0$'$:} One can check that all weights in (\ref{5.Lemma14.1}) are primitive
weights of $P(\mu)$ and all primitive weights of $P(\mu)$ with non-zero $1$-cohomology appear
in (\ref{5.Lemma14.1}).

{\bf Fact 1$'$:} $\mu^{(\ell)}\drar\l_1$: Since
$\mu^{(\ell)}$ is not a strongly primitive weight (as in Fact 2 of Lemma
\ref{lemm6.10}),   $\mu^{(\ell)}
\drar\l$ for some $\l\in P(\mu)\cap P^\vee(\mu^{(\ell)})$ by Remark \ref{rema6.1}.
By using (\ref{5.0.1}), one can check that $ P(\mu)\cap P^\vee(\mu^{(\ell)})=\{\l_1\}$.
(From what to be discussed below, we can see that we must have $\mu^{(\ell)}\rrar\l_1$,
 since if  $\mu^{(\ell)}\rrar\nu\drar\l_1$ for some $\nu$, then
$\nu\in P_+(\mu^{(\ell)})\cap P(\l_1)\cap P(\mu)$,  but there is no such $\nu$.)

{\bf Fact 2$'$:} We do not have $\l_1\drar0$: Otherwise, $\l_1$
is a strongly primitive weight (corresponding to an unlinked code),
thus $U(\l_1)$ is a quotient of $V_{\l_1}$. From the dual Kac module $V^*_{\l_1}$, one
can see (using the same arguments as that given
in the paragraph after (\ref{define-Lambda}))
that in $V_{\l_1}$, we have $\l_1\drar\mu^{(\ell)}\drar0$, and so we must also
have this in $U(\l_1)$. This contradicts the fact that $\mu^{(\ell)}\drar\l_1$.

{\bf Fact 3$'$:} $\mu^{(\ell)}\rrar0$: First we have $\mu^{(\ell)}\drar0$
(as the proof before (\ref{No1.})).
Then consider $U(\mu^{(\ell)})/U(\l_1)$ which is now a quotient of $V_{\mu^{(\ell)}}$
(since there is no primitive weight in $U(\mu^{(\ell)})/U(\l_1)$ with level higher
than that of $\mu^{(\ell)}$ by Fact 1$'$),
we see that $\mu^{(\ell)}\rrar0$ (since in $V_{\mu^{(\ell)}}$ we have this).

{\bf Fact 4$'$:} $0\erar\mu^{(\ell-2)}$: This follows from the same arguments in Fact 2
of Lemma \ref{lemm6.10}.

{\bf Fact 5$'$:} $\l_1\rar\mu^{(\ell-2)}$:
We must have $\l_1\drar\mu^{(\ell-2)}$,
otherwise by Fact 1$'$, we would have $\mu^{(\ell-2)}\in P(U(\mu^{(\ell)})/U(\l_1))
\subset P(\mu^{(\ell)})$, but $P(\mu^{(\ell)})$ does not contain
$\mu^{(\ell-2)}$ by (\ref{5.0.1}). If not $\l_1\rrar\mu^{(\ell-2)}$, then
$\l_1\drar\l\rrar\mu^{(\ell-2)}$ for some $\l$ of $P(\mu)\cap P(\l_1)$ (since
$\l_1$ is strongly primitive). By Lemma \ref{lemm6.6} (if $\l\frar\mu^{(\ell-2)}$)
or Lemma \ref{lemm6.7}(1) and (2) (if $\l\erar\mu^{(\ell-2)}$), we see that such $\l$ does not
exist.

{\bf Fact 6$'$:} As in the proof of (\ref{Lemma5.1}), we do not have
$\mu^{(\ell-2)}\drar-\amin$, but $0\rrar-\amin$ is valid.
Consider the highest weight module $U(\mu^{(\ell)})/U(\l_1)$,
we have the part $\mu^{(\ell)}\rrar0\rrar-\amin$ of the graph, and by (\ref{Lemma5.1}),
we must also have the part $\mu^{(\ell)}\rrar\l_2\rrar-\amin$.

{\bf Fact 7$'$:}
We have proved (\ref{5.Lemma14.1}) except the part $\mu\rrar\mu^{(\ell)}$. But
this must be valid. Otherwise, there is another primitive weight $\l$ with
$\l\rrar\mu^{(\ell)}$, thus we cannot have $\l\frar\mu^{(\ell)}$ since $\l_1$ is the
only weight of higher level linked to $\mu^{(\ell)}$ (Fact 1$'$). We cannot have
$\l\erar\mu^{(\ell)}$ either, since by Lemmas \ref{lemm6.6} and \ref{lemm6.7},
there is no other weight of lower level linked to $\mu^{(\ell)}$.
Thus we have proved (\ref{5.Lemma14.1}).

Now consider the dual Kac module $V^*_\mu$. By Lemma \ref{lemm6.4}, we
have
\equan{5.13.0}
{
\begin{array}{ll}
{\rm dim\,}H^2(\gl,L_\mu)
\!\!\!\!&
={\rm dim\,}H^2(\gl,L_\mu^*)
={\rm dim\,}H^1(\gl,V^*_\mu/L^*_\mu)
\vs{4pt}\\&
\le{\rm dim\,}H^1(\gl,M(\mu^{(\ell-2)}\rrar\l_1^*))
\vs{4pt}\\&\ \ \ \
+\,
{\rm dim\,}H^1(\gl,M(0\rrar\mu^{(\ell)}))
+{\rm dim\,}H^1(\gl,M((-\amin)^*\rrar\l_2^*))
\vs{4pt}\\&
=0,
\end{array}
}
where the last equality follows from the same arguments in the proof of Lemma
\ref{lemm6.11}.
\hspace*{1pt}\vs{-8pt}\end{proof}

\subsubsection{The case with $\mu=\mu^{(\ell)}+\mu^{(2,j)},\,0\le j\le\ell-2\le
n-3$}
\begin{lemma}\label{lemm6.14}
If $\mu=\mu^{(\ell)}+\mu^{(2,j)},\,0\le j\le\ell-2\le n-3,$
then ${\rm dim\,}H^2(\gl,L_\mu)=1$.
\end{lemma}
\begin{proof}
We prove by induction on $\ell+j$. We shall suppose $j\ge1$ and $\ell\le n-2$
since the case $j=0$ or $\ell=n-1$
will become a special case in the following discussion.
We need to prove the existence of the following subgraphs of $P(\mu)$,
\equa{5.15.0}
{
\begin{array}{c}
\l_3\rrar-\amin\llar0\llar\mu^{(\ell)}\llar\mu\rrar\mu^{(j)},
\vs{8pt}\\
\l_1\rrar\mu^{(j-1)},\ \l_2\rrar\mu^{(\ell-1)},
\end{array}
}
where
 $\l_1=\mu^{(\ell)}+\mu^{(2,j-1)},\,
\l_2=\mu^{(\ell-1)}+\mu^{(2,j)},\,\l_3=\mu^{(\ell)}-\amin$ and
$\l_4=\mu^{(\ell-1)}+\mu^{(2,j-1)}$.
To do this, we need to observe the following facts:

{\bf Fact 0$''$:} All
weights appeared in the graphs are primitive weights of $V_\mu$ and
all primitive weights with non-zero $1$-cohomology appear in the above graphs.

{\bf Fact 1$''$:} (a) Since $\mu^{(\ell)}$ cannot be strongly
primitive, we have $\mu^{(\ell)}\drar\l_1$ as in Fact 1$'$ of
Lemma \ref{lemm6.13}, and $P(\mu)\cap
P^\vee(\mu^{(\ell)})=\{\l_1,\mu,\mu^{(\ell)}\}$. (b) By
considering the highest weight modules $U(\l_1)$ and
$U(\mu^{(\ell)})/U(\l_1)$, we see that every primitive weight
derived from $\mu^{(\ell)}$ is in $P(\mu^{(\ell)})\cup P(\l_1)$.

{\bf Fact 2$''$:} We have $\mu^{(\ell)}\drar0\drar-\amin$ as in
the proof of Fact 3 and Fact 6$'$ of previous lemmas. Thus we must
have
$\mu^{(\ell)}{\,}^{\dis\rrar\l_3\rrar}_{\dis\rrar0\rrar}-\amin$.


{\bf Fact 3$''$:} Similar to Fact 1$''$(a), we have $\mu^{(j)}\drar\l_2$
since $P(\mu)\cap P^\vee(\mu^{(j)})=$\linebreak[4] $\{\l_2,\mu,\mu^{(j)}\}.$

{\bf Fact 4$''$:} (a) We have $\mu\rrar\mu^{(\ell)}$. Otherwise suppose
$\mu\drar\nu\rrar\mu^{(\ell)}$, then by Fact 1$''$,
$\nu\in P(\mu)\cap P_+(\mu^{(\ell)})$. Take
\equan{Fact5''}
{
\L=\mu^{(\ell)}+\mu^{(2,j)}+\mu^{(3,j-1)}+\cdots+\mu^{(j+1,0)}.
}
One can check that $\mu$ and $\mu^{(\ell)}$ are strongly primitive weight
of $P(\L)$, and $\mu\drar\mu^{(\ell)}$ in $P(\L)$. In fact in $V_\L$,
the $\gl_0$-highest weight for each of them is unique, and one can easily
construct the corresponding $\gl_0$-highest weight vector of $L_{\mu^{(\ell)}}$
from that of $L_\mu$ (cf.~arguments in the paragraph of (\ref{No1.})).
Define a homomorphism $\phi$ from $V_\mu\to V_\L$ by sending $v_\mu$ to the
corresponding primitive vector in $V_\L$. Note that any element of
$P_+(\mu^{(\ell)})$ is either a strongly primitive weight in $P(\L)$ or not
a primitive weight in $P(\L)$. Thus in $P(\L)$ we do not have $\nu\erar\mu^{(\ell)}$.
So in $P(\mu)$ we do not have this either. Thus $\mu\rrar\mu^{(\ell)}$.

(b) We have $\mu\rrar\mu^{(j)}$. Otherwise,
there exists some $\nu\in P_+(\mu^{(j)})$ such that $\mu\drar\nu\rrar\mu^{(j)}$.
Note that none of the elements of $P_+(\mu^{(j)})\bs\{\mu^{(j)}\}$ is a primitive
weight of $P(\L)$ (thus the homomorphism $\phi$ defined in (a) maps a primitive vector
$v_\nu$ in $V_\mu$ of weight $\nu$
to zero). The map $\phi$ sends a primitive vector $v_{\l_2}$ in $V_\mu$ of weight $\l_2$
to the corresponding primitive vector $v'_{\l_2}$ of the same weight $\l_2$ in
$V_\L$. This can be seen from the following arguments:
Note that $\l_2$ is a strongly primitive
weight in $P(\L)$ such that its primitive vector can be
constructed from that of $\mu$ in $V_\L$. This is because $\l_2$
is 2 levels lower than $\mu$,  and in $V_\L$ the atypical roots of $\L$
corresponding to $\mu$
(i.e., atypical roots corresponding to non-zero columns of the code
for $\L$ which defines $\mu$,
cf.~\cite{HKV, SHK})
are all disconnected from (in sense of \cite{VZ})
or not $c$-related to (cf.~Definition \ref{nqc-relation}) the other 2 atypical
roots  corresponding to $\l_2$.
Therefore, by Fact 3$''$,
$\phi$ must map a primitive vector $v_{\mu^{(j)}}$ in $V_\mu$ of weight $\mu^{(j)}$ to
a primitive vector $v'_{\mu^{(j)}}$ in $V_\L$ which can generate $v'_{\l_2}$.
In particular,
$\phi(v_{\mu^{(j)}})\ne0$. But $\nu\rrar\mu^{(j)}$ means that $\phi(v_{\mu^{(j)}})$
is generated by $\phi(v_\nu)=0$. This contradiction shows that
 we cannot have $\mu\drar\nu\rrar\mu^{(j)}$.

{\bf Fact 5$''$:} We have $\l_2\drar\l_4$, since both are strongly
primitive weights and $\l_4$ is two levels lower than $\l_2$, one
can construct $v_{\l_4}$ from $v_{\l_2}$.

{\bf Fact 6$''$:} Applying Facts 5$''$, 1$''$(a) and
4$''$(a) to the case $\ell-1$ shows that
$\l_2\rrar\mu^{(\ell-1)}\drar\l_4$. Similarly, $\l_1\drar\l_4$,
from this and Facts 3$''$ and 4$''$(b), we have
$\l_1\rrar\mu^{(j-1)}$. This proves (\ref{5.15.0}).

Now in the dual Kac module $V_{\L_\mu}$, take quotient
$V_{\L_\mu}/L_\mu$  and denote $P_1=$\linebreak
$M(0\rrar\mu^{(\ell)})),\,P_2=L_{\mu^{(j)}},\,
P_3=M(\mu^{(j-1)}\rrar\l^*_1)),\,
P_4=M(\mu^{(\ell-1)}\rrar\l^*_2)),\,P_5=$\linebreak[4]
$M((-\amin)^*\rrar\l_3^*))$, from
$H^2(\gl,L_\mu)=H^1(\gl,V_{\L_\mu}/L_\mu)$, by Lemmas
\ref{lemm6.4} and \ref{lemm6.5}, we obtain \equan{15.proof}
{\mbox{$ 1 ={\rm dim\,}H^1(\gl,P_2)\le{\rm dim\,}H^2(\gl,L_\mu)\le
\sum\limits_{i=1}^5{\rm dim\,}H^1(\gl,P_i)=1, $}} where the last
equality follows from the same arguments in the proof of Lemma
\ref{lemm6.11}. \hspace*{1pt}
\end{proof}\vspace*{-10pt}

\subsubsection{The remaining case}
For a weight $\mu$, we denote by
$\KK(\mu)$ and $\LL(\mu)$ the multiplicities of $\gl_0$-highest weight $\ETa$
in $V_\mu$ and in $L_\mu$ respectively.
We need to compute $\KK(\mu)$ and $\LL(\mu)$ for some weights $\mu$. To do this, we use
a {composite Young Diagram} to denote a weight, e.g., if
$\l=(\l^1\,|\,-\l^2)=(4,2,2,1,1\,|\,-1,-1,-3,-5)$, then the composite
Young Diagram of $\l$ is
\\\cl{\small
$
\begin{array}{rl}
(\l^2)^* \ \ \ =\ \ \
\begin{array}{r}
\Box
\vs{-6pt}\\
\Box
\vs{-6pt}\\
\Box\!\Box
\vs{-6pt}\\
\Box\!\Box
\vs{-6pt}\\
\Box\!\Box\!\Box\!\Box\!\Box
\end{array}
\\\!\!\!\!&\!\!\!\!\!\!\!\!
\begin{array}{l}
\Box\!\Box\!\Box\!\Box
\vs{-6pt}\\
\Box\!\Box
\vs{-6pt}\\
\Box\!\Box
\vs{-6pt}\\
\Box
\vs{-6pt}\\
\Box
\end{array}
\ \ \ =\ \ \ \l^1\ .
\end{array}
$}\\
In general, for a weight
$\l=(\l^1\,|\,-\l^2)=(\l_1,...,\l_m\,|\,-\l_{\cp1},...,-\l_{\cp n})$
(in the rest of this section, we use this notation to denote
a weight which is a little different from (\ref{weight1})), all
$\l_i\ge0$,
we place the {Young Diagram} of $\l^1$ on the right side
which consists of boxes in $m$ rows such that the $i$-th
row has $\l_i$'s boxes for $i=1,...,m$,
and we place the {opposite Young Diagram} $(\l^2)^*$
on the top which consists of boxes in $n$ columns such that the $\nu$-th
column has $\l_{\cp \nu}$'s boxes for $\nu=1,...,n$.

Since $V_\mu$ is completely reducible as a $\gl_0$-module, we have
\equan{5.Kac1}
{
\begin{array}{ll}
\KK(\l)
\!\!\!\!&
={\rm dim\,}{\rm Hom}_{\gl_0}(L_{\ETa}^{(0)},U(\gl_{-1})\otimes L_\l^{(0)})
={\rm dim\,}{\rm Hom}_{\gl_0}(U(\gl_{+1})\otimes L_{\ETa}^{(0)},L_\l^{(0)})
\vs{4pt}\\&
=\mbox{ the multiplicity of $\gl_0$-highest weight $\l$ in
$U(\gl_{+1})\otimes L_{\ETa}^{(0)}$}.
\end{array}
}
Note that a $\gl_0$-highest weight $\l$ of $U(\gl_{+1})\otimes L_{\ETa}^{(0)}$
can always be obtained as the sum of a $\gl_0$-highest weight $\mu$ of
$U(\gl_{+1})$ and a weight $\nu$ of $L_{\ETa}^{(0)}$ under the following rule:
Since a weight $\nu$ of $L_{\ETa}^{(0)}$ has the form
\equan{5.2-weight}
{
\nu=(0,...,0,\stl{i}{1},0,...,0,\stl{j}{1},0,...,0\,|\,0,...,
\stl{p}{-1},0,...,0,\stl{q}{-1},0,...,0)
}
(where if $i=j$ the $i$-th coordinate becomes $2$ and if $p=q$
the $p$-th coordinate becomes $-2$), when adding
$\nu$ to $\mu$, we require that $\mu_i\ne\mu_j$ and $\mu_p\ne\mu_q$ and
that the resulting weight $\mu+\nu$ be dominant.
This amounts to adding $2$ boxes to the Young Diagram $\mu^1$
such  that they are not placed in the same column,
and adding $2$ boxes to the opposite Young Diagram $(\mu^2)^*$
such that they are not placed in the same row. Equivalently, $\mu$ is obtained
from $\l$ by removing $2$ boxes of $\l^1$ from different columns and removing
$2$ boxes of $(\l^2)^*$ from different rows.

Since we only need to calculate $\KK(\l)$ for $\l$ being a
$\gl_0$-highest weight
of $U(\gl_{+1})$, $\l$ (and also $\mu$) must has the form (\ref{gl1_highest})
satisfying (\ref{3.2}), thus we can only remove $2$ boxes from different rows
 and different columns of $\l^1$ (and removing the corresponding boxes from
$(\l^2)^*\,$). Thus for $\l$ of the form
(\ref{gl1_highest}) satisfying (\ref{3.2}),
\begin{eqnarray}\label{5.2-w1}
\nonumber
\KK(\l)&=&
\mbox{the number of ways to remove $2$ boxes of $\l^1$ from}
\\\nonumber&&
\mbox{different rows and different columns such that}
\\&&
\mbox{the remaining Young Diagram is still standard}.
\end{eqnarray}
{}From this, it is straightforward to compute
\equa{5.Kac-m1}
{
\KK(\mu^{(\ell)})=\biggl\{\begin{array}{ll}
0&\mbox{if \ }\ell=0,
\vs{4pt}\\
1&\mbox{if \ }1\le \ell\le n-1.
\end{array}
}
Note from (\ref{5.1.1}) that the only possible primitive weights $\l$
of $\mu^{(\ell)}$ with $\LL(\l)\ne0$ are the first $3$ weights in (\ref{5.1.1}).
This and (\ref{5.Kac-m1}) already provide sufficient information to
obtain by induction on $\ell$,
\equa{5.Sim-1}
{
\LL(\mu^{(\ell)}_+)=1,\ \ \
\LL(\mu^{(\ell)})=\LL(\mu_-^{(\ell)})=0.
}
Similarly, we can compute
\equan{5.Kac-m2}
{
\KK(\mu^{(\ell)}+\mu^{(2,j)})=\left\{\begin{array}{ll}
0&\mbox{if \ }j=\ell-1=0,
\vs{4pt}\\
1&\mbox{if \ }1\le j=\ell-1\le n-2.
\vs{4pt}\\
3&\mbox{if \ }j=0,\,2\le \ell\le n-1,
\vs{4pt}\\
6&\mbox{if \ } 1\le j\le\ell-2\le n-3.
\end{array}\right.
}

\begin{lemma}\label{lemm6.15}
Let $\mu=\mu^{(\ell)}+\mu^{(2,j)},\,0\le j\le\ell-2\le n-3$, we have
$ 
{
\LL(\mu)=1.
}$
\end{lemma}
\begin{proof}
By (\ref{5.15.0}), we have a module $V$ with graph $\mu^{(j)}\rrar\mu$.
Let $v'_\amax$ be the $\gl_0$-highest weight vector of weight $\amax$
in $L_{\mu^{(j)}}$ (cf.~Remark \ref{rema5.2}),
and let $v_\mu$ be the highest weight vector in $L_\mu$.
Then we must have $v_\mu=\prod_{\a\in\G}e_\a v'_\amax$ for some subset
$\G\subset\D_1^+$. If $\LL(\mu)=0$, then $e_\amax v'_\amax=0$, so
$\amax\notin\G$. Thus $\G$
is a unique subset such that $\mu=\sum_{\a\in\G\cup\{\amax\}}\a$
(note that the way of writing $\mu$ as a sum of distinct positive odd roots is unique
by observing from (\ref{notations1}) that each $\mu^{(i,j)}$ corresponds to
a unique subset $S_{ij}$ of $\D_1^+$ such that $\mu^{(i,j)}=\sum_{\a\in S_{ij}}\a$).
Then we obtain a non-zero vector
$ 
{ v'_{\mu^{(\ell)}}=\prod_{\a\in\G'}e_\a v'_\amax\in V, }$ where
$\G'=\{\a_{1,\cp n-\ell},...,\a_{1,\cp n-1},\a_{2,\cp
n},..., \a_{m-n+\ell+1,\cp n}\}$,
that is  a $\gl_0$-highest weight
vector of weight $\mu^{(\ell)}$.  The weight $\mu^{(\ell)}$
cannot be a $\fg_0$-highest weight of $L_{\mu^{(j)}}$ since its level is higher than that of
$\mu^{(j)}$ (as $\ell>j$), it cannot be a $\fg_0$-highest weight of $L_\mu$ either since $L_\mu$ does not
contain a $\gl_0$-highest weight $\mu^{(\ell)}$ (it is straightforward to
see that in the Kac module
$V_\mu$, there is up to scalars a unique $\gl_0$-highest weight
vector of weight $\mu^{(\ell)}$, which is in $L_{\mu^{(\ell)}}$). This is
a contradiction against the fact that $V$ has graph
$\mu^{(j)}\rrar\mu$.
\end{proof}
\begin{lemma}\label{lemm6.16}
Let \equa{L17.0-} {\mu=\mu^{(1,j_1)}+\cdots+\mu^{(k,j_k)},\,n-1\ge
j_1>j_2>...>j_k\ge0. } Then for any $\l\in P(\mu)$, $\LL(\l)=1$ if and only if
 \equa{L17.0} {
\l=\mu^{(i)}_+, \quad \mu^{(\ell)}+\mu^{(2,j)},
\, \quad 0\le i\le n-1,\,\, 0\le j\le\ell-2\le n-3. }
\end{lemma}
\begin{proof}
{\it Case (1):} First suppose $j_p\ge j_{p+1}+3$ for $p=1,...,
k-1$ and $j_k\ge1$. It is straightforward to calculate
\equa{L17.1} { \KK(\mu)=(2k-1)+(2k-2)+\cdots+1=k(2k-1). } We can
check that the following $k(2k-1)$ weights are primitive weights
of $V_\mu$,
\begin{eqnarray}\label{L17.2}
\nonumber
&&\mu^{(j_p-1)}_+,\ \
\mu^{(j_p)}+\mu^{(2,j_q)},\ \ \mu^{(j_p)}+\mu^{(2,j_q-1)},
\\&&
\mu^{(j_p-1)}+\mu^{(2,j_q)},\ \
\mu^{(j_p-1)}+\mu^{(2,j_q-1)},
\ \ \ \,1\le p,q\le k,\ \,p<q,
\end{eqnarray}
which have non-zero multiplicities of $\gl_0$-highest
weight $\ETa$ by Lemma \ref{lemm6.15} and (\ref{5.Sim-1}).

{\it Case 2:} Suppose there are $t$ number of $j_i$'s such that
$j_i =j_{i+1}+2$, say the indices are $r_1,...,r_t$. Then we still
have (\ref{L17.1}). However, the $t$ weights
$\mu^{(j_{r_p}-1)}+\mu^{(2,j_{r_p+1})}$, $p=1,...,t$, in the list (\ref{L17.2})
are no longer primitive, but on the other hand, we have $t$ other primitive
weights $\mu^{(j_{r_p}-2)}_+,\,p=1,...,t,$ which should be added to the
list.

{\it Case 3:} Suppose further there are $t'$ number of $j_i$'s
such that $j_i =j_{i+1}+1$, say the indices are
$r'_1,...,r'_{t'}$. Then when using (\ref{5.2-w1}) to compute
$\KK(\mu)$, the removals from rows or columns in
$\mu^{(r'_p,j_{r'_p})},\,p=1,...,t'$, should not be counted.
Thus $\KK(\mu)=(k-t')(2(k-t')-1)$. Similarly, all weights in
(\ref{L17.2}) with indices $j_{r'_p},\,p=1,...,t'$, should be
removed from the list, thus the total number is still $\KK(\mu)$.

{\it Case 4:} Suppose $j_k=0$. Then $\KK(\mu)=\sum_{i=1}^{2k-2}i
=(k-1)(2k-1)$. In this case, there are $1+2(k-1)$ weights
$\mu^{(j_k-1)}_+,\, \mu^{(j_p)}+\mu^{(2,j_k-1)},
\mu^{(j_p-1)}+\mu^{(2,j_k-1)}$ in (\ref{L17.2}) which should not
in the list, the total number is again $\KK(\mu)$.
\end{proof}

Now we prove the main result on the second cohomology groups with
coefficients in finite dimensional irreducible modules.
\begin{theorem}\label{theo6}
Let $L_\mu$ be the finite-dimensional irreducible $\gl$-module with
highest weight $\mu$. Then ${\rm dim\,}H^2(\gl,L_\mu)\le1$. Furthermore,
$H^2(\gl,L_\mu)\ne0$ if and only if $\mu$ is one of the following
weights \equa{Prop4} { \mu^{(\ell)}_+,\,\ \ \mu^{(j)}_-,\,\ \
\mu^{(p)}+\mu^{(2,q)},\,\ \ \mu^{(j-1)}-\amin,\,\ \ \ETb,\,\
\ \ETd, } where  $0\le\ell\le n-1,\ \  1\le j\le n-1,\ \ 0\le p\le
q-2\le n-3,$  and in the case $n=1$, $\ETd$ will not appear.
The total number of weights in $(\ref{Prop4})$ is $\frac{(n+1)(n+2)}{2}-\d_{n,1}$.
\end{theorem}
\begin{proof}
Let $\mu$ be a weight in (\ref{5.0.1})--(\ref{5.1.1})
(recall the statement after (\ref{5.A2}) that in order for an irreducible module
to have non-trivial second cohomology, its highest weight must be in
(\ref{5.0.1})--(\ref{5.1.1})), but we assume that
$\mu$ is not a weight in (\ref{Prop4}), nor a weight already considered in
the previous cases.
Being a weight in  (\ref{5.0.1})--(\ref{5.1.1}),
$\mu$ is a primitive weight in the Kac module with a highest weight
of the form (\ref{L17.0-}), thus $\LL(\mu)=0$.   Furthermore, $\mu$ has
the form \equa{Prop4.1} {\mbox{$
\mu=\sum\limits_{i=1}^k\mu^{(i,j_i)}+\theta_1\a-\theta_2\amin, $}}
where $\theta_1,\theta_2\in\{0,1\}$ such that at most one of
$\theta_1,\,\theta_2$ is non-zero
(cf.~the statement after (\ref{5.0.1})),
and $\a$ is some atypical root. Moreover $k\ge2$, and $k\ge3$ if
$\theta_1=\theta_2=0$.

We shall show that $L_\mu$ has trivial second cohomology
by contradiction. In order for $L_\mu$  to have
non-zero second cohomology, we either have
$\mu^{(j)}\rrar\mu$ for $j\le n-1$, or $-\amin\rrar\mu$.
We only need to consider the former case, as the latter
case is the dual situation for $j=n-1$.

Assume $\mu^{(j)}\rrar\mu$.
Then we have the following module: $V:=M(\mu^{(j)}\rrar\mu)$.
By (\ref{5.A2}), $\mu\in P^\vee(\mu^{(j)})\cup P(\mu^{(j)})$.
Since at least one coordinate of $\mu^{(j)}-\mu$ is $\le0$, we have
$\mu\in P^\vee(\mu^{(j)})$. So $j\le j_1$ and $\theta_2=0$ by (\ref{5.0.1}).
%
Suppose $v_\mu=\prod_{\b\in\G}e_\b v'_\amax$ for some subset $\G\subset\D_1^+$
such that $\amax\notin\G$, where $v_\mu,\,v'_\amax$ are as in Lemma \ref{lemm6.15}.
This means that
$\mu=\sum_{\b\in\G\cup\{\amax\}}\b$ is a sum of distinct positive odd roots.
In this case,  $\G$ must contain a subset $\G'$ such that $\l:=\mu^{(j_1)}+\mu^{(2,j_2)}
=\sum_{\b\in\G'\cup\{\amax\}}\b$, and $v'_\l=\prod_{\b\in\G'}e_\b v'_\amax\ne0$ is
a $\fg_0$-highest weight vector of weight $\l$ in $V$.
However, the weight $\l$ cannot be a $\gl_0$-highest weight of
$L_{\mu^{(j)}}$ (since $\mu^{(j)}-\l$ has negative coordinates) or
$L_\mu$ (since the Kac module $V_\mu$ has up to scalars a unique $\fg_0$-highest
weight vector of weight $\l$ which appears in the composition factor $L_\l$ of $V_\mu$),
thus a contradiction results.

Therefore,  there does not exist any weight $\nu$ with non-zero
$1$-cohomology such that $\nu\rrar\mu$.
This in turn implies that $H^2(\gl,L_\mu)=0.$
\end{proof}

\begin{remark}
The dual weight $\mu^*$ of any weight $\mu$ in $(\ref{Prop4})$
is given by equations $(\ref{4.2.1})$, $(\ref{5.dual1})$--$(\ref{5.dual3})$. From these
equations we can see that $\mu^*$ also belongs to the list
$(\ref{Prop4})$.
\end{remark}

\begin{remark}
Suppose $\mu$ is a weight such that
all the central elements of $U(\gl)$ contained in
$\gl U(\gl)$ act  trivially on $L_\mu$.
Then as in Remark \ref{rema5.2}, by Theorem \ref{theo6}, one can
prove that $H^2(\gl,L_\mu)\ne0$ if and only if $L_\mu$ contains a
copy of irreducible $\gl_0$-module
$L_{\ETa},\,L_{\ETb},\,L_{\ETc}$ or $L_{\ETd}$,
that is, $L_\mu$ contains a copy of irreducible
$\gl_0$-submodule of $\gl_{+1}\wedge\gl_{+1}$ or
$\gl_{-1}\wedge\gl_{-1}$.
\end{remark}
%
\section{Cohomology groups with coefficients in enveloping algebra} %
%
\label{enveloping}
We continue to denote the special linear superalgebra
${\mathfrak{sl}}_{m|n}$ by $\gl$.
The aim of this section is to prove the following
\begin{theorem}\label{theo7.1}
Regarding the universal enveloping algebra $U(\gl)$ as a $\gl$-module
under the adjoint action, we have
$H^1(\gl,U(\gl))\ne0$, and $H^2(\gl,U(\gl))=0$.
\end{theorem}

The fact that $H^1(\gl,U(\gl))\ne0$ follows from Corollary \ref{coro3.2}.
Also, in the case $n=1$, it has been proved in \cite{SZ99} that $H^2(\gl,U(\gl))=0$.

To prove the remaining part of the theorem, we need some preparations.
Recall that $U(\gl)$, regarded as a $\gl$-module, is canonically isomorphic
to the super-symmetric algebra $S(\gl)$ (see \cite{Sc2}). Also, $S(\gl)$
is a direct summand in the $\gl$-module $S({\mathfrak{gl}}(m|n))$. If we
let $\C^{m|n}$ denote the natural ${\mathfrak{gl}}(m|n)$-module, and ${\bar\C}^{m|n}$ its dual,
then ${\mathfrak{gl}}(m|n)\cong \C^{m|n}\otimes {\bar\C}^{m|n}$. Here
${\mathfrak{gl}}(m|n)$ acts on $\C^{m|n}\otimes {\bar\C}^{m|n}$ as the
diagonal subalgebra of ${\mathfrak{gl}}(m|n)\times {\mathfrak{gl}}(m|n)$,
where the latter superalgebra acts on $\C^{m|n}\otimes {\bar\C}^{m|n}$ in the obvious way.
The ${\mathfrak{gl}}(m|n)\times
{\mathfrak{gl}}(m|n)$ action extends uniquely to $S(\C^{m|n}\otimes {\bar\C}^{m|n})$. Now
$S(\C^{m|n}\otimes {\bar\C}^{m|n})$ is isomorphic to the subalgebra
\cite[Definition 3.3]{SZ02} of regular
functions on the general linear supergroup.
A Peter-Weyl type theorem \cite[Proposition 3.1]{SZ02} states that as a ${\mathfrak{gl}}(m|n)\times
{\mathfrak{gl}}(m|n)$-module,
\equa{7.1}
{
S(\C^{m|n}\otimes {\bar\C}^{m|n})=\bigoplus_{\l} L_\l\otimes L^*_\l,
}
where $\l$ runs over all weights
$\l=(\l_1,...,\l_m\,|\,\l_{\cp 1},...,\l_{\cp n})$
such that
\equa{7.2}
{
\l_1\ge\l_2\ge...\ge\l_m\ge0,\ \
\l_{\cp1}\ge\l_{\cp2}\ge...\ge\l_{\cp n}\ge0,\ \
\l_m\ge\#\{\nu\in{\bf I}_2\,|\,\l_{\nu}\ne0\}
}
(thus in this section, $\l$ will not satisfy the condition in (\ref{weight1})).

Consider the action of the diagonal ${\mathfrak{gl}}(m|n)$ subalgebra
of ${\mathfrak{gl}}(m|n)\times {\mathfrak{gl}}(m|n)$ on $S(\C^{m|n}\otimes {\bar\C}^{m|n})$,
and restrict it to an $\gl$-action. Then $S(\gl)$ regarded as a
$\gl$-module under the adjoint action can be embedded as a direct
summand in $S(\C^{m|n}\otimes {\bar\C}^{m|n})$.

The proof of the theorem thus shifts its focus to the $\gl$-submodules
$L_\l\otimes L^*_\l$ of $S(\C^{m|n}\otimes {\bar\C}^{m|n})$.
We divide the proof into a series of technical lemmas.

\begin{lemma}\label{lemm7.1}
If $V=V_\mu\otimes W$, where $\mu$ is an integral dominant  weight and $W$ is any
finite-dimensional $\gl$-module, then there exists a filtration
$($called a Kac flag$)$
\equa{7.3}
{
0=V_0\subset V_1\subset\cdots\subset V_k=V,
}
such that $V_i/V_{i-1}$ is a Kac module for $i=1,...,k$ $($we call such a Kac module
a factor Kac module of $V)$.
\end{lemma}
\begin{proof}
Let $\gl_{\ge 0}=\gl_0+\gl_{+1}$, and denote by $L_\mu^{(0)}$ the
irreducible $\gl_{\ge 0}$-module with highest weight $\mu$.
Then $V=\bigl(U(\gl)\otimes_{U(\gl_{\gl\ge 0})} L_\mu^{(0)}\bigr)\otimes W
=U(\gl)\otimes_{U(\gl_{\gl\ge 0})}(L_\mu^{(0)}\otimes
W)$, where $L_\mu^{(0)}\otimes W$ is regarded as a $\gl_{\ge 0}$-module
with the obvious action. Since the induction functor $U(\gl)\otimes_{U(\gl_{\gl\ge
0})}\!{\sc\!}-{\sc\,}$ is exact, by applying it to any composition series of the finite dimensional
$\gl_{\ge 0}$-module $L_\mu^{(0)}\otimes W$, we produce a Kac flag for $V$.
\end{proof}
\begin{lemma}\label{lemm7.2}
$H^2(\gl,L_\l\otimes L^*_\l)=0$ if $\l$ satisfying $(\ref{7.2})$ is typical
$($i.e., $(\l+\rho,\a)\ne0$ for all $\a\in\D_1^+).$
\end{lemma}
\begin{proof}
If $\l$ is typical then $L^*_\l=V_{\l^*}$ is a Kac module (in this case
$\l^*=2\rho_1-\l^R$), thus we can use Lemma
\ref{lemm7.1}. Note that the highest weight vector of each factor Kac module
of $L_\l\otimes L^*_\l$ has level $\le\level(\l)+\level(\l^*)=\level(2\rho_1)$.
Thus $H^2(\gl,V_i/V_{i-1})=0$ by Theorem \ref{theo4.1}.
Then from the short exact sequence $0\to V_{i-1}\to V_i\to V_i/V_{i-1}\to0$, by
(\ref{long-exact}), we obtain $H^2(\gl,V_i)=0$ by induction on $i$.
\end{proof}
\begin{lemma}\label{lemm7.3}
$H^2(\gl,L_\l\otimes L^*_\l)=0$ if $\l$ satisfying $(\ref{7.2})$ is atypical.
\end{lemma}
\begin{proof}
Recall that
the $(i,\nu)$-entry $A(\l)_{i,\nu}$ of the atypicality matrix $A(\l)$
of $\l$ is  (cf.~(\ref{atypicality-matrix}))
\equa{7.a5}
{
A(\l)_{i,\nu}=(\l+\rho,\a_{i,\cp\nu})=
\l_i+\l_{\cp\nu}+m-i+1-\nu\mbox{ \ for \ }i=1,...,m,\ \nu=1,...,n,
}
(the smallest element is $A(\l)_{m,\cp n}=\l_m+\l_{\cp n}+1-n$), from this we see
that $\l$ is atypical only if $\l_m\le n$ (cf.~(\ref{7.2})).
One observes that the northeast chains $NE_\l$
of $\l$ (cf.~\cite{HKV}, see also Examples \ref{exam5.1} and \ref{exam5.2}) satisfy
\begin{eqnarray}
\label{7.a7}
&&
(m,n)\in NE_\l,
\nonumber\\
&&
(i,j)\in NE_\l\ \LRar\
(k,\ell)\in NE_\l \mbox{ \ if \ }k\ge i\mbox{ \and \ }\ell\ge j.
\end{eqnarray}
We have the short exact sequence
\equan{7.a9}
{
0\to L_\l\otimes L^*_\l\to V_{\L_\l}\otimes L^*_\l\to
 (V_{\L_\l}/L_\l)\otimes L^*_\l\to0.
}
Since $\C$ is a submodule of $L_\l\otimes L^*_\l$ (also of $V_{\L_\l}\otimes L^*_\l$),
we also have another short exact sequence
\equan{7.a10}
{
0\to (L_\l\otimes L^*_\l)/\C\to
 (V_{\L_\l}\otimes L^*_\l)/\C\rb{3pt}{\mbox{$\ ^{\ f}_{\dis-\!\!\!\to}\,$}}
(V_{\L_\l}/L_\l)\otimes L^*_\l\to0.
}
Thus by (\ref{long-exact}),
\equa{7.a11}
{
\begin{array}{cccccc}
&H^1(\gl,(V_{\L_\l}\otimes L^*_\l)/\C)&
\rb{3pt}{\mbox{$\ ^{\ f^1}_{\dis-\!\!\!\to}\,$}}&
H^1(\gl,(V_{\L_\l}/L_\l)\otimes L^*_\l)
\vs{4pt}\\ \to&
H^2(\gl,(L_\l\otimes L^*_\l)/\C)&\to&
H^2(\gl,(V_{\L_\l}\otimes L^*_\l)/\C).
\end{array}
}
Note from Lemma \ref{lemm6.3} that $H^2(\gl, (V_{\L_\l}\otimes L^*_\l)/\C)
=H^2(\gl, V_{\L_\l}\otimes L^*_\l)=0$, where the last equality follows from the
proof of Lemma \ref{lemm7.2} (note that in this case $\l^*=2\rho_1-(\L_\l)^R$
since $-\l^*$ is the lowest weight of $V_{\L_\l}$).
Thus
\equa{7.a12}
{
H^2(\gl,L_\l\otimes L^*_\l)=H^2(\gl,(L_\l\otimes L^*_\l)/\C)
\cong H^1(\gl,(V_{\L_\l}/L_\l)\otimes L^*_\l)/\Im{f^1}.
}
\vskip4pt
{\bf Claim 1.} \ ${\rm dim\,}H^1(\gl,(V_{\L_\l}/L_\l)\otimes L^*_\l)=1$.
\vskip4pt
To prove this claim, we first examine the short exact sequence
\equan{7.a13}
{
0\to (V_{\L_\l}/L_\l)\otimes L^*_\l\to (V_{\L_\l}/L_\l)\otimes V_{\L_{\l^*}}
\to (V_{\L_\l}/L_\l)\otimes (V_{\L_{\l^*}}/L_{\l^*})\to 0,
}
which leads to the exact sequence
\equa{7.a14}
{
H^0(\gl,(V_{\L_\l}/L_\l){\sc\!}\otimes{\sc\!} (V_{\L_{\l^*}}/L_{\l^*}))
{\sc\!}\to{\sc\!}
H^1(\gl,(V_{\L_\l}/L_\l){\sc\!}\otimes{\sc\!} L^*_\l)
{\sc\!}\to{\sc\!}
H^1(\gl,(V_{\L_\l}/L_\l){\sc\!}\otimes{\sc\!} V_{\L_{\l^*}}).
}
\vskip4pt
{\bf Subclaim a)} \ ${\rm dim\,}H^1(\gl,(V_{\L_\l}/L_\l)\otimes V_{\L_{\l^*}})\le1$.
\vskip4pt
To prove this subclaim, by Lemma \ref{lemm7.1} and Theorem \ref{theo3.1}, it
suffices to prove that
\begin{eqnarray}
\label{7.a15}&&
\mbox{$V_{2\rho_1}$ is not a factor Kac module of
$W=(V_{\L_\l}/L_\l)\otimes V_{\L_{\l^*}}$},
\\\label{7.a16}&&
\mbox{$V_{2\rho_1+\amax}$ is a factor Kac module of
$W$ with multiplicity $\le1$}.
\end{eqnarray}
Since $W=U(\gl_{-1})((V_{\L_\l}/L_\l)\otimes L^{(0)}_{\L_{\l^*}})$, suppose
$w_{2\rho_1}\in (V_{\L_\l}/L_\l)\otimes L^{(0)}_{\L_{\l^*}}$ is a primitive vector
of weight $2\rho_1$, then $w_{2\rho_1}$ is a combination of the form
\equa{7.a17}
{
u\otimes v,\ \,u\in V_{\L_\l}/L_\l,\,v\in L^{(0)}_{\L_{\l^*}},
}
and $u$ is in some
$\gl_0$-highest weight submodule $L_\eta^{(0)}$ of $V_{\L_\l}/L_\l$ (for some
weight $\eta$ such that $2\rho_1=\eta'+\L_{\l^*}$, where $\eta'$ is a weight
in $L^{(0)}_\eta$) and the level of $\eta$ is
\equa{7.a18}
{
\level(\eta)=\level(2\rho_1)-\level(\L_{\l^*}).
}
Note that with restriction to $\gl_0$, we have $(2\rho_1)|_{\gl_0}=0$, thus
the fact $2\rho_1=\eta'+\L_{\l^*}$ means that
\equa{7.a19}
{
\eta|_{\gl_0}=\big(\L_{\l^*}|_{\gl_0}\big)^o,
}
where we use $\l^o$ to denote the dual of a weight $\l$ with respect to $\gl_0$
(i.e., in fact $\l^o|_{\gl_0}=-\l^R|_{\gl_0}$).
But we have
\equa{7.a20}
{
\big(\L_{\l^*}|_{\gl_0}\big)^o=-(\L_{\l^*})^R|_{\gl_0}=-(2\rho_1-\l)|_{\gl_0}
=\l|_{\gl_0}.
}
Equations (\ref{7.a18})--(\ref{7.a20}) imply that $\eta=\l$.
Recall that (cf.~(\ref{Lambda_mu}))
\equa{7.a22}
{
\L_\l=\l+\mbox{$\sum\limits_{(i,\nu)\in NE_\l}$}\a_{i,\cp\nu}.
}
By (\ref{7.a7}), we see that when writing $\L_\l-\l$ as a sum of distinct
positive odd roots, there is only one way.
This shows that in $V_{\L_\l}$, there is
only one copy of $\gl_0$-highest weight $\l$ which occurs in $L_\l$.
This contradicts that $\eta=\l$ is a $\gl_0$-highest weight of $V_{\L_\l}/L_\l$.
This proves (\ref{7.a15}).

Next suppose $w_{2\rho_1+\amax}=\sum u\otimes v$ (as in
(\ref{7.a17})) is a primitive vector of weight $2\rho_1+\amax$. As
discussion above, now $\eta$ must be $\l+\a$ for some positive odd
roots $a=\a_{i,\cp\nu}$. Suppose $\a_{i,\cp\nu}\notin NE_\l$.
Then we must have some position of the $i$-th row  and some
position of the $\nu$-th column in $NE_\l$ (say if none of
positions of the $i$-th row is in $NE_\l$, then the $i$-th
coordinate of $\l+\a_{i,\cp\nu}$ is bigger than that of $\L_\l$ by
(\ref{7.a22}), a contradiction). Say $(i,\nu_1),(i_1,\nu)\in
NE_\l$, then by (\ref{7.a7}), we have \equan{7.a23} {
\L_\l-(\l+\a_{i,\cp\nu})=\mbox{$\sum\limits_{(k,\ell)\in
NE_\l\bs\{(i_1,\nu),
(i,\nu_1)\}}$}\a_{k,\cp\ell}+\a_{i_1,\cp\nu_1}, } (where
$\a_{i_1,\cp\nu_1}$ appears twice), which cannot be written as a
sum of distinct positive odd roots. Thus we must have
$\a_{i,\cp\nu}\in NE_\l$. Note from (\ref{7.a5}) and (\ref{7.a7})
that all atypical roots of $\l$ are $c$-related in sense of
\cite{SHK} (cf.~Definition \ref{nqc-relation}) or connected in sense of \cite{VZ}, $\l+\a_{i,\cp\nu}$
is dominant only if $(i,\nu)$ is located at the rightmost and the
topmost position of $NE_\l$, such position is unique by
(\ref{7.a7}). This proves (\ref{7.a16}) and Subclaim a).

\vskip4pt {\bf Subclaim b)} \ $H^0(\gl,(V_{\L_\l}/L_\l)\otimes
(V_{\L_{\l^*}}/L_{\l^*}))=0.$ \vskip4pt We have
\begin{eqnarray*}\label{7.a25}\nonumber
H^0(\gl,(V_{\L_\l}/L_\l)\otimes (V_{\L_{\l^*}}/L_{\l^*}))
&\!\!\!\!=\!\!\!\!&
((V_{\L_\l}/L_\l)\otimes (V_{\L_{\l^*}}/L_{\l^*}))^{\gl}
\\\nonumber
&\!\!\!\!=\!\!\!\!&
{\rm Hom}_{\gl}(\C,(V_{\L_\l}/L_\l)\otimes (V_{\L_{\l^*}}/L_{\l^*}))
\\&\!\!\!\!=\!\!\!\!&
{\rm Hom}_{\gl}(( V_{\L_{\l^*}}/L_{\l^*})^*,V_{\L_\l}/L_\l)
=
{\rm Hom}_{\gl}(M_\l,V_{\L_\l}/L_\l),
\end{eqnarray*}
where $M_\l$ is the maximal proper submodule of
$V_\l$ (note that $(V_{\L_{\l^*}})^*=V_\l$ since the lowest weight of
$V_{\L_{\l^*}}$ is $(\L_{\l^*}-2\rho_1)^R=-\l$).  Suppose $\mu$ is a primitive
weight of $M_\l$, then $\mu=\l-\sum_{\a\in S}\a$, where $S$ is some
subset of $\D_1^+$ such that at least one atypical root $\g$ of $\l$ is in
$S$, 
and all roots
in $S$ are $\le\g$ (see \cite{HKV}; this fact can also be proved by
\cite[Conjecture 4.1]{VZ} and \cite[Main Theorem]{Br}).
By (\ref{7.a22}), we have
$\L_\l-\mu=\sum_{\a\in S\cup NE_\l}\a$. From
the property (\ref{7.a7}) of $NE_\l$ and the fact that
$\g\in S\cap NE_\l$, one can easily see that
$\sum_{\a\in S\cup NE_\l}\a\,$
cannot be written
as a distinct sum of positive odd roots. This means
that $\mu$ is not a primitive weight of $V_{\L_\l}$, that is,
$M_\l$ and $V_{\L_\l}$ do not have a common primitive weight, which implies
that ${\rm Hom}_{\gl}(M_\l,V_{\L_\l}/L_\l)=0$.
\vskip4pt
{\bf
Subclaim c)} \ ${\rm dim\,}H^1(\gl,(V_{\L_\l}/L_\l)\otimes
L^*_\l)\ge1$ (this together with (\ref{7.a14}) and Subclaims a)
and b) proves Claim 1). \vskip4pt Note that the top
$\gl_0$-highest weight submodule of $V_{\L_\l}/L_\l$ is
$\gl_0$-dual of $L_{\l^*}^{(0)}$. There exists a $\gl_0$-highest
weight vector \equan{7.a26} { v_{2\rho_1}=v_{\L_\l}\otimes
w+\cdots\mbox{ \ ($w$ is the lowest weight vector of
$L_{\l^*}^{(0)}$)} } of weight $2\rho_1$. This vector must be
strongly primitive since there is no vector with level higher than
that of $2\rho_1$. Let $W=U(\gl)w_{2\rho_1}$. Take \equa{7.a27}
{\mbox{$ v=\prod\limits_{i>1,\,\nu<n}f_{\a_{i,\cp\nu}}v_{2\rho_1}
=\big(\prod\limits_{i>1,\,\nu<n}f_{\a_{i,\cp\nu}}v_{\L_\l}\big)\otimes
w+\cdots. $}} Note that in $V_{\L_\l}$, the primitive vector
$v_\l$ of $L_\l$ has the form \equan{7.a28} {\mbox{$
v_\l=\prod\limits_{\a\in NE_\l}f_\a v_{\L_\l}+\cdots. $}} By
(\ref{7.a7}), we see that
$NE_\l\not\subset\{(i,\nu)\,|\,i>1,\,\nu<n\}$, thus \equan{7.a29}
{\mbox{$ \prod\limits_{i>1,\,\nu<n}f_{\a_{i,\cp\nu}}v_{\L_\l}\ne0
\mbox{ \ (as a vector in }V_{\L_\l}/L_\l).
$}} Thus the first term of (\ref{7.a27}) is non-zero, clearly this
term cannot be cancelled by any other terms in (\ref{7.a27}),
i.e., $v\ne0$.

Note that $W$ is a quotient module of the Kac module
$V_{2\rho_1}$, and in $V_{2\rho_1}$ the primitive vector of weight
$\mu^{(n-1)}=(n,1,...,1\,|\,-1,...,-1,-m)$ is precisely defined by
(\ref{7.a27}) (cf.~\cite[Theorem 6.6]{SHK}). Thus $v$ is a
primitive vector of weight $\mu^{(n-1)}$. We claim that $U(\gl)v$
is an irreducible submodule $L_{\mu^{(n-1)}}$ of $W$: First $W$
does not contain a trivial submodule since $(V_{\L_\l}/L_\l)\otimes
L^*_\l$ does not contain a trivial submodule as the proof of
Subclaim b), and in $V_{2\rho_1}$, $\mu^{(n-1)}$ is the only
primitive weight links to primitive weight $0$ since $\mu^{(n-1)}$
is the only primitive weight of $V_{2\rho_1}$ with non-zero
$1$-cohomology. Thus $U(\gl)v=L_{\mu^{(n-1)}}$ is an irreducible
submodule of $W$ (and of $(V_{\L_\l}/L_\l)\otimes L^*_\l$), this
shows that $H^1(\gl,(V_{\L_\l}/L_\l)\otimes L^*_\l)\ne0$ by Lemma
\ref{lemm6.5}. This proves Subclaim c) and Claim 1.
\vskip4pt {\bf
Claim 2.} \ In (\ref{7.a11}), the map $f^1$ is surjective, thus by
(\ref{7.a12}), $H^2(\gl,L_\l\otimes L^*_\l)=0$. \vskip4pt Similar
to the proof of Subclaim c), $(V_{\L_\l}\otimes L^*_\l)/\C$ has an
irreducible submodule $V_1=L_{\mu^{(n-1)}}$ which maps under $f$
to the irreducible submodule $V_2=L_{\mu^{(n-1)}}$ of
$(V_{\L_\l}/L_\l)\otimes L^*_\l$, thus induces the map
$f^1:H^1(\gl,V_1)$ onto $H^1(\gl,V_2)$ by (\ref{long-exact}).
Hence $f^1$ is onto. This proves Lemma \ref{lemm7.2}.
\end{proof}

Finally we return to the proof of Theorem \ref{theo7.1}.
By (\ref{7.1}) and Lemmas \ref{lemm7.1} and \ref{lemm7.2}, we have
\\\cl{$\dis
H^2(\gl,U(\gl))=H^2(\gl,S(\gl))\subset
\bigoplus_{\l} H^2(\gl,L_\l\otimes L^*_\l)=0.
$}\\
This completes the proof.

%
\section{Cohomology groups of the Lie superalgebra $C(n)$} %
%
\label{last-section}
\def\fg{C(n)}
In this last section, we generalize the results of previous sections to the
other type I classical Lie superalgebra $C(n)={\mathfrak{osp}}_{2|2n-2}$, which is
a $\Z$-graded subalgebra of ${\mathfrak{sl}}_{2|2n-2}$
(see for example \cite{V}) such that
$C(n)=\fg_{-1}\oplus\fg_0\oplus\fg_{+1}$ with
\\[4pt]\hspace*{1ex}$
\begin{array}{llllll}
\fg_0\!\!\!\!&=\!\!\!\!&\biggl\{
\biggl(\ \,
\rb{3pt}{\mbox{$
\rb{3pt}{\mbox{${}^{^{\sc \a}{\,}_{\sc -\a}}$}}
\rb{-8pt}{\mbox{${\ \,}_{^{\sc \b\ \,\g}_{\sc \d\ -\b^T}}$}}
\put(-58,0){$\line(1,0){55}$}\put(-31,-20){$\line(0,1){40}$}
$}}
\ \biggr)\in{\mathfrak{sl}}_{2|2n-2}
\,\biggl|\ \a\!\in\!\C,\,\b,\g,\d\!\in\!{\mathfrak{gl}}_{n-1},\,
\g^T\!=\!\g,\,\d^T\!=\!\d\biggr\}
\\
\!\!\!\!&\cong\!\!\!\!&\C\oplus{\mathfrak{sp}}_{2n-2},
\end{array}$\\[4pt]\hspace*{1ex}$
\begin{array}{lllll}
\fg_{-1}\!\!\!\!&=\!\!\!\!&\biggl\{
\biggl(\ \,
\rb{6pt}{\mbox{$
\rb{-8pt}{\mbox{${}_{^{\sc -\eta^T\,0}_{\sc \,\xi^T\ \,0}}$}}
\rb{3pt}{\mbox{${\ \,}^{^{\sc 0\ 0}_{\sc \xi\ \eta}}$}}
\put(-46,0){$\line(1,0){55}$}\put(-18,-20){$\line(0,1){40}$}
$}}
\ \ \biggr)\in{\mathfrak{sl}}_{2|2n-2}\,\biggl|\
\xi,\eta\mbox{ are row vectors of dimension }n-1\biggr\},
\end{array}$\\[4pt]\hspace*{1ex}$
\begin{array}{lllll}
\fg_{+1}\!\!\!\!&=\!\!\!\!&\biggl\{
\biggl(\ \,
\rb{6pt}{\mbox{$
\rb{-8pt}{\mbox{${}_{^{\sc 0\,-\eta^T}_{\sc 0\ \,\xi^T}}$}}
\rb{3pt}{\mbox{${\ \,}^{^{\sc \xi\ \eta}_{\sc 0\ 0}}$}}
\put(-46,0){$\line(1,0){55}$}\put(-18,-20){$\line(0,1){40}$}
$}}
\ \ \biggr)\in{\mathfrak{sl}}_{2|2n-2}\,\biggl|\
\xi,\eta\mbox{ are row vectors of dimension }n-1\biggr\},
\end{array}$\\[4pt]
where the up-index ``$T$'' stands for the transpose of a matrix or vector.

Denote $\es=\es_1,\,\d_i=\es_{i+2}\,i=1,2,...,n-1.$ Then we have
$$
\D_0^+=\{\d_i-\d_j,\, \d_i+\d_k\,|\,1\!\le\! i,j,k\!\le\! n\!-\!1,\,i\!<\!j\},\ \
\D_1^+=\{\es\pm\d_j\,|\,1\le j\le n-1\}.
$$
Denote\\\cl{$
\begin{array}{lll}
&&
h_1=E_{11}-E_{22}+E_{33}-E_{n+2,n+2},
\\&&
h_i=E_{i+1,i+1}-E_{i+2,i+2}-E_{n+i,n+i}
+E_{n+i+1,n+i+1}\,(2\le i\le n-1),
\\&&
h_n=E_{n+1,n+1}-E_{2n,2n},
\end{array}$}\\
which forms a basis of the Cartan subalgebra $\fh$.
A weight $\l=\l_0\es+\sum_{i=1}^{n-1}\l_i\d_i\in\fh^*$ can be written as
$$
\l=(\l_0|\l_1,...,\l_{n-1})=[a_1;a_2,...,a_n],\mbox{ \ where }a_i=\l(h_i),
$$
and $\,a_i=\l_{i-1}+\l_i\,(1\le i\le n-1),\,a_n=\l_{n-1}$.
We have
\\[2pt]\cl{\mbox{$
\rho_0=\sum\limits_{i=1}^{n-1}(n-i)\d_i,\ \ \rho_1=(n-1)\es,
$}}\\[2pt]
and $\amin=\es-\d_1,\,\amax=\es+\d_1$ are respectively the minimal and maximal
positive odd roots.
\begin{theorem}\label{last-theo}
We have
\begin{eqnarray}
\label{theo8.1}&&
{\rm dim\,}H^1(\fg,V_\l)=\biggl\{\begin{array}{ll}1&\mbox{if \ }\l
=2\rho_1,2\rho_1+\amax,\\0&\mbox{otherwise},\end{array}
\\ 
\label{theo8.2}&&
{\rm dim\,}H^2(\fg,V_\l)=\biggl\{\begin{array}{ll}1&\mbox{if \ }\l
=2\rho_1+\amax,2\rho_1+2\amax,\\0&\mbox{otherwise},\end{array}
\\ 
\label{theo8.3}&&
{\rm dim\,}H^1(\fg,L_\l)=\biggl\{\begin{array}{ll}1&\mbox{if \ }\l
=-\amin,2\rho_1,\\0&\mbox{otherwise},\end{array}
\\
\label{theo8.4}&&
{\rm dim\,}H^2(\fg,L_\l)=\biggl\{\begin{array}{ll}1&\mbox{if \ }\l
=2\rho_1+\amax,-2\amin,\\0&\mbox{otherwise},\end{array}
\end{eqnarray}
and
\begin{equation}
\label{theo8.5}
H^1(\fg,U(\fg))\ne0\mbox{ \ \ and \ \ }H^2(\fg,\U(\fg))=0.
\end{equation}
\end{theorem}
\begin{proof}
Note that the Kac module $V_\l$ over $C(n)$ is at most singly atypical
(for some discussions of Kac module over $C(n)$, see for example \cite{V}).
The proofs of (\ref{theo8.1}) and (\ref{theo8.2}) are exactly similar to those
of Theorem~\ref{theo3.1} and Theorem~\ref{theo4.1} for the special case
${\mathfrak{sl}}_{m|1}$.

Suppose $H^1(\fg,L_\l)\ne0$, then we have $\l\llar0$
(recall~Definition \ref{defi6.1}).
If $\l$ has level lower than $0$, then $\l\llar0$ is a highest weight module
of highest weight $0$ and thus is
the Kac module $V_0$. In this case $\l$ is the only
other primitive weight of $V_0$, i.e., $\l=-\amin$.
If $\l$ has level higher than $0$, then we have a module $\l\rrar0$, which is now
the Kac module $V_\l$ such that $0$ is a primitive weight, thus $\l=2\rho_1$.
Hence,
\\[2pt]\cl{$
H^1(\fg,L_\l)\ne0\ \ \LRar\ \ \l=-\amin,\,2\rho_1,
$}\\[2pt]
from this one immediate obtains (\ref{theo8.3}) using (\ref{long-exact}).

Now suppose $H^2(\fg,L_\l)\ne0$, then we have $\l\llar\mu$ for some $\mu$ with
$H^1(\fg,L_\mu)$. Since either $\l\llar\mu$ or $\l\rrar\mu$ must be a Kac module,
from this we obtain $\l=2\rho_1+\amax,\,-2\amin$.
From this one immediate obtains (\ref{theo8.4}) using (\ref{long-exact}).

The proof of the first equation of (\ref{theo8.5}) is again similar as for the case of the
special linear superalgebra. As for the second equation, note that regarded as a $\fg$-module
under the adjoint action, $U(\fg)$ does not have any weight with level higher than that of $2\rho_1$,
thus $L_{2\rho_1+\amax}$ is not a composition factor of $U(\fg)$. Similarly, $L_{-2\amin}$ cannot appear
as a composition factor of $U(\fg)$ either, since the lowest $\fg_{\bar0}$-highest weight of
$L_{-2\amin}$ is $-2\rho_1-\amin$, which is not a weight of $U(\fg)$.
Thus no composition factors of $U(\fg)$ has non-zero 2-cohomology,
which implies $H^2(\fg,U(\fg))=0$. 
\ \hspace*{1pt}\end{proof}

\appendix
\section{Atypicality}\label{Appendix}
We briefly recall the definitions \cite{HKV, SHK} of atypical roots
$\g_1,..,\g_r$, atypicality matrix $A(\mu)$, $nqc$-relationship
of atypical roots, and northeast chains $NE_\mu$, of any weight $\mu$,
and also illustrate the concepts by some examples. We consider the
special linear superalgebra $\gl={\mathfrak{sl}}_{m|n}$ only.
Results presented here are used in Section \ref{1-irr-modules}.
\begin{definition}\label{atypicality-matrix}
We define the atypicality matrix of $\mu$
to be the $m\times n$ matrix  $A(\mu)=(A(\mu)_{i,\eta})_{m\times n}$, where
\equa{Aty-matrix}{A(\mu)_{i,\eta}=(\mu+\rho,\a_{i,\cp \eta})=
\mu_i+\mu_{\cp \eta}+m-i+1-\eta,\ 1\le i\le m,\ 1\le \eta\le n.}
An odd root $\a_{i,\cp\eta}$ is an atypical root of $\mu$ if
$A(\mu)_{i,\eta}=0$.
We label the atypical roots of an $r$-fold atypical
dominant integral weight $\mu$ by
$\g_1<...<\g_r$.
\end{definition}
\begin{definition}\label{nqc-relation}
{\rm(1)} For $1\le s,\,t\le r$, let $x_{st}$
be the entry in $A(\mu)$ at the intersection of the column containing
the $\g_s$ zero with the row containing the $\g_t$ zero, then
$x_{st}\in\Z_+\bs\{0\}$ for $s < t$ and $x_{ts}=- x_{st}$.
Denote $h_{st}$ the hook length between the zeros corresponding to
$\g_s,\g_t$, i.e., the number of steps
to go from the $\g_s$ zero via $x_{st}$ to the $\g_t$ zero with the zeros
themselves included in the count.
\[
\begin{array}{rlll}
\mbox{\rm (i)}&
\mbox{
$\g_s,\g_t$ are normally related {\rm($n$)}
}&\Lra& x_{st}>h_{st}-1;
\\
\mbox{\rm (ii)}&
\mbox{
$\g_s,\g_t$ are quasi-critically related {\rm($q$)}
}&\Lra&x_{st}=h_{st}-1;
\\
\mbox{\rm (iii)}&
\mbox{
$\g_s,\g_t$ are critically related {\rm($c$)}
}&\Lra&x_{st}<h_{st}-1.
\end{array}
\]
It is straightforward to show that the $q$-relation is transitive,
i.e., if $\g_s,\g_t$ are $q$-related and $\g_t,
\g_u$ are $q$-related, then $\g_s,\g_u$ are $q$-related.

{\rm(2)} The $nqc$-type {\rm(}atypicality type{\rm)}
of an $r$-fold atypical $\mu$ is a triangular array
\\\cl{
$
nqc(\mu)=\ \begin{array}{l}
d_{1r}\cdots d_{sr}\cdots d_{tr}\cdots 0\\[-4pt]
\stackrel{\vdots}{d_{1t}}\cdots \stackrel{\vdots}{d_{st}}\cdots
\stackrel{\vdots}0{\ }^{^{\sc\cdot\,\,^{^{\sc\cdot\,\,^{^{\sc\cdot}}}}}}\\[-4pt]
\stackrel{\vdots}{d_{1s}}\cdots \stackrel{\vdots}
0{\ }^{^{\sc\cdot\,\,^{^{\sc\cdot\,\,^{^{\sc\cdot}}}}}}\\[-4pt]
\stackrel{\vdots}0{\ }^{^{\sc\cdot\,\,^{^{\sc\cdot\,\,^{^{\sc\cdot}}}}}}
\end{array}
$}\\
where the zeros correspond to $\{\g_1,\cdots,\g_r\}$ and
$d_{st}=n,q,c\ \Lra\ \g_s,\g_t$ are $n$-,
$q$-, $c$-related respectively.

\end{definition}
\begin{definition}\label{ne-chain}
Denote by
$D=\{(i, \eta)\,|\,1\le i\le m,\,1\le \eta\le n\}$
the set of positions of $A(\mu)$.
For $1\le s\le r$, let $(b_s,c_s)$ be the position corresponding
to $\g_s$. The east chain $E_\mu(s)$ emanating from
$(b_s,c_s)$ is a sequence of positions in $D$,
whose entry are $<0$ except the entry of $(b_s,c_s)$,
starting at $(b_s,c_s)$ and extending in an easternly or north-easternly
direction until it reaches the last column
or the position defferent from $(b_s,c_s)$ whose entry is $\ge0$,
or
it cannot extend further
without leaving $A(\mu)$ by passing above its first row. For all $\eta$
with $c_s\le\eta\le n$,
$E_\mu(s)$ has at most one element in the $\eta$-th
column; if $c_s\le \eta< n$, the row of the position in the $(\eta+1)$-th
column is $a_{\cp\eta}$ rows above the row of the position in
the $\eta$-th column, where $a_{\cp\eta}$ is a Dynkin label of $\mu$;
if this is not possible, i.e.,
if this row would be the $M$-th row where $M<1$, or
the entry of this position is $\ge0$,
then $E_\mu(s)$ ends in the $\eta$-th column,
i.e., has no position to the right of the $\eta$-th
column.
Thus $E_\mu(s)$ is the maximal set
of positions $(b_s,c_s)=(i_0,\eta_0),(i_1,
\eta_1),...,(i_k,\eta_k)$ satisfying
\begin{eqnarray}\label{E-chain}
\eta_{p+1}=\eta_p+1\le n,\
i_{p+1}=i_p-a_{\cp \eta_p}\ge1,\
A(\mu)_{i_{p+1},\eta_{p+1}}<0,\ \ p=0,1,...,k-1.
\end{eqnarray}
Similarly, the nouth chain $N_\mu( s)$ emanating from $(b_s,c_s)$
is the maximal set
of positions $(b_s,c_s)=(i_0,\eta_0),
(i_1,\eta_1),...,(i_k,\eta_k)$ satisfying
\begin{eqnarray}\label{N-chain}
i_{p+1}=i_p-1\ge 1,\
\eta_{p+1}=\eta_p+a_{i_{p+1}}\le n,\
A(\mu)_{i_{p+1},\eta_{p+1}}>0,\ \ p=0,1,...,k-1.
\end{eqnarray}
Set $NE_\mu(s)=N_\mu(s)\cup E_\mu(s)$,
and $NE_\mu=\cup_{s=1}^r NE_\mu(s)$,
called the set of northeast chains of $\mu$.
\end{definition}

We shall
illustrate the northeast chains
by giving two examples below. First
note from
\cite[Theorem 3.2]{S} or \cite[Conjecture 4.1]{VZ}
that
\equa{Lambda_mu} {\mbox{$ \L_\mu=\mu+\sum\limits_{\a\in
NE_\mu}\a. $}}

\begin{example}\label{exam5.1}
$\mu=(3, 2,  2,  1,  0,  0\, |\, 0,  0,  -1,  -3,  -4)=[1,0,1,1,0;0;0,1,2,1]$.
{\small\equan{ex1}
{
A(\mu)=
\put(40,-45){$0\ \,1\ \,2\ \,1$}
\begin{array}{l}1\\0\\1\\1\\0\end{array}
\left(\begin{array}{ccccc}
   8\!&\! 7\!&\! 5\!&\! 2\!&\! 0\\
             6\!&\! 5\!&\! 3\!&\! 0\!&\!\bar2 \\
             5\!&\! 4\!&\! 2\!&\!\bar1\!&\!\bar3 \\
             3\!&\! 2\!&\! 0\!&\!\bar3\!&\!\bar5 \\
             1\!&\! 0\!&\!\bar2\!&\!\bar5\!&\!\bar7 \\
             0\!&\!\bar1\!&\!\bar3\!&\!\bar6\!&\!\bar8
\end{array}\right), \ \
NE_\mu=
\put(40,-50){$0\ \,1\ \,2\ \,1$}
\put(59,-14){$\ssc\vector(1,2){9}$}
\begin{array}{l}1\\[2pt]0\\[2pt]1\\[2pt]1\\[2pt]0\end{array}
\left(\begin{array}{ccccc}*\!&\! *\!&\! *\!&\! 1\!&\! 5  \\[-3pt]
             *\!&\! *\!&\! 1\put(0,10){$\ssc\nearrow$}
             \!&\! 4\!&\! 1  \\[-3pt]
             *\!&\! *\!&\! 1\put(-5,10){$\ssc\uparrow$}
             \!&\! 1\put(0,10){$\ssc\nearrow$}
             \!&\! *  \\[-3pt]
             *\!&\! 1\put(0,10){$\ssc\nearrow$}
             \!&\! 3\!&\! *\!&\! *  \\[-3pt]
         1\put(0,10){$\ssc\nearrow$}\!&\! 2\!&\!
         1\!&\! *\!&\! *  \\[-3pt]
             1\put(0,0){$\ssc\rar$}\put(-5,10){$\ssc\uparrow$}\!&\!
             1\put(0,10){$\ssc\nearrow$}\!&\! *\!&\! *\!&\! *
\end{array}\right),
\ \
nqc(\mu)=
\begin{array}{ccccc}
\!c\!&\!q\!&\!q\!&\!q\!&\!0\\\!c\!&\!q\!&\!q\!&\!0\\\!c\!&\!q\!&\!0\\\!c\!&\!0\\\!0
\end{array},
}}where,
we have placed the Dynkin labels $a_1,...,a_{m-1}$ to the left of
the first column, and in between the rows, of $A(\mu)$,
likewise, $a_{\cp1} ,...,a_{\cp n-1}$ are placed below the
bottom row, and in between the columns, of $A(\mu)$;
and where the
 numbered positions in $NE_\mu$ are located in the
northeast chains {\rm(}here and below we identify a matrix position
$(a,b)$ with the corresponding positive odd root $\a_{a,\cp b}${\rm)}.
The atypical roots are $\g_1=\a_{6,\cp 1}=\amin,\g_2=\a_{5,\cp 2},
\g_3=\a_{4,\cp 3},\g_4=\a_{2,\cp 4},\g_5=\a_{1,\cp 5}=\amax$, and
\equan{ex1.1} {\mbox{$ \L_\mu=\mu+\sum\limits_{\a\in
NE_\mu}\a=(5,5,4,3,3,2\,|\,-2,-3,-5,-6,-6). $}}
\end{example}

\begin{example}\label{exam5.2}
$\mu=( 5,  2,  2,  1,  1,  1\,|\,  -1,  -1,  -1,  -3,  -6 )=
[3,0,1,0,0;0;0,0,2,3]$.
{\small
\equan{ex2}
{
A(\mu)=
\put(40,-45){$0\ \,0\ \,2\ \,3$}
\begin{array}{l}3\\0\\1\\0\\0\end{array}
\left(\begin{array}{ccccc}   9\!&\! 8\!&\! 7\!&\! 4\!&\! 0 \\
             5\!&\! 4\!&\! 3\!&\! 0\!&\!\bar4 \\
             4\!&\! 3\!&\! 2\!&\!\bar1\!&\!\bar5\\
             2\!&\! 1\!&\! 0\!&\!\bar3\!&\!\bar7 \\
             1\!&\! 0\!&\!\bar1\!&\!\bar4\!&\!\bar8\\
             0\!&\!\bar1\!&\!\bar2\!&\!\bar5\!&\!\bar9
\end{array}\right), \ \
NE_\mu=
\put(40,-50){$0\ \,0\ \,2\ \,3$}
\put(59,-10){$\ssc\vector(1,2){8}$}
\put(59,-26){$\ssc\vector(1,2){9}$}
\begin{array}{l}3\\[2pt]0\\[2pt]1\\[2pt]0\\[2pt]0\end{array}
\left(\begin{array}{ccccc}*\!&\! *\!&\! *\!&\! *\!&\! 5  \\
         *\!&\! 1\!&\! 2\!&\! 4\!&\! *  \\[-3pt]
             *\!&\! 1\put(-5,10){$\ssc\uparrow$}
             \!&\! 2\put(-5,10){$\ssc\uparrow$}
             \!&\! 2\!&\! *\\[-3pt]
             1\put(0,10){$\ssc\nearrow$}
             \!&\! 2\put(0,10){$\ssc\nearrow$}
             \!&\! 3\!&\! 1\!&\! *  \\[-3pt]
             1\put(-5,10){$\ssc\uparrow$}
             \!&\! {\ssc\!}2\put(0,0){$\ssc\rar$}\put(-5,10){$\ssc\uparrow$}
             \!&\! 2\!&\! *\!&\! *  \\[-3pt]
             1\put(0,0){$\ssc\rar$}\put(-5,10){$\ssc\uparrow$}
             \!&\! 1\put(0,0){$\ssc\rar$}
             \!&\! 1\!&\! *\!&\! *
\end{array}\right), \ \
nqc(\mu)=\begin{array}{ccccc}
\!q\!&\!n\!&\!n\!&\!n\!&\!0\!\\\!c\!&\!c\!&\!q\!&\!0\!\\\!c\!&\!c\!&\!0\!
\\\!c\!&\!0\!\\\!0\!
\end{array}.
}}
The atypical roots are $\g_1=\amin,\g_2=\a_{5,\cp 2},\g_3=\a_{4,\cp 3},
\g_4=\a_{2,\cp 4},\g_5=\a_{1,\cp 5}=\amax$,
and
\equan{ex2.2}
{\mbox{$
\L_\mu=\mu+\sum\limits_{\a\in NE_\mu}\a=(6,5,5,5,4,4\,|\,-4,-6,-6,-6,-7).
$}}
\end{example}

For a weight $\mu$ of the form (\ref{mu_form}), its atypical roots are
\equan{atypical_roots}
{
\g_1=\amin=\a_{m,\cp 1},...,\g_{n-k}=\a_{m-n+k+1,\cp n- k},\ \
\g_{n-k+1}=\a_{k,\cp n- k+ 1},...,\g_n=\a_{1,\cp n}=\amax.
}
We denote
\equa{P_mu}
{
P_\mu
=\{\a_{i,\cp\nu}\in\D_1^+\,|\,
\nu\ge \mu_i-\mu_m\mbox{ and }i\le \mu_{\cp 1}-\mu_{\cp\nu}\}.
}
Note that
\begin{eqnarray}\label{NumP_mu}\nonumber
\#P_\mu
\!\!\!\!&=&\!\!\!\!mn-\big(\sum\limits_{i=1}^m(\mu_i-\mu_m)
+\sum\limits_{\nu=1}^{n}(\mu_{\cp 1}-\mu_{\cp\nu})
\big)
\\&=&\!\!\!\!
m{\sc\,}n-2\,\top+m{\sc\,}\mu_m-n{\sc\,}\mu_{\cp 1}
=m{\sc\,}n-2\,\top+(m+n)\mu_m.
\end{eqnarray}
{}From the definition of $NE_\mu$, one can easily obtain

\begin{lemma}
\label{lemm5.1}
 \ {\rm(1)} $P_\mu\subset NE_\mu$; \ {\rm(2)} $P_\mu=NE_\mu\ \Lra\
\g_1$ is $c$- or $q$-related to any other atypical roots.
\hfill$\Box$
\end{lemma}

\vskip10pt \noindent{\bf Acknowledgements}. Both authors
gratefully acknowledge financial support from the Australian
Research Council. Su is also supported by NSF grant 10171064 of
China, EYTP and TCTPFT grants of Ministry of Education of China.

\end{document}